**Some series and integrals involving the Riemann zeta function, binomial coefficients and the harmonic numbers**

**Volume VI**

Donal F. Connon

## Abstract

In this paper, the last of a series of seven, predominantly by means of elementary analysis, we establish a number of identities related to the Riemann zeta function, including the following:

$$B_n = \sum_{k=0}^{n} \frac{(-1)^k k!}{k+1} S(n,k)$$

$$\int_0^1 \frac{x^{n-1}}{1+x^n} \log\log\left(\frac{1}{x}\right) dx = -\frac{\log(2)\log(2n^2)}{2n}$$

$$\int_0^1 \frac{1}{1+x} \log^{q-1}\left(\frac{1}{x}\right)\left[\log\log\left(\frac{1}{x}\right)\right]^2 dx = \Gamma''(q)\varsigma_a(q) + 2\Gamma'(q)\varsigma_a'(q) + \Gamma(q)\varsigma_a''(q)$$

$$\gamma + \frac{1}{2n} - \sum_{k=1}^{2N+1} \frac{B_{2k}}{2kn^{2k}} < H_n - \log n < \gamma + \frac{1}{2n} - \sum_{k=1}^{2N} \frac{B_{2k}}{2kn^{2k}}$$

$$\lim_{n\to\infty}\left[\frac{1}{2}\left(H_n^{(1)}\right)^2 - \gamma\log n - \frac{1}{2}\log^2 n\right] = \frac{1}{2}\gamma^2$$

$$\lim_{n\to\infty}\left[\sum_{k=1}^{n} \frac{\left(H_k^{(1)}\right)^2}{k} + \sum_{k=1}^{n} \frac{H_k^{(2)}}{k} - H_n^{(1)}H_n^{(2)}\right] = \frac{4}{3}\varsigma(3)$$

$$\lim_{n\to\infty}\left[\frac{1}{6}\left(H_n^{(1)}\right)^3 + \frac{1}{2}H_n^{(1)}H_n^{(2)} - \frac{1}{6}\log^3 n - \frac{\gamma}{2}\log^2 n - \frac{1}{2}\left(\varsigma(2)+\gamma^2\right)\log n\right] = \frac{1}{2}\varsigma(2)\gamma + \frac{1}{6}\gamma^3$$

$$\log\Gamma(x) + \log x + \gamma x = \sum_{k=2}^{\infty} \frac{(-1)^k}{k}\varsigma(k)x^k$$

Whilst this paper is mainly expository, some of the formulae reported in it are believed to be new, and the paper may also be of interest specifically due to the fact that most of the various identities have been derived by elementary methods.

**CONTENTS OF VOLUMES I TO VI:** **Volume/page**

**SECTION:**

**APPENDICES (Volume VI):**

## APPENDIX A

## SOME PROPERTIES OF THE BERNOULLI NUMBERS AND THE BERNOULLI POLYNOMIALS

As stated in (1.9) of Volume I, the Bernoulli polynomials $B_n(x)$ are defined by the series

(A.1) $$\frac{te^{tx}}{e^t - 1} = \sum_{n=0}^{\infty} B_n(x) \frac{t^n}{n!} \qquad , (|t| < 2\pi)$$

The Bernoulli numbers $B_n$ are given by the generating function

(A.2) $$\frac{t}{e^t - 1} = \sum_{n=0}^{\infty} B_n \frac{t^n}{n!} \qquad , (|t| < 2\pi)$$

and it is readily seen from (A.1) that $B_n(0) = B_n$.

Using the Cauchy product [90, p.146] of two infinite series, equation (A.1) may be written as follows

$$\sum_{n=0}^{\infty} B_n(x) \frac{t^n}{n!} = \left(\frac{t}{e^t - 1}\right) e^{tx} = \left(\sum_{n=0}^{\infty} B_n \frac{t^n}{n!}\right)\left(\sum_{n=0}^{\infty} x^n \frac{t^n}{n!}\right)$$

$$= \sum_{n=0}^{\infty} t^n \left(\sum_{k=0}^{n} \frac{B_k x^{n-k}}{k!(n-k)!}\right)$$

and equating coefficients of $t^n$ we have

(A.3) $$B_n(x) = \sum_{k=0}^{n} \binom{n}{k} B_k x^{n-k}$$

This identity could also be obtained in a less rigorous manner by formally differentiating (A.1) with respect to $x$. It is therefore evident from (A.3) that $B_n(x)$ is a polynomial of degree $n$ and that the coefficient of $x^n$ is 1 (i.e., $B_n(x)$ is a monomial because $B_0 = 1$).

Using the identity

$$\frac{te^{(1+x)t}}{e^t - 1} - \frac{te^{xt}}{e^t - 1} = te^{xt}$$

and (A.1) we have

$$\sum_{n=0}^{\infty} B_n(1+x)\frac{t^n}{n!} - \sum_{n=0}^{\infty} B_n(x)\frac{t^n}{n!} = \sum_{n=0}^{\infty} x^n \frac{t^n}{n!}$$

Equating coefficients of $t^n$ we obtain

(A.4) $$B_n(1+x) - B_n(x) = nx^{n-1} \text{ for } n \geq 1.$$

Letting $x = 0$ in (A.4) we get

(A.5) $$B_n(1) = B_n(0) = B_n \text{ for } n \geq 2.$$

Now, substituting $x = 1$ in (A.3) we have for $n \geq 2$

(A.6) $$B_n = \sum_{k=0}^{n} \binom{n}{k} B_k$$

Using L'Hôpital's rule we have

$$\lim_{x \to o} \frac{x}{e^x - 1} = \lim_{x \to o} \frac{1}{e^x} = 1$$

and we therefore conclude that $B_0 = 1$.

Simple algebra shows that

$$\frac{x}{e^x - 1} = \frac{x}{2}\left(\frac{e^x + 1}{e^x - 1} - 1\right) = -\frac{x}{2} + \frac{x}{2}\left(\frac{e^x + 1}{e^x - 1}\right) = -\frac{x}{2} + \frac{x}{2}\left(\frac{e^{x/2} + e^{-x/2}}{e^{x/2} - e^{-x/2}}\right)$$

and therefore using (A.2) we have

(A.7) $$\frac{x}{2} + \sum_{n=0}^{\infty} B_n \frac{x^n}{n!} = \frac{x}{2}\coth(x/2)$$

The right hand side of (A.7) is an even function of $x$ and hence

$$B_{2n+1} = 0 \text{ for all } n \geq 1 \text{ and } B_1 = -\frac{1}{2}$$

It would be possible to determine $B_1$ using the Maclaurin series expansion of (A.2): however, using the expansion to determine other coefficients is difficult because of the complexity introduced by computing the higher derivatives. We can however use the recursion formula (A.6) to obtain succeeding values of $B_n$ as follows:

(A.8) $$B_0 = 1$$

$$B_1 = -\frac{1}{2}$$

$$B_2 = \frac{1}{6}$$

$$B_4 = -\frac{1}{30}$$

$$B_6 = \frac{1}{42}$$

From the method of recursive calculation shown by (A.6), it can immediately be seen that all of the $B_n$ are rational numbers. Since $\varsigma(2n)$ is always positive, we can see from (1.7) that the numbers $B_{2n}$ alternate in sign, i.e.

(A.9) $\qquad (-1)^{n+1} B_{2n} > 0$

From (1.1) we have

(A.10) $\qquad \varsigma(n) = \dfrac{1}{1-2^{1-n}} \displaystyle\sum_{k=1}^{\infty} \frac{(-1)^{k-1}}{k^n} = \frac{\varsigma_a(n)}{1-2^{1-n}}$

Since $\varsigma_a(n)$ is an alternating series, we have the inequality

$$1 - \frac{1}{2^n} < \varsigma_a(n) < 1$$

and hence we have the bounds for $\varsigma(n)$

$$\frac{1-2^{-n}}{1-2^{1-n}} < \varsigma(n) < \frac{1}{1-2^{1-n}}$$

Therefore, it is apparent that $\varsigma(n) \to 1$ as $n \to \infty$.

In [118a] Sasvári gave an elementary proof of Binet's formula for the gamma function

$$\Gamma(x+1) = \left(\frac{x}{e}\right)^x \sqrt{2\pi\, x}.e^{\vartheta(x)}$$

where $\qquad \vartheta(x) = \displaystyle\int_0^{\infty} \left(\frac{1}{e^t - 1} - \frac{1}{t} + \frac{1}{2}\right) \frac{e^{-xt}}{t}\, dt$

Since $\lim_{x\to\infty} \vartheta(x) = 0$, we immediately obtain Stirling's asymptotic formula

(A.11) $$n! = \Gamma(n+1) \sim \left(\frac{n}{e}\right)^n \sqrt{2\pi n}$$

and we have

$$(-1)^{n+1} B_{2n} \sim 4\pi\sqrt{e}\left(\frac{n}{\pi e}\right)^{2n+\frac{1}{2}}$$

Hence the absolute value of $B_{2n}$ grows rapidly in size. In fact, $B_{48} \sim 1.20866 \times 10^{23}$ and this is of the same order of magnitude as Avogadro's number N ~ $6.022 \times 10^{23}$. Indeed $B_{102}$ is approximately equal to the number of baryons ($10^{80}$) in the observable universe [106, p.728].

Since the left hand side of (A.1), with $x = 1/2$, is an even function of $t$, we deduce

(A.12) $B_{2n+1}(1/2) = 0$ for $n \geq 0$.

Differentiating (A.3) it is seen that

(A.13) $B_n'(x) = nB_{n-1}(x)$

We also have

(A.14) $B_n(1-x) = (-1)^n B_n(x)$

$$B_n(1+x) - B_n(x) = nx^{n-1}$$

Since $B_{2n}(1-x) = B_{2n}(x)$, we have

(A.14a) $B_{2n}(1) = B_{2n}(0) = B_{2n}$.

Similarly, $B_{2n+1}(1-x) = -B_{2n+1}(x)$ implies that

(A.14b) $B_{2n+1}(1) = -B_{2n+1}(0) = -B_{2n+1} = 0$.

With knowledge of $B_n$ we can compute the Bernoulli polynomials $B_n(x)$ using (A.3). The first few are:

(A.15) $B_0(x) = 1$

$$B_1(x) = x - \frac{1}{2}$$

$$B_2(x) = x^2 - x + \frac{1}{6}$$

$$B_3(x) = x^3 - \frac{3}{2}x^2 + \frac{1}{2}x = \frac{1}{4}x(x-1)(2x-1)$$

$$B_4(x) = x^4 - 2x^3 + x^2 - \frac{1}{30}$$

Since $B_{2n+1} = 0$, equation (A.7) can be written as

(A.16) $$\frac{x}{2}\coth(x/2) = \sum_{n=0}^{\infty} B_{2n} \frac{x^{2n}}{(2n)!} \qquad , (|x| < 2\pi)$$

and, since $\cot x = i \coth ix$, $(i = \sqrt{-1})$, replacing $x$ by $2ix$ we have the identity which was employed in Volume V

(A.17) $$x\cot x = \sum_{n=0}^{\infty} (-1)^n \frac{2^{2n} B_{2n}}{(2n)!} x^{2n} \qquad , (|x| < \pi)$$

The following identity is easily derived

(A.18) $$\cot x - \tan x = \frac{\cos x}{\sin x} - \frac{\sin x}{\cos x} = \frac{\cos^2 x - \sin^2 x}{\sin x \cos x} = 2\cot 2x$$

and this enables us to write

(A.19) $$\tan x = \sum_{n=1}^{\infty} (-1)^{n+1} \frac{2^{2n}(2^{2n}-1)B_{2n}}{(2n)!} x^{2n-1} \qquad , (|x| < \pi/2)$$

Since

(A.20) $$\cot x + \tan x = \frac{2}{\sin 2x}$$

we can also easily show that

(A.21) $$\frac{x}{\sin x} = \sum_{n=0}^{\infty} (-1)^{n+1} \frac{(2^{2n}-2)B_{2n}}{(2n)!} x^{2n}$$

Note that identity (A.20) is incorrectly recorded in Knopp's book [90, p.208].

We also have [90, p.239]

(A.22) $$\frac{1}{\cos x} = \sum_{n=0}^{\infty} (-1)^n \frac{E_{2n}}{(2n)!} x^{2n}$$

where $E_{2n}$ are the Euler numbers defined in (A.26) below.

Since $\cot(x/2) \to \infty$ as $x \to 2\pi$, one would suspect that the radius of convergence of (A17) would be $2\pi$. This is shown to be the case in [90, p.237] and [66, p.586]. The radius of convergence is equal to $2\pi$ because the nearest singularities of $z/(e^z - 1)$ are $\pm 2\pi i$.

We also have an explicit formula for $B_n$ given by F. Lee Cook in [47]: see also Gould's paper [73a], Rademacher's book [110a, p.9] and the short paper published by Rzadkowski [116a] in 2004. In what follows, I have applied Lee Cook's approach to the Bernoulli polynomial, rather than to $B_n$ .

Using (A.1) and the Maclaurin expansion we have

$$B_n(x) = \left[ \frac{d^n}{dt^n} \left( \frac{te^{tx}}{e^t - 1} \right) \right]_{t=0}$$

Since $t = \log\left[1 - (1 - e^t)\right]$ and $\log(1-u) = -\sum_{m=1}^{\infty} \frac{u^m}{m}$ for $|u| < 1$, we have

$$t = -\sum_{m=1}^{\infty} \frac{(1-e^t)^m}{m}$$

Therefore, dividing by $(e^t - 1)$ we get

$$\frac{te^{tx}}{e^t - 1} = \sum_{m=1}^{\infty} \frac{e^{tx}(1-e^t)^{m-1}}{m} = \sum_{k=0}^{\infty} \frac{e^{tx}(1-e^t)^k}{k+1} \quad \text{for } |t| < \log 2$$

Hence we have

$$B_n(x) = \sum_{k=0}^{\infty} \frac{1}{k+1} \left[ \frac{d^n}{dt^n} \left(1 - e^t\right) e^{tx} \right]_{t=0}$$

Using the Leibniz rule for the derivative of a product we note that the $n$ th derivative vanishes at $t = 0$ for $k \geq n+1$. Therefore, using the binomial theorem we have

$$B_n(x) = \sum_{k=0}^{\infty} \frac{1}{k+1} \sum_{j=0}^{k} (-1)^j \binom{k}{j} \left[ \frac{d^n}{dt^n} (e^{(j+x)t}) \right]_{t=0}$$

Hence we have

$$\text{(A.23)} \qquad B_n(x) = \sum_{k=0}^{n} \frac{1}{k+1} \sum_{j=0}^{k} (-1)^j \binom{k}{j} (x+j)^n$$

A different proof, using the Hurwitz-Lerch zeta function, was recently given by Guillera and Sondow [75aa]. They also noted that

$$\sum_{j=0}^{k}(-1)^j\binom{k}{j}(x+j)^n = 0 \text{ for } k > n = 0,1,2,...$$

and we therefore have

(A.23i) $$B_n(x) = \sum_{k=0}^{\infty}\frac{1}{k+1}\sum_{j=0}^{k}(-1)^j\binom{k}{j}(x+j)^n$$

When $x=0$ we obtain

(A.23aa) $$B_n(0) = B_n = \sum_{k=0}^{n}\frac{1}{k+1}\sum_{j=0}^{k}(-1)^j\binom{k}{j}j^n = \sum_{k=0}^{\infty}\frac{1}{k+1}\sum_{j=0}^{k}(-1)^j\binom{k}{j}j^n$$

Note the structural similarity of (A.23aa) with the Stirling numbers of the second kind which we briefly referred to in (3.97) et seq. of Volume I. In particular (3.100) shows that

$$S(n,k) = \frac{1}{k!}\sum_{j=0}^{k}(-1)^{k-j}\binom{k}{j}j^n$$

We therefore conclude that

(A.23a) $$B_n = \sum_{k=0}^{n}\frac{(-1)^k k!}{k+1}S(n,k)$$

Surprisingly, this simple relationship does not appear in the book "Concrete Mathematics" [75]: instead, the authors report a much more complex identity involving both kinds of the Stirling numbers [75, p.289]. I subsequently discovered that the formula (A.23a) is reported in the Wolfram Mathworld website dealing with the Bernoulli numbers. Furthermore, this identity was also proved by Kaneko in 2000 in his paper "The Akiyama-Tanigawa algorithm for Bernoulli numbers" [82a], albeit his method was somewhat less direct. Interestingly, in that paper Kaneko reports that the Bernoulli numbers were independently discovered by Takakazu Seki (1642-1708) one year prior to Jakob Bernoulli. A proof was also given by Akiyama and Tanigawa in [6ai]. That paper also mentions that

(A.23b) $$a_k = \sum_{j=1}^{k}s(k,j)b_j \quad \Leftrightarrow \quad b_k = \sum_{j=1}^{k}S(k,j)a_j$$

and hence we obtain

(A.23c) $$\sum_{r=1}^{k}s(k,r)B_r = \frac{(-1)^k k!}{k+1}$$

In fact, the Bernoulli numbers were originally introduced by Johann Faulhaber (1580-1635) in his book Academia Algebrae published in 1631. To his credit, Jakob

Bernoulli (1654-1705) did in fact give priority to Faulhaber in his treatise on probability theory, Ars Conjectandi, which was published posthumously in 1713.

Euler invented the Euler numbers [62] to study the sums

(A.24) $$T_k(n) = \sum_{j=0}^{n-1} (-1)^j j^k$$

The Euler Polynomials $E_n(x)$ are defined in [126, p.63] by means of the generating function

(A.25) $$\frac{2e^{xt}}{e^t + 1} = \sum_{n=0}^{\infty} E_n(x) \frac{t^n}{n!} \qquad ,(|t| < \pi)$$

and the Euler numbers are defined by

(A.26) $$E_n = 2^n E_n(1/2)$$

The Dirichlet Beta function is defined by

(A.27) $$\beta(s) = \sum_{n=0}^{\infty} \frac{(-1)^n}{(2n+1)^s}$$

and for $s = 2$, $\beta(2)$ is known as Catalan's constant $G$

$$G = \sum_{n=0}^{\infty} \frac{(-1)^n}{(2n+1)^2} = 0.915965...$$

It is still not known if $G$ is irrational (or indeed transcendental).

We have the relationship for the Dirichlet beta function

(A.28) $$\beta(2k+1) = \frac{(-1)^k (\pi/2)^{2k+1} E_{2k}}{2(2k)!}$$

The Dirichlet beta function was employed in the recent paper by Dalai [51].

Other miscellaneous identities are set out below for ease of reference.

(A.29) $$\sec x = \sum_{n=0}^{\infty} \frac{|E_{2n}|}{(2n)!} x^{2n} \qquad |x| < \pi/2$$

and using (A.28) we have

(A.30) $$\pi \sec(\pi x) = \sum_{n=0}^{\infty} 4^{n+1} \beta(2n+1) x^{2n}$$

## APPENDIX B

## A WELL-KNOWN INTEGRAL

We recall the famous dictum [133, p.123] of the renowned Irish mathematician and physicist, Lord Kelvin (William Thomson (1824-1907)), that

"A mathematician is a person to whom $\int_0^\infty e^{-x^2}\,dx = \frac{\sqrt{\pi}}{2}$ is as obvious as 1+1 = 2".

A very elegant proof of the above identity, presented by H.F. Sandham [117] in 1946, is set out below.

Let

$$I = \int_0^\infty e^{-t^2}\,dt$$

and, using the substitution $t = xy$, we obtain

$$I = \int_0^\infty e^{-x^2y^2}\,y\,dx$$

Therefore we have

$$2e^{-y^2} I = \int_0^\infty e^{-y^2(1+x^2)}\,2y\,dx$$

Integrating again from 0 to $\infty$ with respect to $y$

$$2I\int_0^\infty e^{-y^2}\,dy = \int_0^\infty dy\int_0^\infty e^{-y^2(1+x^2)}\,2y\,dx$$

and, using Fubini's theorem to interchange the order of integration, we have

$$= \int_0^\infty dx\int_0^\infty e^{-y^2(1+x^2)}\,2y\,dy = \int_0^\infty \frac{dx}{1+x^2} = \tan^{-1} x\Big|_0^\infty = \frac{\pi}{2}$$

Therefore we get $2I^2 = \pi/2$, and hence we have the important identity

$$\int_0^\infty e^{-x^2}\,dx = \frac{\sqrt{\pi}}{2}$$

I subsequently discovered that a similar proof was given by J.-A. Serret in his 1900 tome "Cours de Calcul Différentiel et Intégral", Vol. 2, Calcul Intégral (p.132), a copy of which may be viewed on the internet at The Cornell Library Historical Mathematics Monographs website.

## APPENDIX C

## EULER'S REFLECTION FORMULA AND RELATED MATTERS

**Theorem:** Euler's reflection formula for the gamma function states that

(C.1) $$\Gamma(x)\Gamma(1-x) = \frac{\pi}{\sin \pi x}$$

**Proof:**

The following is based on a proof provided by Dedekind [52] in 1853 (as outlined in an exercise in the recent book "Special Functions" [8a, p.49] by Andrews et al.). Various intermediate steps of this proof give rise to interesting identities in their own right.

Let us put

(C.2) $$\phi(x) = \int_0^\infty \frac{t^{x-1}}{1+t}\, dt$$

where we require $0 < x < 1$ for convergence (throughout this proof, I have found it useful to keep in mind the fact that we will eventually show that $\phi(x) = \frac{\pi}{\sin \pi x}$).

Using the substitution $t = sy$ it is easily seen that

(C.3) $$s^{-x}\phi(x) = \int_0^\infty \frac{t^{x-1}}{1+st}\, dt$$

Now let $s = 1/p$ in (C.3) to obtain

(C.4) $$p^{x-1}\phi(x) = \int_0^\infty \frac{t^{x-1}}{p+t}\, dt$$

Or equivalently

(C.5) $$s^{x-1}\phi(x) = \int_0^\infty \frac{t^{x-1}}{s+t}\, dt$$

Dividing (C.5) by $(s+1)$ and integrating we obtain

(C.6) $$\phi(x)\int_0^\infty \frac{s^{x-1}}{s+1}\, ds = \int_0^\infty \frac{1}{s+1}\left(\int_0^\infty \frac{t^{x-1}}{s+t}\, dt\right) ds$$

Using (C.2) we have

(C.7) $$\phi^2(x) = \int_0^\infty \frac{1}{s+1}\left(\int_0^\infty \frac{t^{x-1}}{s+t}\,dt\right)ds$$

Changing the order of integration we obtain

(C.8) $$\phi^2(x) = \int_0^\infty t^{x-1}dt \int_0^\infty \frac{1}{(s+1)(s+t)}\,ds$$

It is readily seen that

(C.9) $$\int_0^\infty \frac{1}{(s+1)(s+t)}\,ds = \int_0^\infty \left[\frac{1}{(1-t)(s+t)} - \frac{1}{(1-t)(s+1)}\right]ds$$

(C.10) $$= \frac{1}{1-t}\log\left(\frac{s+t}{s+1}\right)\Bigg|_0^\infty = \frac{1}{t-1}\log t$$

From (C.8) we therefore have

(C.11) $$\phi^2(x) = \int_0^\infty \frac{t^{x-1}\log t}{t-1}\,dt = \frac{d}{dx}\int_0^\infty \frac{t^{x-1}}{t-1}\,dt$$

We have the elementary integral $\int t^x dy = t^x / \log t$ and hence integrating (C.11) with respect to $x$ we have

(C.12) $$\int_{1-y}^{y} \phi^2(x)\,dx = \int_{1-y}^{y} dx \int_0^\infty \frac{t^{x-1}\log t}{t-1}\,dt$$

(C.13) $$= \int_0^\infty \frac{t^{y-1} - t^{-y}}{t-1}\,dt$$

Now, subtracting (C.3) from (C.5) we obtain

(C.14) $$\frac{(s^{x-1} - s^{-x})}{s-1}\ \phi(x) = \int_0^\infty \frac{t^{x-1}(t-1)}{(st+1)(s+t)}\,dt$$

We now integrate (C.14) with respect to $s$ over $(0,\infty)$ and use (C.13) to obtain

(C.15) $$\phi(x)\int_0^\infty \frac{(s^{x-1} - s^{-x})}{s-1}\,ds = \phi(x)\int_{1-x}^{x} \phi^2(t)dt = \int_0^\infty ds \int_0^\infty \frac{t^{x-1}(t-1)}{(st+1)(s+t)}\,dt$$

(C.16) $$= \int_0^\infty (t-1)t^{x-1}dt \int_0^\infty \frac{1}{(st+1)(s+t)}\,dt$$

Using partial fractions it is easily shown that

(C.17) $$\int_0^\infty \frac{1}{(st+1)(s+t)}\,ds = \int_0^\infty \frac{1}{(1-t^2)(s+t)} - \frac{t}{(1-t^2)(st+1)}\,ds$$

(C.18) $$= \frac{1}{1-t^2}\log\left(\frac{s+t}{st+1}\right)\Bigg|_0^\infty = \frac{2}{t^2-1}\log t$$

Substituting (C.18) in (C.16) we have

(C.19) $$\phi(x)\int_{1-x}^{x} \phi^2(t)\,dt = 2\int_0^\infty \frac{t^{x-1}\log t}{1+t}\,dt$$

(C.20) $$= 2\phi'(x)$$, using the definition of $\phi(x)$ in (C.1).

Substituting $t = 1/u$ in (C.2) it is easily seen that

(C.21) $$\phi(x) = \int_0^\infty \frac{t^{x-1}}{1+t}\,dt = \int_0^\infty \frac{t^{-x}}{1+t}\,dt = \phi(1-x)$$

In addition, using the substitution $t = 1-u$ we have

(C.22) $$\int_{1-x}^{x} \phi^2(t)\,dt = \int_{1-x}^{1/2} \phi^2(t)dt + \int_{1/2}^{x} \phi^2(t)\,dt = 2\int_{1/2}^{x} \phi^2(t)\,dt$$

Therefore, using (20) we obtain

(C.23) $$\phi(x)\int_{1/2}^{x} \phi^2(t)\,dt = \phi'(x)$$

Differentiating (C.23) gives us the differential equation

(C.24) $$\phi\phi'' - (\phi')^2 = \phi^4$$

With $\phi = 1/u$ we obtain

(C.24a) $$uu'' = (u')^2 - 1$$

and the solution of (C.24) with the initial conditions $\phi(1/2) = \pi$ and $\phi'(1/2) = 0$ is

(C.25) $$\phi(x) = \frac{\pi}{\sin \pi x}$$

However, without actually knowing the solution, it is not immediately obvious to me how one would integrate (C.24) from first principles.

Now making the substitution $t = s/(s+1)$ in the beta integral

$$B(x,y) = \int_0^1 t^{x-1}(1-t)^{y-1}dt$$

we obtain

(C.26) $$B(x,y) = \int_0^\infty \frac{s^{x-1}}{(1+s)^{x+y}}\, ds = \frac{\Gamma(x)\Gamma(y)}{\Gamma(x+y)}$$

Letting $y = 1 - x$ in (C.26) we have

(C.27) $$\phi(x) = \int_0^\infty \frac{s^{x-1}}{1+s} ds = \frac{\Gamma(x)\Gamma(1-x)}{\Gamma(1)}$$

We have therefore proved the theorem

(C.28) $$\Gamma(x)\Gamma(1-x) = \frac{\pi}{\sin \pi x}$$

Since $\Gamma(x+1) = x\Gamma(x)$ it is easily seen that this identity is in fact valid [128, p.149] for all $x \in \mathbf{C}$ except for $x = 0$ and all of the integers (positive and negative).

A number of interesting results automatically arise from various steps in the proof: for example, by differentiating (C.2) we obtain

(C.29) $$\frac{\pi^2 \cos \pi x}{\sin^2 \pi x} = -\int_0^\infty \frac{t^{x-1}\log t}{1+t} dt$$

By differentiating (C.5) with respect to $s$ we obtain

(C.30) $$\frac{\pi(x-1)s^{x-2}}{\sin \pi x} = -\int_0^\infty \frac{t^{x-1}}{(s+t)^2} dt$$

Alternatively, by differentiating (C.5) with respect to $x$ we obtain

(C.31) $$s^{x-1}\left(\frac{\pi \log s}{\sin \pi x} - \frac{\pi^2 \cos \pi x}{\sin^2 \pi x}\right) = \int_0^\infty \frac{t^{x-1}\log t}{s+t} dt$$

From (C.11) we have

(C.32) $$\frac{\pi^2}{\sin^2 \pi x} = \int_0^\infty \frac{t^{x-1}\log t}{t-1}\,dt$$

and differentiating (C.32) with respect to $x$ we obtain

(C.33) $$\frac{2\pi^3 \cos \pi x}{\sin^3 \pi x} = -\int_0^\infty \frac{t^{x-1}\log^2 t}{t-1}\,dt$$

We have

$$\phi(\alpha) = \int_0^\infty \frac{t^{\alpha-1}}{1+t}\,dt = \int_0^1 \frac{t^{\alpha-1}}{1+t}\,dt + \int_1^\infty \frac{t^{\alpha-1}}{1+t}\,dt$$

and using the substitution $t \to 1/t$ in the last integral we obtain

$$\int_0^\infty \frac{t^{\alpha-1}}{1+t}\,dt = \int_0^1 \frac{t^{\alpha-1}}{1+t}\,dt + \int_0^1 \frac{t^{-\alpha}}{1+t}\,dt$$

$$= \int_0^1 \frac{t^{\alpha-1} + t^{\beta-1}}{1+t}\,dt \quad \text{, where } \beta = 1-\alpha$$

If $0 < \alpha < 1$ then we also have $0 < \beta < 1$. Let us now consider the integral

$$\phi(\alpha) = \int_0^1 \frac{t^{\alpha-1} + t^{-\alpha}}{1+t}\,dt$$

Then, upon integration we have

$$\int_a^x \phi(\alpha)\,d\alpha = \int_a^x d\alpha \int_0^1 \frac{t^{\alpha-1} + t^{-\alpha}}{1+t}\,dt$$

and, interchanging the order of integration, we obtain

$$= \int_0^1 \frac{t^{x-1} - t^{-x}}{(1+t)\log t}\,dt - \int_0^1 \frac{t^{a-1} - t^{-a}}{(1+t)\log t}\,dt$$

This immediately reminded me of equation (4.4.112g) which was originally derived by Kummer and is valid for $\alpha > 0$, $\beta > 0$.

$$\int_0^1 \frac{t^{\alpha-1} - t^{\beta-1}}{(1+t)\log t} dt = \log \frac{\Gamma\left(\frac{1+\alpha}{2}\right)\Gamma\left(\frac{\beta}{2}\right)}{\Gamma\left(\frac{1+\beta}{2}\right)\Gamma\left(\frac{\alpha}{2}\right)}$$

As a result we have with $\beta = 1-\alpha$

(C.33a) $$\int_0^1 \frac{t^{x-1} - t^{-x}}{(1+t)\log t} dt = \log \frac{\Gamma\left(\frac{1+x}{2}\right)\Gamma\left(\frac{1-x}{2}\right)}{\Gamma\left(\frac{2-x}{2}\right)\Gamma\left(\frac{x}{2}\right)}$$

We therefore have

$$\int_a^x \phi(\alpha)d\alpha = \log\Gamma\left(\frac{1+x}{2}\right) + \log\Gamma\left(\frac{1-x}{2}\right) - \log\Gamma\left(\frac{2-x}{2}\right) - \log\Gamma\left(\frac{x}{2}\right)$$

$$-\log\Gamma\left(\frac{1+a}{2}\right) - \log\Gamma\left(\frac{1-a}{2}\right) + \log\Gamma\left(\frac{2-a}{2}\right) + \log\Gamma\left(\frac{a}{2}\right)$$

In (C.25) $\phi(\alpha)$ was defined as

$$\phi(\alpha) = \frac{\pi}{\sin \pi\alpha}$$

and we have the elementary integral

$$\int \frac{dx}{\sin \pi x} = \int \frac{dx}{2\sin(\pi x/2)\cos(\pi x/2)} = \int \frac{\sec^2(\pi x/2)}{2\tan(\pi x/2)} dx = \frac{1}{\pi}\log\tan(\pi x/2)$$

Hence we have

$$\int_a^x \phi(\alpha)d\alpha = \pi\int_a^x \frac{d\alpha}{\sin \pi\alpha} = \log\tan(\pi x/2) - \log\tan(\pi a/2)$$

and therefore we get

(C.34) $\log\tan(\pi x/2) - \log\tan(\pi a/2) =$

$$\log\Gamma\left(\frac{1+x}{2}\right) + \log\Gamma\left(\frac{1-x}{2}\right) - \log\Gamma\left(\frac{2-x}{2}\right) - \log\Gamma\left(\frac{x}{2}\right)$$

$$-\log\Gamma\left(\frac{1+a}{2}\right) - \log\Gamma\left(\frac{1-a}{2}\right) + \log\Gamma\left(\frac{2-a}{2}\right) + \log\Gamma\left(\frac{a}{2}\right)$$

where we require $0 < a < x < 1$ , or more generally $n < a < x < n+1$ with $n$ being a strictly positive integer to ensure that $\sin \pi\alpha \neq 0 \ \forall \alpha \in [a, x]$.

We can write

$$\log\Gamma\left(\frac{1+x}{2}\right)+\log\Gamma\left(\frac{1-x}{2}\right)-\log\Gamma\left(\frac{2-x}{2}\right)-\log\Gamma\left(\frac{x}{2}\right)$$

$$=\log\Gamma\left(\frac{1}{2}+\frac{x}{2}\right)+\log\Gamma\left(\frac{1}{2}-\frac{x}{2}\right)-\log\Gamma\left(1-\frac{x}{2}\right)-\log\Gamma\left(\frac{x}{2}\right)$$

$$=\log\Gamma(1-y)+\log\Gamma(y)-\log\Gamma(1-z)-\log\Gamma(z)$$

where $y=\frac{1}{2}-\frac{x}{2}$ and $z=\frac{x}{2}$.

$$=\log\left[\Gamma(1-y)\log\Gamma(y)\right]-\log\left[\Gamma(1-z)\log\Gamma(z)\right]$$

$$=\log\frac{\pi}{\sin\pi y}-\log\frac{\pi}{\sin\pi z}$$

$$=-\log\sin\pi y+\log\sin\pi z$$

Since $\sin\pi y=\sin\pi\left[\frac{1}{2}-\frac{x}{2}\right]=\cos\frac{\pi x}{2}$ we obtain

(C.35) $$\log\Gamma\left(\frac{1+x}{2}\right)+\log\Gamma\left(\frac{1-x}{2}\right)-\log\Gamma\left(\frac{2-x}{2}\right)-\log\Gamma\left(\frac{x}{2}\right)=\log\tan\left(\frac{\pi x}{2}\right)$$

and therefore

(C.36) $$\frac{\Gamma\left(\frac{1+x}{2}\right)\Gamma\left(\frac{1-x}{2}\right)}{\Gamma\left(\frac{2-x}{2}\right)\Gamma\left(\frac{x}{2}\right)}=\tan\left(\frac{\pi x}{2}\right)$$

The above result may be obtained more directly from (C.34) by letting $a=1/2$.

We then see from (C.33a) that

(C.36a) $$\int_0^1\frac{t^{x-1}-t^{-x}}{(1+t)\log t}dt=\log\tan\left(\frac{\pi x}{2}\right)$$

Letting $x=3/2$ we have

$$\frac{\Gamma\left(1+\frac{1}{4}\right)\Gamma\left(-\frac{1}{4}\right)}{\Gamma\left(\frac{1}{4}\right)\Gamma\left(1-\frac{1}{4}\right)} = \tan\left(\frac{3\pi}{4}\right)$$

Using Euler's functional equation we have $\Gamma(1-x) = -x\Gamma(-x)$ and hence

$$\Gamma\left(1-\frac{1}{4}\right) = -\frac{1}{4}\Gamma\left(-\frac{1}{4}\right)$$

$$\frac{\Gamma\left(1+\frac{1}{4}\right)}{\Gamma\left(\frac{1}{4}\right)} = -\frac{1}{4}\tan\left(\frac{3\pi}{4}\right) = \frac{1}{4}$$

This is a particular case of the more general identity

$$\Gamma\left(n+\frac{1}{4}\right) = \frac{1}{4^n}\Gamma\left(1+\frac{1}{4}\right)\prod_{k=1}^{n}(4k-3)$$

From (C.36), and using Euler's reflection formula in the denominator, we have (as recorded in [74, p.887])

(C.37) $$\Gamma\left(\frac{1+x}{2}\right)\Gamma\left(\frac{1-x}{2}\right) = \pi\sec\left(\frac{\pi x}{2}\right)$$

and letting $x = 1/2$ we have

(C.37a) $$\Gamma\left(\frac{3}{4}\right)\Gamma\left(\frac{1}{4}\right) = \pi\sec\left(\frac{\pi}{4}\right) = \pi\sqrt{2}$$

(this is in numerical agreement with the approximations contained in [1, p.255] for $\Gamma\left(\frac{3}{4}\right)$ and $\Gamma\left(\frac{1}{4}\right)$). According to Havil [78, p.55], no closed form expression for $\Gamma\left(\frac{1}{4}\right)$ is yet known.

We have Legendre's relation [66, p.493] for $x > 0$

(C.37b) $$\Gamma\left(\frac{x}{2}\right)\Gamma\left(\frac{1+x}{2}\right) = \frac{\sqrt{\pi}}{2^{x-1}}\Gamma(x)$$

and for $x = 1/2$ we get the same result as (C.37a). With $x \to 2x$ we get

$$\Gamma(x)\Gamma\left(x+\frac{1}{2}\right)=\frac{\sqrt{\pi}}{2^{2x-1}}\Gamma(2x)$$

Rudin [115, p.194] gives a short proof of Legendre's relation using the Bohr-Mollerup theorem [25, p.187]. Harald Bohr (1887-1951), who was the younger brother of the renowned Danish physicist Niels Bohr, was one of the most prominent Danish mathematicians in the first half of the 20th century (as well as being an international footballer who played in the 1908 Olympics for Denmark, beating France by the remarkable score of 17 to 1!). It is also possible to prove the Euler reflection formula by applying the Bohr-Mollerup theorem to the function

$$f(x)=\frac{\pi}{\Gamma(1-x)\sin\pi x}$$

From (C.13) we have

$$\int_{1-\alpha}^{\alpha}\phi^2(x)\,dx=\int_0^{\infty}\frac{t^{\alpha-1}-t^{-\alpha}}{t-1}\,dt$$

Since $\int\frac{1}{\sin^2\pi x}dx=-\frac{1}{\pi}\cot\pi x$ we have

$$\int_{1-\alpha}^{\alpha}\frac{\pi^2}{\sin^2\pi x}dx=\pi\left[\cot\pi(1-\alpha)-\cot\pi\alpha\right]=\int_0^{\infty}\frac{t^{\alpha-1}-t^{-\alpha}}{t-1}dt$$

We have

$$\int_0^{\infty}\frac{t^{\alpha-1}-t^{-\alpha}}{t-1}dt=\int_0^{1}\frac{t^{\alpha-1}-t^{-\alpha}}{t-1}dt+\int_1^{\infty}\frac{t^{\alpha-1}-t^{-\alpha}}{t-1}dt$$

and using the substitution $t\to 1/t$ in the second integral we obtain

$$\int_0^{\infty}\frac{t^{\alpha-1}-t^{-\alpha}}{t-1}dt=\int_0^{1}\frac{t^{\alpha-1}-t^{-\alpha}}{t-1}dt-\int_0^{1}\frac{t^{-\alpha}-t^{\alpha-1}}{t-1}dt$$

$$=2\int_0^{1}\frac{t^{\alpha-1}-t^{-\alpha}}{t-1}dt$$

and therefore we have

(C.38) $$\pi\left[\cot\pi(1-\alpha)-\cot\pi\alpha\right]=2\int_0^{1}\frac{t^{\alpha-1}-t^{-\alpha}}{t-1}dt$$

From the comments following (6.42) we know that $\cot \pi\alpha = -\cot \pi(1-\alpha)$ and hence we obtain (G&R [74, p.318])

(C.39) $$\pi \cot \pi\alpha = \int_0^1 \frac{t^{\alpha-1} - t^{-\alpha}}{1-t} dt$$

Since $\dfrac{t^{\alpha-1} - t^{-\alpha}}{1-t} = \sum_{n=0}^{\infty} \left( t^{n+\alpha-1} - t^{n-\alpha} \right)$ we have

(C.40) $$\int_0^1 \frac{t^{\alpha-1} - t^{-\alpha}}{1-t} dt = \sum_{n=0}^{\infty} \left( \frac{1}{n+\alpha} - \frac{1}{n+1-\alpha} \right)$$

Therefore we obtain the well-known series expansion

(C.41) $$\pi \cot \pi\alpha = \sum_{n=0}^{\infty} \left( \frac{1}{n+\alpha} - \frac{1}{n+1-\alpha} \right)$$

Letting $\pi\alpha = x$ and $\pi\alpha = \pi/2 - x$ we get the series

(C.42a) $$\cot x = \frac{1}{x} - \left( \frac{1}{\pi - x} - \frac{1}{\pi + x} \right) - \left( \frac{1}{2\pi - x} - \frac{1}{2\pi + x} \right) - \ldots = \frac{1}{x} - 2\sum_{n=1}^{\infty} \frac{x}{(n\pi)^2 - x^2}$$

(C.42b) $$\tan x = \left( \frac{1}{\frac{\pi}{2} - x} - \frac{1}{\frac{\pi}{2} + x} \right) + \left( \frac{1}{\frac{3\pi}{2} - x} - \frac{1}{\frac{3\pi}{2} + x} \right) + \ldots$$

Using the identity $\csc x = \dfrac{1}{\sin x} = \dfrac{1}{2}\tan\dfrac{x}{2} + \dfrac{1}{2}\cot\dfrac{x}{2}$ we have

(C.42c) $$\csc x = \frac{1}{x} + \left( \frac{1}{\pi - x} - \frac{1}{\pi + x} \right) - \left( \frac{1}{2\pi - x} - \frac{1}{2\pi + x} \right) + \ldots$$

Letting $x \to \pi/2 - x$ we get

(C.42d) $$\sec x = \left( \frac{1}{\frac{\pi}{2} - x} - \frac{1}{\frac{\pi}{2} + x} \right) - \left( \frac{1}{\frac{3\pi}{2} - x} - \frac{1}{\frac{3\pi}{2} + x} \right) + \ldots$$

We may write these equations as

(C.42e) $$\csc x = \frac{1}{x} + 2x\sum_{n=1}^{\infty}\frac{(-1)^{n+1}}{(n\pi)^2 - x^2}$$

(C.42f) $$\sec x = \pi\sum_{n=1}^{\infty}\frac{(-1)^{n+1}(2n-1)}{\left(\frac{(2n-1)\pi}{2}\right)^2 - x^2}$$

Differentiating (C.39) we obtain

$$-\frac{\pi^2}{\sin^2 \pi\alpha} = \int_0^1 \frac{\left(t^{\alpha-1} + t^{-\alpha}\right)\log t}{1-t}\,dt$$

As set out below, this is equivalent to the result which we obtained in (C.32)

$$-\frac{\pi^2}{\sin^2 \pi\alpha} = \int_0^\infty \frac{t^{\alpha-1}\log t}{1-t}\,dt$$

because

$$\int_0^\infty \frac{t^{\alpha-1}\log t}{1-t}\,dt = \int_0^1 \frac{t^{\alpha-1}\log t}{1-t}\,dt + \int_1^\infty \frac{t^{\alpha-1}\log t}{1-t}\,dt$$

$$= \int_0^1 \frac{t^{\alpha-1}\log t}{1-t}\,dt + \int_0^1 \frac{t^{-\alpha}\log t}{1-t}\,dt$$

(where we have used the substitution $t \to 1/t$ in the second integral).

Integrating (C.39) we obtain

$$\int_0^u \pi\cot\pi\alpha\, d\alpha = \int_0^1 \frac{t^{u-1} + t^{-u} - t^{-1} + 1}{(1-t)\log t}\,dt$$

In (C.21) we showed that $\phi'(x) = -\phi'(1-x)$ and hence $\phi'(1/2) = 0$.

We have therefore proved that

$$\int_0^\infty \frac{\log t}{\sqrt{t}(1+t)}\,dt = 0$$

and hence we have

$$\int_0^1 \frac{\log t}{\sqrt{t}(1+t)}\,dt = -\int_1^\infty \frac{\log t}{\sqrt{t}(1+t)}\,dt$$

More generally we have $\phi^{(2n+1)}(x) = -\phi^{(2n+1)}(1-x)$ and therefore

$$\phi^{(2n+1)}(1/2) = \int_0^\infty \frac{[\log t]^{2n+1}}{\sqrt{t}(1+t)} dt = 0$$

This is also proved in an entirely different manner by Srivastava and Choi in their book [126, p.268].

By differentiating Euler's reflection formula

$$\Gamma(x)\Gamma(1-x) = \frac{\pi}{\sin \pi x}$$

we obtain

$$\Gamma'(x)\Gamma(1-x) - \Gamma(x)\Gamma'(1-x) = -\frac{\pi^2 \cos \pi x}{\sin^2 \pi x}$$

and dividing by $\Gamma(x)\Gamma(1-x)$ we get

$$\frac{\Gamma'(x)}{\Gamma(x)} - \frac{\Gamma'(1-x)}{\Gamma(1-x)} = -\pi \cot \pi x$$

Or, alternatively, this equates to the well-known formula [126, p.14]

(C.42g) $$\psi(x) - \psi(1-x) = -\pi \cot \pi x$$

We also see that

$$\psi^{(n)}(x) - (-1)^n \psi^{(n)}(1-x) = -\pi \frac{d^n}{dx^n} \cot \pi x$$

and hence we have

$$2\psi^{(2n+1)}(x) = -\pi \frac{d^{2n+1}}{dx^{2n+1}} \cot \pi x$$

This gives us

$$2\psi^{(2n+1)}\left(\frac{1}{2}\right) = -\pi \frac{d^{2n+1}}{dx^{2n+1}} \cot \pi x \Big|_{x=1/2}$$

and therefore we have

$$\psi^{(2n+1)}\left(\frac{1}{2}\right) = (2n+1)!(2^{2n+2}-1)\varsigma(2n+2)$$

which then enables us to compute $\varsigma(2n)$. Derivatives of $\cot \pi x$ may be easily obtained by noting that if $y = \cot \pi x$ then

$$\tan^{-1}\left(\frac{1}{y}\right) = \pi x$$

and subsequently we obtain

$$\frac{dy}{dx} = -\pi(y^2 - 1)$$

(and we then apply the Leibniz rule).

We saw in (C.42g) that

$$\psi(x) - \psi(1-x) = -\pi \cot \pi x$$

and note for rational $x \in (0,1)$ that $\cot \pi x$ is algebraic; this then proves that $\psi(x) - \psi(1-x)$ is transcendental (and irrational) for all rational $x \in (0,1)$.

We have from the definition of the digamma function

$$\frac{d}{dx}\log\left[\Gamma(x)\Gamma(1-x)\right] = \psi(x) - \psi(1-x)$$

and obviously this implies that

$$\log\Gamma(x) + \log\Gamma(1-x) - \log\Gamma(a) - \log\Gamma(1-a) = \int_a^x \left[\psi(t) - \psi(1-t)\right] dt$$

$$= -\pi \int_a^x \cot \pi t \, dt$$

$$= \log\sin\pi a - \log\sin\pi x$$

Consider the integral

$$\text{(C.43)} \qquad u(x) = \int_x^{x+1} \log\Gamma(z) dz$$

Therefore we have

$$u'(x) = \log\Gamma(x+1) - \log\Gamma(x)$$

$$= \log x$$

Hence we get

$$u(x) = x\log x - x + c$$

Accordingly

$$u(0) = \int_0^1 \log\Gamma(z)dz = c$$

Letting $z \to 1-z$ we see that

$$\int_0^1 \log\Gamma(z)dz = \int_0^1 \log\Gamma(1-z)dz = c$$

and therefore we have

$$\int_0^1 \log\left[\Gamma(z)\Gamma(1-z)\right]dz = 2c$$

Then, using Euler's reflection formula $\Gamma(z)\Gamma(1-z) = \dfrac{\pi}{\sin\pi z}$ we get

$$\int_0^1 \log\left[\frac{\pi}{\sin\pi z}\right]dz = 2c$$

and thus

$$c = \frac{1}{2}\log\pi - \frac{1}{2}\int_0^1 \log\sin\pi z\,dz$$

We have from (3.2) $\int_0^{\pi/2} \log\sin x\,dx = -\dfrac{\pi}{2}\log 2$ and hence

$$\int_0^{\pi} \log\sin x\,dx = \int_0^{\pi/2} \log\sin x\,dx + \int_{\pi/2}^{\pi} \log\sin x\,dx = -\pi\log 2\,.$$

Therefore we have

$$c = \frac{1}{2}\log 2\pi$$

And accordingly we have

(C.43a) $$u(x) = \int_x^{x+1} \log\Gamma(z)dz = \frac{1}{2}\log 2\pi + x\log x - x$$

(C.43b) $$u(0) = \int_0^1 \log\Gamma(z)dz = \frac{1}{2}\log 2\pi$$

These identities were first shown by M. Raabe in Crelle's Journal [135, p.261] and the formula can be used to derive Stirling's asymptotic formula for $\Gamma(z)$ or $n!$.

The following analysis is taken from page 269 of "Exposé de la théorie, propriétés, des formules de transformation, et méthodes d'évaluation des intégrales définies" by David Bierens de Haan (1822-1895), a copy of which is available on the internet at the University of Michigan Historical Mathematics Collection.

Let us designate the following integrals as follows

$$I(x) = \int_0^1 \log\Gamma(x)\cos p\pi x\, dx$$

$$K(x) = \int_0^1 \log\Gamma(x)\sin p\pi x\, dx$$

Then we have by letting $x \to 1-x$

(C.44a) $$I(x) = \cos p\pi \int_0^1 \log\Gamma(1-x)\cos p\pi x\, dx + \sin p\pi \int_0^1 \log\Gamma(1-x)\sin p\pi x\, dx$$

$$= \cos p\pi I(1-x) + \sin p\pi K(1-x)$$

(C.44b) $$K(x) = \sin p\pi \int_0^1 \log\Gamma(1-x)\cos p\pi x\, dx + \cos p\pi \int_0^1 \log\Gamma(1-x)\sin p\pi x\, dx$$

$$= \sin p\pi I(1-x) - \cos p\pi K(1-x)$$

We have Euler's reflection formula

$$\log\Gamma(x) + \log\Gamma(1-x) = \log(2\pi) - \log(2\sin\pi x)$$

and, multiplying this equation by $\cos p\pi x$ and integrating, we obtain

$$\int_0^1 \log\Gamma(x)\cos p\pi x\, dx + \int_0^1 \log\Gamma(1-x)\cos p\pi x\, dx = \log(2\pi)\int_0^1 \cos p\pi x\, dx - \int_0^1 \log(2\sin\pi x)\cos p\pi x\, dx$$

Then, using the Fourier series (7.8) $\log\left[2\sin\left(x/2\right)\right] = -\sum_{n=1}^{\infty}\frac{\cos nx}{n}$ we get

$$I(x)+I(1-x)=\frac{\sin p\pi}{p\pi}\log(2\pi)+\sum_{n=1}^{\infty}\frac{1}{2n}\int_0^1\left[\cos\{(2n+p)\pi x\}+\cos\{(2n-p)\pi x\}\right]dx$$

$$=\frac{\sin p\pi}{p\pi}\log(2\pi)+\sum_{n=1}^{\infty}\frac{1}{2n}\left[\frac{\sin\{(2n+p)\pi\}}{(2n+p)\pi}+\frac{\sin\{(2n-p)\pi\}}{(2n-p)\pi}\right]$$

$$=\frac{\sin p\pi}{p\pi}\log(2\pi)-\frac{\pi\sin p\pi}{\pi}\sum_{n=1}^{\infty}\frac{1}{n}\frac{1}{4n^2-p^2}$$

Similarly, multiplying by $\sin p\pi x$ and integrating, we obtain

$$\int_0^1\log\Gamma(x)\sin p\pi x\,dx+\int_0^1\log\Gamma(1-x)\sin p\pi x\,dx=K(x)+K(1-x)$$

$$=\frac{1-\cos p\pi}{p\pi}\log(2\pi)+\sum_{n=1}^{\infty}\frac{1}{2n}\int_0^1\left[\sin\{(2n+p)\pi x\}-\sin\{(2n-p)\pi x\}\right]dx$$

$$=\frac{1-\cos p\pi}{p\pi}\log(2\pi)+\sum_{n=1}^{\infty}\frac{1}{2n}\left[\frac{1-\cos p\pi\{(2n+p)\pi\}}{(2n+p)\pi}+\frac{1-\cos p\pi\{(2n-p)\pi\}}{(2n-p)\pi}\right]$$

$$=\frac{1-\cos p\pi}{p\pi}\log(2\pi)-\frac{1-\cos p\pi}{\pi}\sum_{n=1}^{\infty}\frac{1}{n}\frac{1}{4n^2-p^2}$$

Denoting

$$P=\log(2\pi)-p\sum_{n=1}^{\infty}\frac{1}{n}\frac{1}{4n^2-p^2}$$

we have

(C.45a) $$I(x)+I(1-x)=\frac{\sin p\pi}{p\pi}P$$

(C.45b) $$K(x)+K(1-x)=\frac{1-\cos p\pi}{p\pi}P$$

Letting $p=2a$ where $a$ is an integer we get

$$I(x)+I(1-x)=\frac{\sin 2a\pi}{p\pi}\log(2\pi)+\sum_{n=1}^{\infty}\frac{1}{2n}\left[\frac{\sin\{(n+a)2\pi\}}{(n+p)2\pi}\right]+\sum_{n=1}^{a-1}\frac{\sin\{(a-n)2\pi\}}{(a-n)2\pi}$$

$$+\frac{1}{2a}\int_0^1\cos(0)\,dx+\sum_{n=a+1}^{\infty}\frac{\sin\{(a-n)2\pi\}}{(a-n)2\pi}$$

$$=\frac{1}{2a}$$

$$K(x)+K(1-x)=\frac{1-\cos p\pi}{p\pi}\log(2\pi)+\sum_{n=1}^{\infty}\frac{1}{2n}\left[\frac{1-\cos\{(n+a)2\pi\}}{(n+a)2\pi}\right]-\sum_{n=1}^{a-1}\frac{1}{2n}\left[\frac{1-\cos\{(a-n)2\pi\}}{(a-n)2\pi}\right]$$

$$-\frac{1}{2a}\int_0^1\sin(0)dx-\sum_{n=a+1}^{\infty}\frac{1}{2n}\left[\frac{1-\cos\{(a-n)2\pi\}}{(a-n)2\pi}\right]$$

$$=0$$

Accordingly, we obtain

$$I(x)=I(1-x)\quad\text{and}\quad \mathrm{K}(x)=-K(1-x)$$

Hence we have

(C.46) $$I(x)=\int_0^1\log\Gamma(x)\cos 2a\pi x\,dx=\int_0^1\log\Gamma(1-x)\cos 2a\pi x\,dx=\frac{1}{4a}$$

The following analysis is based primarily on Adamchik's 1997 paper, "A class of logarithmic integrals" [2a]. By using the binomial theorem we have (as mentioned previously in (4.4.112e))

(C.47) $$\lambda\int_0^z\frac{x^{p-1}}{\lambda+x^n}dx=\sum_{k=0}^{\infty}\frac{(-1)^k}{\lambda^k}\int_0^z x^{nk+p-1}dx\quad\left(\text{provided }\left|\frac{x^n}{\lambda}\right|<1\right)$$

$$=\sum_{k=0}^{\infty}\frac{(-1)^k z^{nk+p}}{\lambda^k(nk+p)}$$

With $z=1$ we have

(C.48) $$\lambda\int_0^1\frac{x^{p-1}}{\lambda+x^n}dx=\sum_{k=0}^{\infty}\frac{(-1)^k}{\lambda^k(nk+p)}$$

Therefore, from (C.48) we have

(C.49) $$\int_0^1\frac{x^{p-1}}{1+x}dx=\sum_{k=0}^{\infty}\frac{(-1)^k}{k+p}$$

Differentiating (C.48) $q-1$ times with respect to $p$ we obtain

(C.50) $$\int_0^1 \frac{x^{p-1}}{\lambda + x^n} \log^{q-1}\left(\frac{1}{x}\right) dx = \frac{\Gamma(q)}{\lambda} \sum_{k=0}^{\infty} \frac{(-1)^k}{\lambda^k (nk+p)^q}$$

In what seems to me a rather novel twist, Adamchik now regards $q$ as a continuous variable, rather than as discrete, and differentiating both sides of the new identity with respect to $q$ we obtain

(C.51)
$$\int_0^1 \frac{x^{p-1}}{\lambda + x^n} \log^{q-1}\left(\frac{1}{x}\right) \log\log\left(\frac{1}{x}\right) dx = \frac{\Gamma'(q)}{\lambda} \sum_{k=0}^{\infty} \frac{(-1)^k}{\lambda^k (nk+p)^q} - \frac{\Gamma(q)}{\lambda} \sum_{k=0}^{\infty} \frac{(-1)^k \log(nk+p)}{\lambda^k (nk+p)^q}$$

Letting $q = 1$ we get

(C.52) $$\int_0^1 \frac{x^{p-1}}{\lambda + x^n} \log\log\left(\frac{1}{x}\right) dx = -\frac{\gamma}{\lambda} \sum_{k=0}^{\infty} \frac{(-1)^k}{\lambda^k (nk+p)^q} - \frac{1}{\lambda} \sum_{k=0}^{\infty} \frac{(-1)^k \log(nk+p)}{\lambda^k (nk+p)}$$

Then substituting (C.50) we get

(C.53) $$\int_0^1 \frac{x^{p-1}}{\lambda + x^n} \log\log\left(\frac{1}{x}\right) dx = -\gamma \int_0^1 \frac{x^{p-1}}{\lambda + x^n} dx - \frac{1}{\lambda} \sum_{k=0}^{\infty} \frac{(-1)^k \log(nk+p)}{\lambda^k (nk+p)}$$

Curiously, at this stage, Adamchik comments that the infinite series in (C.53) is divergent when $\lambda = 1$ : as we shall shortly see the series can in fact be convergent.

We have previously referred to the Hurwitz-Lerch zeta function $\Phi(z,s,x)$ in (4.4.82): this is defined by [126, p.121] as

$$\Phi(z,s,x) = \sum_{n=0}^{\infty} \frac{z^n}{(n+x)^s}$$

We have

$$\frac{e^{-y(x-1)}}{e^y - z} = \sum_{n=0}^{\infty} z^n e^{-(x+n)y} \text{ and}$$

$$\int_0^{\infty} e^{-(x+n)y} y^{s-1} dy = \frac{1}{(x+n)^s} \int_0^{\infty} e^{-t} t^{s-1} dt = \frac{\Gamma(s)}{(x+n)^s}$$

and we have the integral representation [126, p.121]

$$\Phi(z,s,x) = \frac{1}{\Gamma(s)} \int_0^{\infty} \frac{y^{s-1} e^{-y(x-1)}}{e^y - z} dy$$

Adamchik then notes that

(C.54) $$\sum_{k=0}^{\infty}\frac{(-1)^k \log(nk+p)}{\lambda^k (nk+p)} = -\lim_{s\to 1}\frac{\partial}{\partial s}\sum_{k=0}^{\infty}\frac{(-1)^k}{\lambda^k (nk+p)^s} = -\lim_{s\to 1}\frac{\partial}{\partial s}\left[\frac{\Phi\left(-\frac{1}{\lambda}, s, \frac{p}{n}\right)}{n^s}\right]$$

With $\lambda = 1$ we have

$$\sum_{k=0}^{\infty}\frac{(-1)^k \log(nk+p)}{(nk+p)} = -\lim_{s\to 1}\frac{\partial}{\partial s}\left[\frac{\Phi\left(-1, s, \frac{p}{n}\right)}{n^s}\right]$$

We have

$$\Phi\left(-1, s, \frac{p}{n}\right) = \frac{1}{2^s}\left[\varsigma\left(s, \frac{p}{2n}\right) - \varsigma\left(s, \frac{1}{2}+\frac{p}{2n}\right)\right]$$

and we therefore get with $\lambda = 1$

$$\sum_{k=0}^{\infty}\frac{(-1)^k \log(nk+p)}{(nk+p)} = -\lim_{s\to 1}\frac{\partial}{\partial s}\frac{1}{(2n)^s}\left[\varsigma\left(s, \frac{p}{2n}\right) - \varsigma\left(s, \frac{1}{2}+\frac{p}{2n}\right)\right]$$

$$= -\lim_{s\to 1}\left\{\frac{1}{(2n)^s}\left[\varsigma'\left(s, \frac{p}{2n}\right) - \varsigma'\left(s, \frac{1}{2}+\frac{p}{2n}\right)\right] - \frac{\log(2n)}{(2n)^s}\left[\varsigma\left(s, \frac{p}{2n}\right) - \varsigma\left(s, \frac{1}{2}+\frac{p}{2n}\right)\right]\right\}$$

Using the relationship $\varsigma(s,a) \to -\gamma - \psi(a) + \varsigma(s)$ as $s \to 1$ (which is proved in (4.3.228c) in Volume II(b)) we obtain

(C.55)
$$\sum_{k=0}^{\infty}\frac{(-1)^k \log(nk+p)}{(nk+p)} = -\frac{1}{2n}\left[\varsigma'\left(1, \frac{p}{2n}\right) - \varsigma'\left(1, \frac{1}{2}+\frac{p}{2n}\right)\right] - \frac{\log(2n)}{2n}\left[\psi\left(\frac{p}{2n}\right) - \psi\left(\frac{1}{2}+\frac{p}{2n}\right)\right]$$

$$= -\frac{1}{2n}\left[\varsigma'\left(1, \frac{p}{2n}\right) - \varsigma'\left(1, \frac{n+p}{2n}\right)\right] - \frac{\log(2n)}{2n}\left[\psi\left(\frac{p}{2n}\right) - \psi\left(\frac{n+p}{2n}\right)\right]$$

We also need to consider the other integral in (C.53) which we may transform by the obvious substitution $y = x^n$

(C.56) $$\int_0^1 \frac{x^{p-1}}{1+x^n}dx = \frac{1}{n}\int_0^1 \frac{y^{\frac{p}{n}-1}}{1+y}dy = \frac{1}{2n}\left[\psi\left(\frac{n+p}{2}\right) - \psi\left(\frac{p}{2n}\right)\right]$$

Adamchik finally arrives at the rather fearsome integral which we have referred to previously in (4.4.44g)

(C.57)

$$\int_0^1 \frac{x^{p-1}}{1+x^n}\log\log\left(\frac{1}{x}\right)dx = \frac{\gamma+\log(2n)}{2n}\left[\psi\left(\frac{p}{2n}\right)-\psi\left(\frac{n+p}{2n}\right)\right]+\frac{1}{2n}\left[\varsigma'\left(1,\frac{p}{2n}\right)-\varsigma'\left(1,\frac{n+p}{2n}\right)\right]$$

where $\operatorname{Re}(p) > 0$, $\operatorname{Re}(n) > 0$. This equation was used in (6.107rvii) in Volume V.

Several integrals are then evaluated by Adamchik using this formula; for example we have

(C.58) $$\int_0^1 \frac{x^{n-1}}{1+x^n}\log\log\left(\frac{1}{x}\right)dx = -\frac{\log(2)\log(2n^2)}{2n}$$

Hence, with $n=1$ we have

(C.59) $$\int_0^1 \frac{1}{1+x}\log\log\left(\frac{1}{x}\right)dx = -\frac{\log^2 2}{2}$$

Referring back to (C.53) with $n = p = 1$ we get

(C.60) $$\int_0^1 \frac{1}{1+x}\log\log\left(\frac{1}{x}\right)dx = -\gamma\int_0^1 \frac{1}{1+x}dx - \sum_{k=0}^{\infty}\frac{(-1)^k\log(k+1)}{(k+1)}$$

$$= -\gamma\log 2 + \sum_{k=1}^{\infty}\frac{(-1)^k\log k}{k}$$

This integral is contained in G&R [74, p.565] together with several other related integrals.

Therefore, we have the well-known result

(C.61) $$\sum_{k=1}^{\infty}\frac{(-1)^k\log k}{k} = \log 2\left[\gamma-\frac{\log 2}{2}\right] = \varsigma_a'(1)$$

We also note another result from Adamchik [2a]

$$\int_0^1 \frac{x^{p-1}}{\left(1+x^n\right)^2}\log\log\left(\frac{1}{x}\right)dx = \frac{(n-p)[\gamma+\log(2n)]}{2n^2}\left[\psi\left(\frac{p}{2n}\right)-\psi\left(\frac{n+p}{2n}\right)\right]$$

$$-\frac{1}{2n}\left[\gamma+\log(2n)-2\log\frac{\Gamma\left(\frac{p}{2n}\right)}{\Gamma\left(\frac{n+p}{2n}\right)}\right]+\frac{(n-p)}{2n^2}\left[\varsigma'\left(1,\frac{p}{2n}\right)-\varsigma'\left(1,\frac{n+p}{2n}\right)\right]$$

It may be noted that Roach [110b] provided alternative proofs of Adamchik's formulae in 2005.

Differentiating (C.56) results in

$$\int_0^1 \frac{x^{p-1}\log x}{1+x^n}dx = \frac{1}{4n}\left[\psi'\left(\frac{n+p}{2}\right)-\psi'\left(\frac{p}{2n}\right)\right]$$

The integral $\int_0^\infty \frac{\log x}{e^x+1}dx$ was posed as a problem by Sandham [118aa] in 1950 and the solution provided by Trost proceeded as follows. Consider the integral

$$I(a,\varepsilon) = \int_\varepsilon^\infty \frac{\log x}{e^{ax}+1}dx = -\frac{1}{a}\int_\varepsilon^\infty \frac{(-ae^{-ax})\log x}{1+e^{-ax}}dx$$

and integration by parts produces

$$= \frac{1}{a}\log\varepsilon\log\left[1+e^{-a\varepsilon}\right]+\frac{1}{a}\int_\varepsilon^\infty \frac{\log\left[1+e^{-ax}\right]}{x}dx$$

Since $\log[1+x] = \log\left[1-x^2\right]-\log[1-x]$ we have

$$K = \int_\varepsilon^\infty \frac{\log\left[1+e^{-ax}\right]}{x}dx = \int_\varepsilon^\infty \frac{\log\left[1-e^{-2ax}\right]}{x}dx - \int_\varepsilon^\infty \frac{\log\left[1-e^{-ax}\right]}{x}dx$$

$$= \int_{2\varepsilon}^\infty \frac{\log\left[1-e^{-ax}\right]}{x}dx - \int_\varepsilon^\infty \frac{\log\left[1-e^{-ax}\right]}{x}dx$$

$$= -\int_\varepsilon^{2\varepsilon} \frac{\log\left[1-e^{-ax}\right]}{x}dx$$

$$= -\int_1^2 \frac{\log\left[1-e^{-a\varepsilon x}\right]}{x}dx$$

We now write

$$\log\left[1-e^{-a\varepsilon x}\right] = \log\left[1-\left(1-\frac{a\varepsilon x}{1!}+\frac{(a\varepsilon x)^2}{2!}-\ldots\right)\right] = \log\left[a\varepsilon x\left(1-\frac{a\varepsilon x}{2!}+\frac{(a\varepsilon x)^2}{3!}-\ldots\right)\right]$$

and the integral becomes

$$= -\int_1^2 \frac{\log\left[a\varepsilon x\right]}{x} dx - \int_1^2 \frac{\log\left[1 - \frac{a\varepsilon x}{2!} + \frac{(a\varepsilon x)^2}{3!} - \ldots\right]}{x} dx$$

We have

$$-\int_1^2 \frac{\log\left[a\varepsilon x\right]}{x} dx = -\int_1^2 \frac{\log x}{x} dx - \log\left[a\varepsilon\right] \int_1^2 \frac{dx}{x}$$

$$= -\frac{1}{2}\log^2 2 - \log[a\varepsilon]\log 2$$

Using the Maclaurin expansion of $\log(1-u)$ we readily see that

$$-\int_1^2 \frac{\log\left[1 - \frac{a\varepsilon x}{2!} + \frac{(a\varepsilon x)^2}{3!} - \ldots\right]}{x} dx = \delta(\varepsilon)$$

where $\delta(\varepsilon) \to 0$ as $\varepsilon \to 0$.

Therefore we have

$$K = -\log\varepsilon \log 2 - \frac{1}{2}\log^2 2 - \log a \log 2 - \delta(\varepsilon)$$

We now have

$$I(a,\varepsilon) = \frac{1}{a}\log\varepsilon \log\left[\left(1+e^{-a\varepsilon}\right)/2\right] - \frac{1}{2a}\log 2 \log^2(2a) - \frac{1}{a}\delta(\varepsilon)$$

Then as $\varepsilon \to 0$ we get

(C.61a) $\qquad I(a,0) = -\frac{1}{2a}\log 2 \log(2a^2)$ and $I(1,0) = -\frac{1}{2}\log^2 2$

See also a related problem posed by Klamkin [86a] in 1954.

From (4.4.42) we have

$$\varsigma_a(s) = \frac{1}{\Gamma(s)}\int_0^\infty \frac{u^{s-1}}{e^u+1} du$$

Differentiating with respect to $s$ we obtain

$$\varsigma_a'(s) = -\frac{\Gamma'(s)}{\left[\Gamma(s)\right]^2}\int_0^\infty \frac{u^{s-1}}{e^u+1} du + \frac{1}{\Gamma(s)}\int_0^\infty \frac{u^{s-1}\log u}{e^u+1} du$$

Therefore we get

$$\varsigma_a'(1) = \gamma\int_0^\infty \frac{du}{e^u+1} + \int_0^\infty \frac{\log u}{e^u+1}du$$

and since

$$\int_0^\infty \frac{du}{e^u+1} = -\log\left(1+e^{-u}\right)\Big|_0^\infty = \log 2$$

we end up with a proof of (C.61).

We also have

$$\varsigma_a''(s) = -\frac{\left[\Gamma(s)\right]^2\Gamma''(s) - 2\left[\Gamma'(s)\right]^2\Gamma(s)}{\left[\Gamma(s)\right]^2}\int_0^\infty \frac{u^{s-1}}{e^u+1}du$$

$$-2\frac{\Gamma'(s)}{\left[\Gamma(s)\right]^2}\int_0^\infty \frac{u^{s-1}\log u}{e^u+1}du + \frac{1}{\Gamma(s)}\int_0^\infty \frac{u^{s-1}\log^2 u}{e^u+1}du$$

and hence we get using (E.16d)

$$\varsigma_a''(1) = \left[\gamma^2 - \varsigma(2)\right]\int_0^\infty \frac{du}{e^u+1} + 2\gamma\int_0^\infty \frac{\log u}{e^u+1}du + \int_0^\infty \frac{\log^2 u}{e^u+1}du$$

$$= \left[\gamma^2 - \varsigma(2)\right]\log 2 - \gamma\log^2 2 + \int_0^\infty \frac{\log^2 u}{e^u+1}du$$

where we have used (C.61). This then becomes

(C.61b) $$\varsigma_a''(1) = \left[\gamma^2 - \varsigma(2) - \gamma\log 2\right]\log 2 + \int_0^\infty \frac{\log^2 u}{e^u+1}du$$

With the substitution $x = e^{-u}$ this becomes (C.68) which was derived in a rather different manner.

This time we consider the integral

$$J(\varepsilon) = \int_\varepsilon^\infty \frac{\log^2 x}{e^x+1}dx = -\int_\varepsilon^\infty \frac{(-e^{-x})\log^2 x}{1+e^{-x}}dx$$

and integration by parts produces

$$= \log^2\varepsilon\log\left[1+e^{-\varepsilon}\right] + 2\int_\varepsilon^\infty \frac{\log\left[1+e^{-x}\right]\log x}{x}dx$$

We have

$$L = \int_{\varepsilon}^{\infty} \frac{\log\left[1+e^{-x}\right]}{x} \log x\,dx = \int_{\varepsilon}^{\infty} \frac{\log\left[1-e^{-2x}\right]}{x} \log x\,dx - \int_{\varepsilon}^{\infty} \frac{\log\left[1-e^{-x}\right]}{x} \log x\,dx$$

$$= \int_{2\varepsilon}^{\infty} \frac{\log\left[1-e^{-x}\right]}{x} \log(x/2)\,dx - \int_{\varepsilon}^{\infty} \frac{\log\left[1-e^{-x}\right]}{x} \log x\,dx$$

$$= -\int_{\varepsilon}^{2\varepsilon} \frac{\log\left[1-e^{-x}\right]}{x} \log x\,dx - \log 2 \int_{2\varepsilon}^{\infty} \frac{\log\left[1-e^{-x}\right]}{x}\,dx$$

$$= -\int_{1}^{2} \frac{\log\left[1-e^{-\varepsilon x}\right]}{x} \log x\,dx - \log\varepsilon \int_{1}^{2} \frac{\log\left[1-e^{-\varepsilon x}\right]}{x}\,dx - \log 2 \int_{2\varepsilon}^{\infty} \frac{\log\left[1-e^{-x}\right]}{x}\,dx$$

As previously, we note that

$$\lambda(\varepsilon) = \int_{1}^{2} \frac{\log\left[1-\frac{\varepsilon x}{2!}+\frac{(\varepsilon x)^2}{3!}-...\right]}{x} \log x\,dx \to 0 \;\text{ as }\; \varepsilon \to 0$$

We easily see that

$$-\int_{1}^{2} \frac{\log[\varepsilon x]}{x} \log x\,dx = -\frac{1}{3}\log^3 2 - \frac{1}{2}\log\varepsilon \log^2 2$$

We have already shown that

$$-\log\varepsilon \int_{1}^{2} \frac{\log\left[1-e^{-\varepsilon x}\right]}{x}\,dx = -\log\varepsilon\left[-\log\varepsilon\log 2 - \frac{1}{2}\log^2 2 - \delta(\varepsilon)\right]$$

and, as before, we have

$$-\int_{1}^{2} \frac{\log\left[1-e^{-\varepsilon x}\right]}{x} \log x\,dx = -\int_{1}^{2} \frac{\log[\varepsilon x]}{x} \log x\,dx - \int_{1}^{2} \frac{\log\left[1-\frac{\varepsilon x}{2!}+\frac{(\varepsilon x)^2}{3!}-...\right]}{x} \log x\,dx$$

Integration by parts results in

$$\int_{2\varepsilon}^{\infty} \frac{\log x}{e^x-1}\,dx = -\log x \log\left[1-e^{-x}\right]\Big|_{2\varepsilon}^{\infty} + \int_{2\varepsilon}^{\infty} \frac{\log\left[1-e^{-x}\right]}{x}\,dx$$

and therefore we get

$$\int_{2\varepsilon}^{\infty} \frac{\log\left[1-e^{-x}\right]}{x} dx = \int_{2\varepsilon}^{\infty} \frac{\log x}{e^x - 1} dx - \log 2\varepsilon \log\left[1-e^{-2\varepsilon}\right]$$

We note from (4.4.44n) that

$$\varsigma_a''(1) = 2\gamma_1 \log 2 - \gamma \log^2 2 + \frac{1}{3}\log^2 2$$

It is not clear to me whether further work in this area will produce any meaningful results.

From (C.50) we have with $\lambda = p = 1,\ q = 2$ and $n = 3$

(C.62) $$\int_0^1 \frac{\log x}{1+x^3} dx = \sum_{k=0}^{\infty} \frac{(-1)^{k+1}}{(3k+1)^2}$$

$$= -\frac{1}{1^2} + \frac{1}{4^2} - \frac{1}{7^2} + \frac{1}{11^2} - ....$$

Differentiating (C.56) gives us

(C.62a) $$\int_0^1 \frac{x^{p-1} \log x}{1+x^n} dx = \frac{1}{4n}\left[\psi'\left(\frac{n+p}{2}\right) - \psi'\left(\frac{p}{2n}\right)\right]$$

and with $p = 1$ and $n = 3$ we get

(C.62b) $$\int_0^1 \frac{\log x}{1+x^3} dx = \frac{1}{12}\left[\psi'(2) - \psi'\left(\frac{1}{6}\right)\right]$$

Recall that in (8.40) of Volume V it was reported that

(C.63) $$\frac{1}{1^2} + \frac{1}{2^2} - \frac{1}{4^2} - \frac{1}{5^2} + \frac{1}{7^2} + \frac{1}{8^2} - ... = -\frac{2}{27}\pi^2 - 2\int_0^1 \frac{\log x}{1+x^3} dx$$

and we therefore get

(C.64) $$\frac{1}{1^2} - \frac{1}{2^2} - \frac{1}{4^2} + \frac{1}{5^2} - \frac{1}{7^2} - \frac{1}{8^2} + ... = \frac{2}{27}\pi^2$$

We also note from (6.30y) that

$$\frac{1}{1^2} + \frac{1}{3^2} - \frac{1}{5^2} - \frac{1}{7^2} + \frac{1}{9^2} + ... = \sqrt{2}\phi\left(\sqrt{2}-1\right) + \frac{\pi}{4\sqrt{2}} \log\left(1+\sqrt{2}\right)$$

but nothing particularly miraculous currently springs to mind!

Let us now differentiate (C.51) with respect to $q$ to obtain

(C.65)

$$\int_0^1 \frac{x^{p-1}}{1+x^n}\log^{q-1}\left(\frac{1}{x}\right)\left[\log\log\left(\frac{1}{x}\right)\right]^2 dx = \Gamma''(q)\sum_{k=0}^{\infty}\frac{(-1)^k}{(nk+p)^q} - 2\Gamma'(q)\sum_{k=0}^{\infty}\frac{(-1)^k\log(nk+p)}{(nk+p)^q}$$

$$+\Gamma(q)\sum_{k=0}^{\infty}\frac{(-1)^k\left[\log(nk+p)\right]^2}{(nk+p)^q}$$

with $p = n = 1,$ we get

(C.66)

$$\int_0^1 \frac{1}{1+x}\log^{q-1}\left(\frac{1}{x}\right)\left[\log\log\left(\frac{1}{x}\right)\right]^2 dx = \Gamma''(q)\sum_{k=0}^{\infty}\frac{(-1)^k}{(k+1)^q} - 2\Gamma'(q)\sum_{k=0}^{\infty}\frac{(-1)^k\log(k+1)}{(k+1)^q}$$

$$+\Gamma(q)\sum_{k=0}^{\infty}\frac{(-1)^k\left[\log(k+1)\right]^2}{(k+1)^q}$$

and this simplifies to the pretty formula

(C.67)

$$\int_0^1 \frac{1}{1+x}\log^{q-1}\left(\frac{1}{x}\right)\left[\log\log\left(\frac{1}{x}\right)\right]^2 dx = \Gamma''(q)\varsigma_a(q) + 2\Gamma'(q)\varsigma_a'(q) + \Gamma(q)\varsigma_a''(q)$$

Letting $q = 1$ and using (4.2.4) and (E.16d) we obtain

(C.68) $$\int_0^1 \frac{1}{1+x}\left[\log\log\left(\frac{1}{x}\right)\right]^2 dx = \left[\gamma^2 + \varsigma(2)\right]\log 2 - 2\gamma\varsigma_a'(1) + \varsigma_a''(1)$$

$$= \left[\gamma^2 + \varsigma(2)\right]\log 2 - 2\gamma\log 2\left[\gamma - \frac{\log 2}{2}\right] + \varsigma_a''(1)$$

$$= \left[-\gamma^2 + \varsigma(2) + \gamma\log 2\right]\log 2 + \varsigma_a''(1)$$

Adamchik [2a] also showed that

(C.69) $$\int_0^1 x^{n-1}\log\log\left(\frac{1}{x}\right)dx = -\frac{\log n+\gamma}{n}$$

and completing the summation we obtain

$$\sum_{n=1}^{\infty}\frac{1}{n^s}\int_0^1 x^{n-1}\log\log\left(\frac{1}{x}\right)dx = -\sum_{n=1}^{\infty}\frac{\log n}{n^{s+1}}-\gamma\varsigma(s+1)=\varsigma'(s+1)-\gamma\varsigma(s+1)$$

Hence we have

(C.70) $$\int_0^1 Li_s(x)\log\log\left(\frac{1}{x}\right)\frac{dx}{x}=\varsigma'(s+1)-\gamma\varsigma(s+1)$$

and with the substitution $x=e^{-u}$ this becomes

(C.70) $$\int_0^{\infty} Li_s(e^{-u})\log u\,du = \varsigma'(s+1)-\gamma\varsigma(s+1)$$

(which was proved in a more elementary manner in (4.4.187).

With $n=2\pi p$ in (C.69) we note that

$$\int_0^1 x^{2\pi p-1}\log\log\left(\frac{1}{x}\right)dx = -\frac{\log 2p\pi+\gamma}{2p\pi}$$

and we will see in (E.46) that

$$\int_0^1 \log\Gamma(x)\sin 2p\pi x\,dx = \frac{\log 2p\pi+\gamma}{2p\pi}$$

and therefore we have (coincidentally?)

(C.71) $$\int_0^1 x^{2\pi p-1}\log\log\left(\frac{1}{x}\right)dx = -\int_0^1 \log\Gamma(x)\sin 2p\pi x\,dx$$

With regard to (C.69) we could also make the summation

$$\sum_{n=1}^{\infty}(-1)^n\int_0^1 x^{n-1}\log\log\left(\frac{1}{x}\right)dx = -\sum_{n=1}^{\infty}(-1)^n\frac{\log n}{n}+\gamma\log 2 = -\varsigma_a'(1)+\gamma\log 2$$

which gives us back (C.59)

(C.72) $$\int_0^1 \frac{1}{1+x}\log\log\left(\frac{1}{x}\right)dx = \varsigma_a'(1) - \gamma\log 2 = -\frac{1}{2}\log^2 2$$

## APPENDIX D

## A VERY ELEMENTARY PROOF OF $\frac{\pi^2}{8} = \sum_{n=0}^{\infty}\frac{1}{(2n+1)^2}$

We can easily verify the following identity by multiplying the numerator and denominator by the complex conjugate $(1+e^{-i\alpha x})$

(D.1) $$\frac{1}{1+e^{i\alpha x}} = \frac{1+e^{-i\alpha x}}{2+e^{i\alpha x}+e^{-i\alpha x}} = \frac{1+\cos\alpha x - i\sin\alpha x}{2(1+\cos\alpha x)} = \frac{1}{2} - \frac{i}{2}\frac{\sin\alpha x}{1+\cos\alpha x}$$

where we have Euler's identity $e^{i\alpha x} = \cos\alpha x + i\sin\alpha x$.

Using the geometric series we have

(D.2) $$\frac{1}{1-y} = \sum_{n=0}^{N} y^n + \frac{y^{N+1}}{1-y}$$

and letting $y = -e^{i\alpha x}$ we obtain

(D.3) $$\frac{1}{1+e^{i\alpha x}} = \sum_{n=0}^{N}(-1)^n (e^{i\alpha x})^n + \frac{(-1)^{N+1}(e^{i\alpha x})^{N+1}}{1+e^{i\alpha x}}$$

(D.4) $$= \sum_{n=0}^{N}(-1)^n e^{i\alpha nx} + \frac{(-1)^{N+1}e^{i\alpha(N+1)x}}{1+e^{i\alpha x}}$$

We have

(D.5) $$\frac{(-1)^{N+1}e^{i\alpha(N+1)x}}{1+e^{i\alpha x}} = \frac{(-1)^{N+1}e^{i\alpha(N+1)x}}{e^{i\alpha x/2}\left\{e^{-i\alpha x/2}+e^{i\alpha x/2}\right\}} = \frac{(-1)^{N+1}e^{i\alpha(N+1)x}e^{-i\alpha x/2}}{2\cos(\alpha x/2)}$$

(D.6) $$= \frac{(-1)^{N+1}e^{i\alpha\left(N+\frac{1}{2}\right)x}}{2\cos(\alpha x/2)} = \frac{(-1)^{N+1}\left\{\cos\alpha\left(N+\frac{1}{2}\right)x + i\sin\alpha\left(N+\frac{1}{2}\right)x\right\}}{2\cos(\alpha x/2)}$$

Separating the real and imaginary parts of (D.1) and (D.4) we obtain

(D.7a) $$\frac{1}{2} = \sum_{n=0}^{N} (-1)^n \cos \alpha n x + \frac{(-1)^{N+1} \cos \alpha \left( N + \frac{1}{2} \right) x}{2 \cos(\alpha x / 2)}$$

(D.7b) $$-\frac{1}{2} \frac{\sin \alpha x}{1 + \cos \alpha x} = \sum_{n=0}^{N} (-1)^n \sin \alpha n x + \frac{(-1)^{N+1} \sin \alpha \left( N + \frac{1}{2} \right) x}{2 \cos(\alpha x / 2)}$$

With $\alpha = 2$, we now multiply (D.7a) by $(\pi / 2 - x)$ and integrate over the interval $[0, \pi / 2]$ to obtain

(D.8)

$$\int_0^{\pi/2} \frac{(\pi / 2 - x)}{2} dx = \int_0^{\pi/2} \sum_{n=0}^{N} (-1)^n (\pi / 2 - x) \cos 2nx \, dx + \int_0^{\pi/2} \frac{(-1)^{N+1} (\pi / 2 - x) \cos(2N+1) x}{2 \cos x} dx$$

A straightforward integration by parts gives

(D.9) $$\int_0^{\pi/2} (\pi / 2 - x) \cos 2nx \, dx = \frac{\pi \sin 2nx}{2n} - \frac{x \sin 2nx}{2n} - \frac{\cos 2nx}{4n^2} \Bigg|_0^{\pi/2} \qquad , n \geq 1$$

Accordingly we obtain

(D.10) $$\frac{\pi^2}{16} = \frac{1}{4} \sum_{n=1}^{N} \frac{1}{n^2} - \frac{1}{4} \sum_{n=1}^{N} \frac{(-1)^n}{n^2} + \frac{\pi^2}{8} + \int_0^{\pi/2} \frac{(-1)^{N+1} (\pi / 2 - x) \cos(2N+1) x}{2 \cos x} dx$$

having remembered to include the term for $n = 0$.

We now consider the nature of the remainder term.

(D.11) $$\int_0^{\pi/2} \frac{(-1)^{N+1} (\pi / 2 - x) \cos(2N+1) x}{2 \cos x} dx = \int_0^{\pi/2} g(x) \cos(2N+1) x \, dx = R_N$$

where $$g(x) = \frac{(-1)^{N+1} (\pi / 2 - x)}{2 \cos x}$$

Using integration by parts we have

(D.12) $$\int_0^{\pi/2} g(x) \cos(2N+1) x \, dx = g(x) \frac{\sin(2N+1) x}{2N+1} \Bigg|_0^{\pi/2} - \frac{1}{2N+1} \int_0^{\pi/2} g'(x) \sin(2N+1) x \, dx$$

From L'Hôpital's theorem it is easily seen that

$$\lim_{x\to\pi/2}\frac{\pi/2-x}{\cos x}=\lim_{x\to\pi/2}\frac{-1}{\sin x}=-1$$

and hence the integrated term in (D.12) is finite. We have

(D.13) $$I_N=-\frac{1}{2N+1}\int_0^{\pi/2} g'(x)\sin\left(2N+1\right)x\,dx$$

Upon differentiation we have

$$g'(x)=\frac{(-1)^{N+1}\left[-\cos x+(\pi/2-x)\sin x\right]}{4\cos^2 x}$$

and substituting this in (D.13) we obtain

(D.14) $$I_N=-\frac{1}{2N+1}\int_0^{\pi/2}\frac{(-1)^{N+1}\left[-\cos x+(\pi/2-x)\sin x\right]}{4\cos^2 x}\sin\left(2N+1\right)x\,dx$$

(D.14a) $$=-\frac{(-1)^{N+1}}{4(2N+1)}\int_0^{\pi/2}\left[-1+\frac{(\pi/2-x)\sin x}{\cos x}\right]\frac{\sin\left(2N+1\right)x}{\cos x}dx$$

Using L'Hôpital's theorem we can see that

$$\lim_{x\to\pi/2}\frac{\sin\left(2N+1\right)x}{\cos x}=\lim_{x\to\pi/2}\frac{\left(2N+1\right)\cos\left(2N+1\right)x}{\sin x}=0$$

and hence the integrand in (D.14a) is finite throughout the interval $[0,\pi/2]$. Therefore the integral must be finite and accordingly we have

$$\lim_{N\to\infty}R_N=0$$

From (D.10) we therefore conclude that

$$\frac{\pi^2}{4}=\lim_{N\to\infty}\left[\sum_{n=1}^{N}\frac{1}{n^2}+\sum_{n=1}^{N}\frac{(-1)^{n+1}}{n^2}+R_N\right]$$

$$=\sum_{n=1}^{\infty}\frac{1}{n^2}+\sum_{n=1}^{\infty}\frac{(-1)^{n+1}}{n^2}$$

$$=\sum_{n=1}^{\infty}\left[\frac{1}{n^2}+\frac{(-1)^{n+1}}{n^2}\right]$$

$$= 2\sum_{n=0}^{\infty} \frac{1}{(2n+1)^2}$$

Having been presented with this example, a class could be requested to investigate the topic further by reference to (i) different intervals of integration and (ii) different appropriate functions for $p(x)$ in the style of Section 6 of Volume V. An inquisitive student will certainly enjoy the results of her research.

## APPENDIX E

## SOME ASPECTS OF EULER'S CONSTANT $\gamma$ AND THE GAMMA FUNCTION

Consider the sequence $u_n$ defined by

(E.1) $$u_n = \int_0^1 \frac{x}{n(n+x)} dx = \int_0^1 f_n(x)\, dx$$

Since $f_n(x) > 0 \;\forall n \geq 1$ and $x \in (0,1]$, it is clear that $u_n > 0$. By differentiation we see that

$$f_n'(x) = \frac{1}{(n+x)^2} > 0$$

and hence $f_n(x)$ is a monotonic increasing function. The maximum value of $f_n(x)$ in the interval $x \in [0,1]$ is therefore equal to $1/n(n+1)$ and hence

$$0 < u_n = \int_0^1 \frac{x}{n(n+x)} dx < \frac{1}{n(n+1)}$$

We may integrate (E.1) using partial fractions

(E.2) $$u_n = \int_0^1 \frac{x}{n(n+x)} dx = \int_0^1 \left( \frac{1}{n} - \frac{1}{n+x} \right) dx = \frac{1}{n} - \log \frac{n+1}{n}$$

The following finite sum and inequality are easily derived

(E.3) $$\sum_{n=1}^{N} u_n = H_N - \log\left(N+1\right) < \sum_{n=1}^{N} \frac{1}{n(n+1)}$$

We have the familiar telescoping series

(E.4) $$\sum_{n=1}^{N}\frac{1}{n(n+1)}=\sum_{n=1}^{N}\left(\frac{1}{n}-\frac{1}{(n+1)}\right)=1-\frac{1}{N+1}$$

and accordingly we have

(E.5) $$0<\sum_{n=1}^{N}u_n<1-\frac{1}{N+1}$$

Therefore $\sum_{n=1}^{\infty}u_n$ converges because $\sum_{n=1}^{N}u_n$ is monotonic increasing and is bounded above.

We have the obvious identity

(E.5a) $$H_N-\log N=\sum_{n=1}^{N}u_n+\log\frac{N+1}{N}$$

and therefore we get

(E.6) $$\lim_{N\to\infty}\left(H_N-\log N\right)=\lim_{N\to\infty}\left(\sum_{n=1}^{N}u_n+\log\frac{N+1}{N}\right)=\sum_{n=1}^{\infty}u_n$$

The limit is the well-known Euler's constant which is usually designated by $\gamma$ (it is however designated by $C$ in Gradshteyn and Ryzhik's "Tables of Integrals, Series and Products" [74], and in some other works). Sometimes it is referred to as the Euler-Mascheroni constant [78, p.89] (in honour of the Italian mathematician Mascheroni (1750-1800) who computed 32 digits of $\gamma$ in 1790: it was later ascertained in 1809 that only the first 19 digits were correct!).

Euler's constant is therefore defined as

(E.6) $$\gamma=\lim_{n\to\infty}\left(H_n-\log n\right)=\sum_{n=1}^{\infty}\int_0^1\frac{x}{n(n+x)}dx$$

Using (E.3) it is also easily seen that

(E.6a) $$\gamma=\lim_{n\to\infty}\left(H_n-\log[n+1]\right)=\lim_{n\to\infty}\left(H_n-\log[n+\alpha]\right)$$

From (E.5) we know that $\gamma$ satisfies the inequality $0<\gamma<1$: in 2003 Chao-Ping Chen and Feng Qi [43a] showed that

(E.6b) $$\frac{1}{2n+\beta}\le H_n-\log n-\gamma<\frac{1}{2n+\alpha}$$

where $\alpha = \frac{1}{1-\gamma} - 2$ , $\beta = \frac{1}{3}$ , and the constants $\alpha$ and $\beta$ are the best possible.
By squaring the terms in (E.6b) and applying the squeeze theorem, it is easily seen that

(E.6d) $$\lim_{n\to\infty}\left[\left(H_n\right)^2 - 2H_n \log n + \log^2 n\right] = \gamma^2$$

Such limits are considered in more detail later in this Appendix.

We have the approximation

(E.6e) $$\gamma = 0.5772157....$$

It is not yet known whether $\gamma$ is irrational or transcendental. By computing a large number of digits of $\gamma$ and using a continued fraction expansion, T. Papanikolaou [25, p.185] showed in 1997 that if $\gamma$ is a rational number then its denominator must have at least 242,080 digits. I suspect that it is fairly safe to assume that $\gamma$ is irrational!

From (E.2) we see that

$$\frac{1}{n} > \log\frac{n+1}{n}$$

and therefore $H_N > \log(n+1)$ (and this easily shows that the harmonic series diverges).

In connection with (E.6) above, in 2002 Yingying [139a] showed that

(E.6f) $$\sum_{n=k+1}^{\infty}\int_0^1 \frac{x}{n(n+x)}dx = \sum_{n=1}^{\infty}\frac{a_n}{(k+1)...(k+n)}$$

where

(E.6g) $$a_1 = \frac{1}{2} \qquad a_n = \frac{1}{n}\int_0^1 t(1-t)...(n-1-t)dt = \frac{(-1)^{n+1}}{n}\sum_{j=1}^{n}\frac{s(n,j)}{j+1}$$

and $s(n, j)$ are the Stirling numbers of the first kind.

□

We have

$$\sum_{k=1}^{n}\log\left(1+\frac{1}{k}\right) = \sum_{k=1}^{n}[\log(k+1) - \log k]$$

and this telescopes to

$$\sum_{k=1}^{n} \log\left(1+\frac{1}{k}\right) = \log(n+1)$$

We therefore can write

(E.6h) $$\gamma = \lim_{n\to\infty} \sum_{k=1}^{n} \left[\frac{1}{k} - \log\left(1+\frac{1}{k}\right)\right]$$

It is interesting to note that Sondow [123aa] has discovered a similar alternating series

(E.6i) $$\log\frac{4}{\pi} = \lim_{n\to\infty} \sum_{k=1}^{n} (-1)^{k+1} \left[\frac{1}{k} - \log\left(1+\frac{1}{k}\right)\right]$$

$$= \sum_{k=1}^{\infty} \frac{(-1)^{k+1}}{k} - \sum_{k=1}^{\infty} (-1)^{k+1} \log\left(1+\frac{1}{k}\right)$$

$$= \log 2 - \log\frac{\pi}{2} = \log\frac{4}{\pi}$$

where, in the last line, we have employed (4.4.102) from Volume III. Sondow also notes that $\gamma$ and $\log\frac{4}{\pi}$ are related by Euler's formula

(E.6j) $$\gamma - \log\frac{4}{\pi} = 2\sum_{n=2}^{\infty} (-1)^n \frac{\varsigma(n)}{n2^n}$$

and this concurs with the Maclaurin expansion (E.22n) with $x = 1/2$.

Sondow's formula may also be obtained from Lerch's series expansion for the digamma function for $0 < x < 1$ (see for example [71a, p.105] and [104b, p.204])

(E.6k) $$\psi(x)\sin\pi x + \frac{\pi}{2}\cos\pi x + (\gamma + \log 2\pi)\sin\pi x = \sum_{n=1}^{\infty} \sin(2n+1)\pi x . \log\frac{n+1}{n}$$

Letting $x = 1/2$ we obtain

$$\psi(1/2) + \gamma + \log 2\pi = \sum_{n=1}^{\infty} (-1)^n \log\frac{n+1}{n}$$

and, since [126, p.20] $\psi(1/2) = -\gamma - 2\log 2$, (E.6i) follows automatically.

We may also write (E.6k) as

(E.6l) $$\psi(x) + \frac{\pi}{2}\cot\pi x + (\gamma + \log 2\pi) = \sum_{n=1}^{\infty} \frac{\sin(2n+1)\pi x}{\sin\pi x} . \log\frac{n+1}{n}$$

The above formula suggests that integration may be fruitful employing the integral in G&R [74, p.163]

$$\int \frac{\sin(2n+1)x}{\sin x}dx = x + \sum_{k=1}^{n} \frac{\sin 2kx}{k}$$

The integral $\int x^p \psi(x)$ may also produce interesting results.

Using (E.14) we may see that $\lim_{x\to 0}[\psi(x)\sin \pi x] = -\pi$ and hence, curiously, we see that (E.6k) is not true in the limit as $x \to 0$.

We saw in (E.5a) that

$$H_n = \sum_{k=1}^{n} u_k + \log(n+1)$$

and therefore we have

$$\sum_{n=1}^{\infty} \frac{H_n}{n^s} = \sum_{n=1}^{\infty} \frac{1}{n^s} \sum_{k=1}^{n} u_k + \sum_{n=1}^{\infty} \frac{\log(n+1)}{n^s}$$

We have

$$\sum_{n=1}^{\infty} \frac{1}{n^s} \sum_{k=1}^{n} u_k = \sum_{n=1}^{\infty} \frac{1}{n^s} \sum_{k=1}^{\infty} \int_0^1 \frac{x}{k(k+x)} dx = \int_0^1 \sum_{n=1}^{\infty} \frac{1}{n^s} \sum_{k=1}^{\infty} \frac{x}{k(k+x)} dx$$

$$= \int_0^1 \sum_{n=1}^{\infty} \frac{x}{n(n+x)} \sum_{k=n}^{\infty} \frac{1}{k^s} dx$$

$$= \int_0^1 \sum_{n=1}^{\infty} \frac{x}{n(n+x)} [\varsigma(s) - H_{n-1}^{(s)}] dx$$

$$= \varsigma(s) \int_0^1 \sum_{n=1}^{\infty} \frac{x}{n(n+x)} dx - \int_0^1 \sum_{n=1}^{\infty} \frac{x H_{n-1}^{(s)}}{n(n+x)} dx$$

$$= \gamma\varsigma(s) - \int_0^1 \sum_{n=1}^{\infty} \frac{x H_{n-1}^{(s)}}{n(n+x)} dx$$

$$= \gamma\varsigma(s) - \int_0^1 \sum_{n=1}^{\infty} \frac{x H_n^{(s)}}{n(n+x)} dx + \int_0^1 \sum_{n=1}^{\infty} \frac{x}{n^{s+1}(n+x)} dx$$

$$= \gamma\varsigma(s) - \sum_{n=1}^{\infty} \frac{H_n^{(s)}}{n} \log\frac{n+1}{n} + \sum_{n=1}^{\infty} \frac{1}{n^{s+1}} \log\frac{n+1}{n}$$

$$= \gamma\varsigma(s) - \sum_{n=1}^{\infty} \frac{H_n^{(s)}}{n} \log\frac{n+1}{n} + \sum_{n=1}^{\infty} \frac{\log(n+1)}{n^{s+1}} + \varsigma'(s+1)$$

therefore we obtain

$$\sum_{n=1}^{\infty} \frac{H_n}{n^s} = \gamma\varsigma(s) - \sum_{n=1}^{\infty} \frac{H_n^{(s)}}{n} \log\frac{n+1}{n} + \sum_{n=1}^{\infty} \frac{\log(n+1)}{n^{s+1}} + \varsigma'(s+1) + \sum_{n=1}^{\infty} \frac{\log(n+1)}{n^s}$$

Since $H_n - \log n + \varepsilon_n = \gamma$ we obtain

$$\sum_{n=1}^{\infty} \frac{H_n}{n^s} = \sum_{n=1}^{\infty} \frac{\log n}{n^s} + \sum_{n=1}^{\infty} \frac{\varepsilon_n}{n^s} + \gamma\varsigma(s)$$

and with $s = 3$ we get

$$\sum_{n=1}^{\infty} \frac{H_n}{n^3} = 2\varsigma(3) = -\varsigma'(3) + \sum_{n=1}^{\infty} \frac{\varepsilon_n}{n^s} + \gamma\varsigma(3)$$

In (4.17) we showed that

(E.7) $$H_n = \int_0^1 \frac{1-(1-x)^n}{x} dx$$

Using the substitution $x = t/n$ we have

$$H_n = \int_0^n \frac{1-(1-t/n)^n}{t} dt$$

$$= \int_0^1 \frac{1-(1-t/n)^n}{t} dt + \int_1^n \frac{1-(1-t/n)^n}{t} dt$$

$$= \int_0^1 \frac{1-(1-t/n)^n}{t} dt - \int_1^n \frac{(1-t/n)^n}{t} dt + \log n$$

Hence using the definition in (E.6) we have

$$\gamma = \lim_{n\to\infty}\int_0^1 \frac{1-(1-t/n)^n}{t}\,dt - \lim_{n\to\infty}\int_1^n \frac{(1-t/n)^n}{t}\,dt$$

Since $e^{-t} = \lim_{n\to\infty}\left(1-\frac{t}{n}\right)^n$ this suggests that

(E.8) $$\gamma = \int_0^1 \frac{1-e^{-t}}{t}\,dt - \int_1^\infty \frac{e^{-t}}{t}\,dt$$

and a rigorous proof of (E.8) is shown in [135, p.241]. This analysis can be justified by the Lebesgue dominated convergence theorem.

Making the substitution $x = e^{-u}$ in (E.7) we obtain

$$H_n = \int_0^\infty \left[1-(1-e^{-u})^n\right] du$$

and, for example, we have

$$\sum_{n=1}^\infty \frac{H_n}{n^2} = \sum_{n=1}^\infty \frac{1}{n^2}\int_0^\infty \left[1-(1-e^{-u})^n\right] du$$

$$= \sum_{n=1}^\infty \int_0^\infty e^{-ns}ds \int_0^\infty e^{-nt}dt \int_0^\infty 1-(1-e^{-u})^n\, du$$

$$= \int_0^\infty\int_0^\infty\int_0^\infty \sum_{n=1}^\infty e^{-n(s+t)}\left[1-(1-e^{-u})^n\right] ds\,dt\,du$$

$$= \int_0^\infty\int_0^\infty\int_0^\infty \left[\frac{1}{e^{s+t}-1} - \frac{1-e^{-u}}{e^{s+t}-(1-e^{-u})}\right] ds\,dt\,du$$

Let us now consider the integral

$$\int_0^\infty e^{-x}\log x\,dx = \int_0^1 e^{-x}\log x\,dx + \int_1^\infty e^{-x}\log x\,dx$$

$$= -\int_0^1 \frac{d}{dx}\left(e^{-x}-1\right)\log x\,dx - \int_1^\infty \frac{d}{dx}\left(e^{-x}-1\right)\log x\,dx$$

Integration by parts results in

$$= -\int_0^1 \frac{1-e^{-x}}{x} dx + \int_1^\infty \frac{e^{-x}}{x} dx$$

and using (E.8) we obtain

(E.9) $$\int_0^\infty e^{-x} \log x \, dx = -\gamma$$

Euler’s gamma function was defined in (4.3.1) as

(E.10) $$\Gamma(x) = \int_0^\infty t^{x-1} e^{-t} dt \quad , \text{Re}(x) > 0$$

Hence by parametric differentiation we have

(E.10a) $$\Gamma'(x) = \int_0^\infty t^{x-1} e^{-t} \log t \, dt$$

and therefore using (E.9) we have

(E.10b) $$\Gamma'(1) = \int_0^\infty e^{-t} \log t \, dt = -\gamma$$

Euler defined the gamma function as an infinite product and we can easily obtain this representation from Euler’s definition (E.10)

$$\Gamma(x) = \int_0^\infty e^{-t} t^{x-1} dt$$

$$= \int_0^\infty \lim_{n\to\infty} \left(1 - \frac{t}{n}\right)^n t^{x-1} dt$$

$$= \lim_{n\to\infty} \int_0^n \left(1 - \frac{t}{n}\right)^n t^{x-1} dt$$

With the substitution $u = t/n$ we have

$$= \lim_{n\to\infty} n^x \int_0^1 (1-u)^n u^{x-1} du$$

This integral was computed in (4.4.1) and we therefore have

(E.11) $$\Gamma(x) = \lim_{n\to\infty} n^x \frac{n!}{x(1+x)...(n+x)} = \lim_{n\to\infty} n^x \prod_{k=1}^{n} \frac{1}{1+x/k}$$

$$= \frac{1}{x} \lim_{n\to\infty} \frac{1.2...n}{(1+x)...(n+x)} n^x$$

$$= \frac{1}{x} \lim_{n\to\infty} \frac{1}{(1+x)(1+x/2)...(1+x/n)} n^x$$

$$= \frac{1}{x} \lim_{n\to\infty} \frac{1}{(1+x)(1+x/2)...(1+x/n)} \frac{2^x 3^x ... n^x}{1^x 2^x ...(n-1)^x}$$

and, since $\lim_{n\to\infty} \frac{(n+1)^x}{n^x} = 1$, we can multiply by that factor to obtain

$$= \frac{1}{x} \lim_{n\to\infty} \frac{1}{(1+x)(1+x/2)...(1+x/n)} \frac{2^x 3^x ...(n+1)^x}{1^x 2^x ... n^x}$$

$$= \frac{1}{x} \prod_{n=1}^{\infty} \left[ \frac{1}{1+x/n} \frac{(n+1)^x}{n^x} \right]$$

Therefore we have Euler's formula for the gamma function

(E.12) $$\Gamma(x) = \frac{1}{x} \prod_{n=1}^{\infty} \left[ \left(1+\frac{1}{n}\right)^x \left(1+\frac{x}{n}\right)^{-1} \right]$$

It is quite obvious from this formulation that $\Gamma(1) = 1$.

A slightly more concise derivation of (E.12) may be obtained by first of all taking the logarithm of (E.11).We have

$$\log\Gamma(x) = \log \lim_{n\to\infty} n^x \frac{n!}{x(1+x)...(n+x)}$$

$$= \lim_{n\to\infty} \log\, n^x \frac{n!}{x(1+x)...(n+x)}$$

$$= \lim_{n\to\infty} \left[ x\log n + \log n! - \sum_{k=0}^{n} \log(k+x) \right]$$

$$= \lim_{n\to\infty} \left[ x\log n - \log x - \left\{ \sum_{k=1}^{n} \log(k+x) - \sum_{k=1}^{n} \log k \right\} \right]$$

$$= \lim_{n\to\infty}\left[x\log n - \log x - \sum_{k=1}^{n}\log\left(1+\frac{x}{k}\right)\right]$$

Therefore we obtain an equivalent form of (E.12)

(E.12a) $$\log\Gamma(x+1) = \sum_{n=1}^{\infty}\left[x\log\left(1+\frac{1}{n}\right) - \log\left(1+\frac{x}{n}\right)\right]$$

From this formulation it is immediately seen that $\log\Gamma(2) = 0$ and hence $\Gamma(2) = 1$.

We see from (E.12a) that

$$\log\Gamma(x+1) = \sum_{n=1}^{\infty}\left[x\log\left(1+\frac{1}{n}\right) + \sum_{k=1}^{\infty}\frac{(-1)^k x^k}{kn^k}\right]$$

$$= \sum_{n=1}^{\infty}\left[x\left[\log\left(1+\frac{1}{n}\right) - \frac{1}{n}\right] + \sum_{k=2}^{\infty}\frac{(-1)^k x^k}{kn^k}\right]$$

$$= x\sum_{n=1}^{\infty}\left[\log\left(1+\frac{1}{n}\right) - \frac{1}{n}\right] + \sum_{n=1}^{\infty}\sum_{k=2}^{\infty}\frac{(-1)^k x^k}{kn^k}$$

$$= -\gamma x + \sum_{k=2}^{\infty}\frac{(-1)^k \varsigma(k)}{k}x^k$$

Differentiating (E.12a) we obtain

(E.12ai) $$\psi(x+1) = \sum_{n=1}^{\infty}\left[\log\left(1+\frac{1}{n}\right) - \frac{1}{n+x}\right]$$

We may write this as

$$\psi(x+1) = \sum_{n=1}^{\infty}\left[\log\left(1+\frac{1}{n}\right) - \frac{1}{n} - \left(\frac{1}{n+x} - \frac{1}{n}\right)\right]$$

and using (E.14) this becomes

$$\psi(x+1) = \sum_{n=1}^{\infty}\left[\log\left(1+\frac{1}{n}\right) - \frac{1}{n}\right] + \psi(x) + \frac{1}{x} + \gamma$$

Since $\psi(x+1) = \psi(x) + \dfrac{1}{x}$ we therefore have in accordance with (E.6h) (see also the reference to the generalised Euler constant function in (E.43b))

$$\gamma = \sum_{n=1}^{\infty} \left[ \frac{1}{n} - \log\left(1+\frac{1}{n}\right) \right]$$

Integrating (E.12ai) gives us

$$\int_0^u x\psi(x+1)dx = \sum_{n=1}^{\infty} \left[ \frac{1}{2}u^2 \log\left(1+\frac{1}{n}\right) - u + n\log\left(1+\frac{u}{n}\right) \right]$$

and we recall from (4.3.83) in Volume II(a) that

$$\phi_2(u) = \int_0^u x\psi(x+1)dx = \log G(u+1) - \frac{1}{2}u\log(2\pi) + \frac{1}{2}u(u+1)$$

Hence we have a series expansion for the Barnes double gamma function

(E.12aii)

$$\log G(u+1) - \frac{1}{2}u\log(2\pi) + \frac{1}{2}u(u+1) = \sum_{n=1}^{\infty} \left[ \frac{1}{2}u^2 \log\left(1+\frac{1}{n}\right) - u + n\log\left(1+\frac{u}{n}\right) \right]$$

and expanding $\log\left(1+\frac{u}{n}\right)$ as before will result in (4.3.87e) of Volume II(a).

We may write (E.12aii) as

$$G(u+1) = (2\pi)^{\frac{1}{2}u} \exp\left[-\frac{1}{2}u(u+1)\right] \prod_{n=1}^{\infty} \left[ \left(1+\frac{u}{n}\right)^n \left(1+\frac{1}{n}\right)^{\frac{1}{2}u^2} \exp(-u) \right]$$

and compare this with (4.3.81a) from Volume II(a)

$$G(u+1) = (2\pi)^{\frac{1}{2}u} \exp\left[-\frac{1}{2}u(u+1)\right] \exp\left[-\frac{1}{2}\gamma u^2\right] \prod_{n=1}^{\infty} \left[ \left(1+\frac{u}{n}\right)^n \exp\left(\frac{u^2}{2n} - u\right) \right]$$

This implies that

$$\exp\left[-\frac{1}{2}\gamma u^2\right] \prod_{n=1}^{\infty} \left[ \exp\left(\frac{u^2}{2n}\right) \right] = \prod_{n=1}^{\infty} \left[ \left(1+\frac{1}{n}\right)^{\frac{1}{2}u^2} \right]$$

With $u = 1$ in (E.12aii) we obtain

$$1 - \frac{1}{2}\log(2\pi) = \sum_{n=1}^{\infty} \left[ \frac{1}{2}\log\left(1+\frac{1}{n}\right) - 1 + n\log\left(1+\frac{1}{n}\right) \right]$$

and hence we have

(E.12aiii) $$1-\frac{1}{2}\log(2\pi)=\sum_{n=1}^{\infty}\log\frac{1}{e}\left(1+\frac{1}{n}\right)^{n+\frac{1}{2}}$$

In 2005, Müller and Schleicher [103ad] using their method of fractional sums showed that for $a>0$

$$\prod_{n=1}^{\infty}\frac{1}{e}\left(1+\frac{1}{an}\right)^{an+\frac{1}{2}}=\sqrt{\frac{\Gamma\left(1+\frac{1}{a}\right)}{2\pi}}\exp\left[\frac{1}{2}\left(1+\frac{1}{a}\right)-a\left\{\varsigma'\left(-1,1+\frac{1}{a}\right)-\varsigma'(-1)\right\}\right]$$

and with $a=1$ this becomes

$$\prod_{n=1}^{\infty}\frac{1}{e}\left(1+\frac{1}{n}\right)^{n+\frac{1}{2}}=\sqrt{\frac{1}{2\pi}}\exp\left[1-\left\{\varsigma'(-1,2)-\varsigma'(-1)\right\}\right]$$

which is equivalent to (E.12aiii) since $\varsigma'(-1,2)=\varsigma'(-1)$.

We may write (E.12aii) as

$$\log G(u+1)-\frac{1}{2}u\log(2\pi)+\frac{1}{2}u(u+1)=\sum_{n=1}^{\infty}\left[\frac{1}{2}u^2\left\{\log\left(1+\frac{1}{n}\right)-\frac{1}{n}\right\}+\frac{1}{2n}u^2-u+n\log\left(1+\frac{u}{n}\right)\right]$$

and we then obtain

$$\log G(u+1)=\frac{1}{2}u\log(2\pi)-\frac{1}{2}u(u+1)-\frac{1}{2}\gamma u^2=\sum_{n=1}^{\infty}\left[\frac{1}{2n}u^2-u+n\log\left(1+\frac{u}{n}\right)\right]$$

Integrating (E.12aii) results in

$$\int_0^x\log G(1+u)du-\frac{1}{4}x^2\log(2\pi)+\frac{1}{6}x^3+\frac{1}{4}x^2=$$

$$\sum_{n=1}^{\infty}\left[\frac{1}{6}x^3\log\left(1+\frac{1}{n}\right)-\frac{1}{2}x^2-nx+n(n+x)\log\left(1+\frac{x}{n}\right)\right]$$

and we recall from (4.3.87c) that

$$\int_0^x\log G(1+u)du=\left(\frac{1}{4}-2\log A\right)x+\frac{1}{4}x^2\log(2\pi)-\frac{1}{6}x^3+(x-1)\log G(1+x)-2\log\Gamma_3(1+x)$$

Hence we obtain

$$\sum_{n=1}^{\infty}\left[\frac{1}{6}x^3\log\left(1+\frac{1}{n}\right)-\frac{1}{2}x^2-nx+n(n+x)\log\left(1+\frac{x}{n}\right)\right]$$

(E.12aiv)

$$=\left(\frac{1}{4}-2\log A\right)x+\frac{1}{4}x^2+(x-1)\log G(1+x)-2\log\Gamma_3(1+x)$$

and with $x=1$ we get

(E.12av) $$\sum_{n=1}^{\infty}\left[[n(n+1)+1]\log\left(1+\frac{1}{n}\right)-\frac{1}{2}-n\right]=\frac{1}{2}-2\log A$$

(which we may compare with (F.24g)). I subsequently discovered that (E.12av) had been previously reported by Ferreira and López [65b] in 2001.

There are several other expressions for the gamma function and these include the so-called Weierstrass canonical form obtained by Weierstrass in 1856 (this formula was actually first derived by F.W. Newman in 1848: see [135, p.236])

From (E.11) we have

(E.12b) $$\Gamma(x)=\lim_{n\to\infty}n^x\frac{n!}{x(1+x)...(n+x)}$$

$$=\lim_{n\to\infty}\left[x^{-1}\left(1+\frac{x}{1}\right)^{-1}\left(1+\frac{x}{2}\right)^{-1}...\left(1+\frac{x}{n}\right)^{-1}n^x\right]$$

Therefore we obtain

$$\frac{1}{\Gamma(x)}=\lim_{n\to\infty}\left[x\left(1+\frac{x}{1}\right)\left(1+\frac{x}{2}\right)...\left(1+\frac{x}{n}\right)e^{-x\log n}\right]$$

$$\frac{1}{\Gamma(x)}=\lim_{n\to\infty}\left[x\left(1+\frac{x}{1}\right)e^{-x}\left(1+\frac{x}{2}\right)e^{-x/2}...\left(1+\frac{x}{n}\right)e^{-x/n}\exp\left(\left[1+\frac{1}{2}+...+\frac{1}{n}-\log n\right]x\right)\right]$$

Hence we have the Weierstrass canonical form

(E.13) $$\frac{1}{\Gamma(z)}=ze^{\gamma z}\prod_{n=1}^{\infty}\left\{\left(1+\frac{z}{n}\right)e^{-\frac{z}{n}}\right\}$$

Taking logarithms results in

(E.13a) $$\log\Gamma(z) = -\log z - \gamma z - \sum_{n=1}^{\infty}\left[\log\left(1+\frac{z}{n}\right)-\frac{z}{n}\right]$$

By logarithmically differentiating (E.13) it is easily seen that $\Gamma'(1) = -\gamma$ and the following identities are also easily derived.

(E.14) $$\frac{\Gamma'(x)}{\Gamma(x)} = -\gamma - \sum_{k=1}^{\infty}\left(\frac{1}{x+k-1}-\frac{1}{k}\right) = -\gamma - \frac{1}{x} + \sum_{k=1}^{\infty}\left(\frac{1}{k}-\frac{1}{x+k}\right)$$

(E.15) $$\frac{\Gamma'(x)}{\Gamma(x)} - \frac{\Gamma'(1)}{\Gamma(1)} = -\sum_{k=0}^{\infty}\left(\frac{1}{x+k}-\frac{1}{1+k}\right) = (x-1)\sum_{k=0}^{\infty}\frac{1}{(x+k)(1+k)}$$

(E.16) $$\psi'(x) = \frac{d^2}{dx^2}\log\Gamma(x) = \sum_{k=0}^{\infty}\frac{1}{(x+k)^2}$$

(E.16a) $$\psi^{(n)}(x) = \frac{d^{n+1}}{dx^{n+1}}\log\Gamma(x) = (-1)^{n+1}n!\sum_{k=0}^{\infty}\frac{1}{(x+k)^{n+1}}$$

$$= (-1)^{n+1}n!\varsigma(n+1,x)$$

(E.16b) $$\psi^{(n)}(1) = (-1)^{n+1}n!\varsigma(n+1)$$

Rather formally, we may retrace our steps as follows. From (E.6) and (E.14) we have

$$\gamma = \sum_{k=1}^{\infty}\int_0^1 \frac{x}{k(k+x)}dx = \sum_{k=1}^{\infty}\int_0^1\left(\frac{1}{k}-\frac{1}{k+x}\right)dx$$

$$= \int_0^1\left[\psi(x)+\frac{1}{x}+\gamma\right]dx = \gamma$$

From (E.11) we see that

$$\Gamma(x+1) = \lim_{n\to\infty} n^x \frac{n!}{(1+x)...(n+x)}\frac{n}{(n+x+1)}$$

$$= \lim_{n\to\infty} n^x \frac{n!}{(1+x)...(n+x)}$$

$$\log\Gamma(x+1) = \lim_{n\to\infty}\{x\log n + \log n! - [\log(x+1)+\log(x+1)+...+\log(x+n)]\}$$

and differentiation results in

$$\psi(x+1) = \lim_{n\to\infty}\left\{\log n - \left[\frac{1}{x+1}+\frac{1}{x+2}+...+\frac{1}{x+n}\right]\right\}$$

With $x=0$ we see that

$$\psi(1) = \lim_{n\to\infty}\left\{\log n - \left[\frac{1}{1}+\frac{1}{2}+...+\frac{1}{n}\right]\right\} = -\gamma$$

From (E.14) we may note that $\lim_{x\to 0}[x\psi(x)] = -1$.

Differentiating (E.14) we obtain

(E.16c) $$\psi'(x) = \frac{\Gamma(x)\Gamma''(x)-\left[\Gamma'(x)\right]^2}{\Gamma^2(x)} = \sum_{k=0}^{\infty}\frac{1}{(k+x)^2} = \varsigma(2,x)$$

and hence we get

(E.16d) $$\Gamma''(1) = \gamma^2 + \varsigma(2)$$

Differentiating again and solving for $\Gamma^{(3)}(1)$ we have

(E.16e) $$-\Gamma^{(3)}(1) = \gamma^3 + \frac{\gamma\pi^2}{2} + 2\varsigma(3)$$

Using (E.11) and Euler's reflection formula one can easily provide a proof of Euler's infinite product for $\sin x$ (1.6d).

$$\sin x = x\prod_{n=1}^{\infty}\left(1-\frac{x^2}{n^2\pi^2}\right)$$

We have

$$\Gamma(x)\Gamma(1-x) = \left[\lim_{n\to\infty}\frac{n!}{x(1+x)...(n+x)}\right]\left[\lim_{n\to\infty}\frac{n^{1-x}n!}{(1-x)(2-x)...(n+1-x)}\right]$$

$$= \lim_{n\to\infty}\frac{n\left[n!\right]^2}{x(1^2-x^2)(2^2-x^2)...(n^2-x^2)(n+1-x)}$$

$$= \lim_{n\to\infty}\left[x(1^2-x^2/1^2)(2^2-x^2/2^2)...(n^2-x^2/n^2)(n+1-x)/n\right]^{-1}$$

$$= \left[x\prod_{n=1}^{\infty}\left(1-\frac{x^2}{n^2}\right)\right]^{-1}$$

and using Euler's reflection formula (C.1)

$$\Gamma(x)\Gamma(1-x) = \frac{\pi}{\sin \pi x}$$

we deduce that

$$\sin \pi x = \pi x \prod_{n=1}^{\infty} \left( 1 - \frac{x^2}{n^2} \right)$$

Using (4.4.1) and (E.11) we have yet another representation for $\Gamma(x)$

(E.17) $$\Gamma(x) = \lim_{n \to \infty} n^x \sum_{k=0}^{n} \binom{n}{k} \frac{(-1)^k}{k+x}$$

This formula does not appear to be widely reported in the literature. The singularities of $\Gamma(x)$ when $x$ is a negative integer are clearly identifiable from this formula.

Upon differentiation of (E.17) we formally obtain

(E.18) $$\Gamma'(x) = \lim_{n \to \infty} n^x \left( \log n \sum_{k=0}^{n} \binom{n}{k} \frac{(-1)^k}{k+x} - \sum_{k=0}^{n} \binom{n}{k} \frac{(-1)^k}{(k+x)^2} \right)$$

In (4.2.3), (4.2.16) and (4.2.28) we showed that

(E.18a) $$\sum_{k=0}^{n} \binom{n}{k} \frac{(-1)^k}{k+1} = \frac{1}{(n+1)}$$

(E.18b) $$\sum_{k=0}^{n} \binom{n}{k} \frac{(-1)^k}{(1+k)^2} = \frac{H_{n+1}}{n+1}$$

(E.18c) $$\sum_{k=0}^{n} \binom{n}{k} \frac{(-1)^k}{(1+k)^3} = \frac{1}{n+1} \left\{ \frac{1}{2} (H_{n+1}^{(1)})^2 + \frac{1}{2} H_{n+1}^{(2)} \right\}$$

and, hence from (E.17) and (E.18a), we again see that

(E.18d) $$\Gamma(1) = \lim_{n \to \infty} n \frac{1}{(n+1)} = 1$$

With (E.18) and (E.18b) we obtain another proof of the well-known result $\Gamma'(1) = -\gamma$

$$\Gamma'(1) = \lim_{n \to \infty} n \left( \frac{\log n}{n+1} - \frac{H_{n+1}}{n+1} \right)$$

$$= \lim_{n\to\infty} \frac{n}{n+1}(\log n - H_{n+1})$$

(E.18e) $$= \lim_{n\to\infty} \left( \frac{n}{n+1}(\log n - H_n) - \frac{n}{(n+1)^2} \right) = -\gamma$$

We refer again to (E.18)

$$\Gamma'(x) = \lim_{n\to\infty} n^x \left( \log n \sum_{k=0}^{n} \binom{n}{k} \frac{(-1)^k}{k+x} - \sum_{k=0}^{n} \binom{n}{k} \frac{(-1)^k}{(k+x)^2} \right)$$

and using (4.4.73c) from Volume III

$$\sum_{k=0}^{n} \binom{n}{k} \frac{(-1)^k}{(k+x)^2} = -B(x, n+1)\left[\psi(x) - \psi(x+n+1)\right]$$

we may write (E.18) as

$$\Gamma'(x) = \lim_{n\to\infty} n^x B(x, n+1)\left[\log n + \psi(x) - \psi(x+n+1)\right]$$

$$= \lim_{n\to\infty} n^x B(x, n+1) \lim_{n\to\infty} \left[\log n + \psi(x) - \psi(x+n+1)\right]$$

and using (4.4.5) this becomes

$$= \Gamma(x) \lim_{n\to\infty} \left[\log n + \psi(x) - \psi(x+n+1)\right]$$

$$= \Gamma(x)\psi(x) + \Gamma(x) \lim_{n\to\infty} \left[\log n - \psi(x+n+1)\right]$$

and this proves that $\lim_{n\to\infty} \left[\log n - \psi(x+n+1)\right] = 0$.

Differentiating (E.18) we obtain

(E.19) $$\Gamma''(x) = \lim_{n\to\infty} n^x \left( -\log n \sum_{k=0}^{n} \binom{n}{k} \frac{(-1)^k}{(k+x)^2} + 2\sum_{k=0}^{n} \binom{n}{k} \frac{(-1)^k}{(k+x)^3} \right)$$

$$+ \lim_{n\to\infty} n^x \log n \left( \log n \sum_{k=0}^{n} \binom{n}{k} \frac{(-1)^k}{k+x} - \sum_{k=0}^{n} \binom{n}{k} \frac{(-1)^k}{(k+x)^2} \right)$$

and hence

$$\Gamma''(x) = \lim_{n\to\infty} n^x \left( \log^2 n \sum_{k=0}^{n} \binom{n}{k} \frac{(-1)^k}{(k+x)} - 2\log n \sum_{k=0}^{n} \binom{n}{k} \frac{(-1)^k}{(k+x)^2} + 2\sum_{k=0}^{n} \binom{n}{k} \frac{(-1)^k}{(k+x)^3} \right)$$

Hence we deduce

$$\Gamma''(1) = \lim_{n\to\infty} n\left(\frac{\log^2 n}{n+1} - 2\log n\frac{H_{n+1}^{(1)}}{n+1} + \frac{\left(H_{n+1}^{(1)}\right)^2 + H_{n+1}^{(2)}}{n+1}\right)$$

$$= \lim_{n\to\infty} \frac{n}{n+1}\left(\left[\log n - H_{n+1}^{(1)}\right]^2 + H_{n+1}^{(2)}\right)$$

$$= \lim_{n\to\infty} \frac{n}{n+1}\left(\left[\log n - H_n^{(1)}\right]^2 + H_n^{(2)} - 2\frac{(\log n - H_n^{(1)})}{n} + \frac{2}{n^2}\right)$$

This gives us the well-known result

(E.19a) $\Gamma''(1) = \gamma^2 + \varsigma(2)$

which is the expression obtained by Levenson [99] using the integral form of the definition of the $\Gamma(x)$ function. It is clear that this method may be extended, with increasing algebraic effort, to derive expressions for $\Gamma^{(k)}(x)$.

Using (4.4.73d) from Volume III

$$2\sum_{k=0}^{n}\binom{n}{k}\frac{(-1)^k}{(k+x)^3} = B(x,n+1)\left[\psi'(x) - \psi'(x+n+1) + \left\{\psi(x) - \psi(x+n+1)\right\}^2\right]$$

we may write (E.19) as

$$\Gamma''(x) = \lim_{n\to\infty} n^x B(x,n+1)\left(\begin{array}{l}\log^2 n + 2\log n\left[\psi(x) - \psi(x+n+1)\right] \\ \\ +\left[\psi'(x) - \psi'(x+n+1) + \left\{\psi(x) - \psi(x+n+1)\right\}^2\right]\end{array}\right)$$

$$= \Gamma(x)\lim_{n\to\infty}\left(\begin{array}{l}\log^2 n + 2\log n\left[\psi(x) - \psi(x+n+1)\right] \\ \\ +\left[\psi'(x) - \psi'(x+n+1) + \left\{\psi(x) - \psi(x+n+1)\right\}^2\right]\end{array}\right)$$

We have

$$\Gamma''(x) = \Gamma(x)\left[\psi'(x) + \psi^2(x)\right]$$

and hence we have

$$\psi'(x)+\psi^2(x)=\lim_{n\to\infty}\begin{pmatrix}\log^2 n+2\log n\left[\psi(x)-\psi(x+n+1)\right] \\ \\ +\left[\psi'(x)-\psi'(x+n+1)+\left\{\psi(x)-\psi(x+n+1)\right\}^2\right]\end{pmatrix}$$

We note that

$$\log^2 n+2\log n\left[\psi(x)-\psi(x+n+1)\right]+\left[\psi'(x)-\psi'(x+n+1)+\left\{\psi(x)-\psi(x+n+1)\right\}^2\right]=$$

$$\psi'(x)+\psi^2(x)+\log^2 n+2\log n\left[\psi(x)-\psi(x+n+1)\right]$$

$$-\psi'(x+n+1)-2\psi(x)\psi(x+n+1)+\psi^2(x+n+1)$$

and we therefore deduce that

$$\lim_{n\to\infty}\begin{pmatrix}\log^2 n+2\log n\left[\psi(x)-\psi(x+n+1)\right]-\psi'(x+n+1) \\ \\ -2\psi(x)\psi(x+n+1)+\psi^2(x+n+1)\end{pmatrix}=0$$

This may be written as

$$\lim_{n\to\infty}\left(\left[\log n-\psi(x+n+1)\right]^2+2\log n\,\psi(x)-\psi'(x+n+1)-2\psi(x)\psi(x+n+1)\right)=0$$

which results in

$$\lim_{n\to\infty}\left(2\log n\,\psi(x)-\psi'(x+n+1)-2\psi(x)\psi(x+n+1)\right)=0$$

$$\lim_{n\to\infty}\left(2\psi(x)[\log n-\psi(x+n+1)]-\psi'(x+n+1)\right)=0$$

This is obviously correct since $\lim_{n\to\infty}\psi'(x+n+1)=0$.

We now recall the results obtained by Larcombe et al. [95] (as reported in equations (4.4.135) et seq. of Volume IV of this series of papers) where they show that for integers $m\ge 1, n\ge 0$

$$m\binom{m+n}{n}\sum_{k=0}^{n}\binom{n}{k}\frac{(-1)^k}{m+k}=1$$

$$m\binom{m+n}{n}\sum_{k=0}^{n}\binom{n}{k}\frac{(-1)^k}{(m+k)^2}=\sum_{k=m}^{m+n}\frac{1}{k}$$

$$2m\binom{m+n}{n}\sum_{k=0}^{n}\binom{n}{k}\frac{(-1)^k}{(m+k)^3}=\left(\sum_{k=m}^{m+n}\frac{1}{k}\right)^2+\sum_{k=m}^{m+n}\frac{1}{k^2}$$

Therefore, using (E.17) we have

$$\Gamma(m)=\lim_{n\to\infty}n^m\sum_{k=0}^{n}\binom{n}{k}\frac{(-1)^k}{k+m}=\lim_{n\to\infty}\frac{n^m}{m\binom{m+n}{n}}$$

and this accords with (E.12b).

Using (E.18) and $m$ equal to a positive integer we have

$$\Gamma'(m)=\lim_{n\to\infty}n^m\left(\log n\sum_{k=0}^{n}\binom{n}{k}\frac{(-1)^k}{k+m}-\sum_{k=0}^{n}\binom{n}{k}\frac{(-1)^k}{(k+m)^2}\right)$$

$$=\lim_{n\to\infty}\frac{n^m}{m\binom{m+n}{n}}\left(\log n-H_{m+n}^{(1)}+H_{m-1}^{(1)}\right)$$

$$=\Gamma(m)\lim_{n\to\infty}\left(\left[\log n-H_n^{(1)}\right]+H_{m+n}^{(1)}-H_{m-1}^{(1)}+H_n^{(1)}\right)$$

Therefore we have (a subset of (E.14))

(E.20) $$\frac{\Gamma'(m)}{\Gamma(m)}=-\gamma-\sum_{k=0}^{\infty}\left(\frac{1}{m+k}-\frac{1}{k+1}\right)$$

because $H_{m+n}^{(1)}=H_{m-1}^{(1)}+\frac{1}{m}+\frac{1}{m+1}...+\frac{1}{m+n}$.

This gives us for example

(E.20a) $$\psi(1)=\frac{\Gamma'(1)}{\Gamma(1)}=-\gamma \qquad \psi(2)=-\gamma+1$$

and more generally we have [126, p.23]

$$\psi(m)=\frac{\Gamma'(m)}{\Gamma(m)}=-\gamma-\sum_{k=0}^{\infty}\left(\frac{1}{m+k}-\frac{1}{k+1}\right)=H_{m-1}^{(1)}-\gamma$$

With (E.19) and $m$ equal to a positive integer we obtain

$$\Gamma''(m)=\lim_{n\to\infty}n^m\left(\log^2 n\sum_{k=0}^{n}\binom{n}{k}\frac{(-1)^k}{(k+m)}-2\log n\sum_{k=0}^{n}\binom{n}{k}\frac{(-1)^k}{(k+m)^2}+2\sum_{k=0}^{n}\binom{n}{k}\frac{(-1)^k}{(k+m)^3}\right)$$

$$= \lim_{n\to\infty} \frac{n^m}{m\binom{m+n}{n}} \left( \log^2 n - 2\left[ H_{m+n}^{(1)} - H_{m-1}^{(1)} \right] \log n + \left[ H_{m+n}^{(1)} - H_{m-1}^{(1)} \right]^2 + \left[ H_{m+n}^{(2)} - H_{m-1}^{(2)} \right] \right)$$

$$= \Gamma(m) \lim_{n\to\infty} \left( \log^2 n - 2\left[ H_{m+n}^{(1)} - H_{m-1}^{(1)} \right] \log n + \left[ H_{m+n}^{(1)} - H_{m-1}^{(1)} \right]^2 + \left[ H_{m+n}^{(2)} - H_{m-1}^{(2)} \right] \right)$$

$$= \Gamma(m) \lim_{n\to\infty} \left( \left\{ \log n - \left[ H_{m+n}^{(1)} - H_{m-1}^{(1)} \right] \right\}^2 + \left[ H_{m+n}^{(2)} - H_{m-1}^{(2)} \right] \right)$$

Therefore we have in terms of the Hurwitz zeta function

(E.21) $$\frac{\Gamma''(m)}{\Gamma(m)} = \varsigma(2,m) + \lim_{n\to\infty} \left\{ \log n - \sum_{k=0}^{n} \frac{1}{m+k} \right\}^2$$

We know from [126, p.14] that

$$\psi(x) = \frac{\Gamma'(x)}{\Gamma(x)} = \lim_{n\to\infty} \left\{ \log n - \sum_{k=0}^{n} \frac{1}{x+k} \right\}$$

and hence

(E.22) $$\frac{\Gamma''(m)}{\Gamma(m)} = \varsigma(2,m) + \psi^2(m)$$

This is just a particular case of the well-known formula (E16c) (also see for example [119])

(E.22a) $$\psi'(x) = \frac{\Gamma(x)\Gamma''(x) - \left(\Gamma'(x)\right)^2}{\Gamma^2(x)} = \varsigma(2,x)$$

We have $1/r = \int_0^\infty e^{-rx} dx$ ($r > 0$) and, integrating that expression, we obtain Frullani's integral

$$\int_1^n \frac{dr}{r} = \int_1^n dr \int_0^\infty e^{-rx} dx = \int_0^\infty dx \int_1^n e^{-rx} dr$$

which implies that

(E.22aa) $$\log n = \int_0^\infty \frac{e^{-x} - e^{-nx}}{x} dx$$

Therefore we have

$$a_n = H_n - \log n = \int_0^\infty \sum_{r=1}^n e^{-rx} dx - \int_0^\infty \frac{e^{-x} - e^{-nx}}{x} dx$$

$$= \int_0^\infty \frac{e^{-x} - e^{-(n+1)x}}{1 - e^{-x}} dx - \int_0^\infty \frac{e^{-x} - e^{-nx}}{x} dx$$

$$= \int_0^\infty \left[ \frac{1}{e^x - 1} - \frac{e^{-x}}{x} \right] dx + \int_0^\infty \left[ \frac{e^{-nx}}{x} - \frac{e^{-nx}}{e^x - 1} \right] dx$$

Now let $t = nx$ in the second integral to obtain

$$a_n = \int_0^\infty \left[ \frac{1}{e^x - 1} - \frac{e^{-x}}{x} \right] dx + \int_0^\infty \left[ \frac{e^{-x}}{x} - \frac{e^{-x}}{n(e^{x/n} - 1)} \right] dx$$

$$= \int_0^\infty \left[ \frac{1}{e^x - 1} - \frac{e^{-x}}{n(e^{x/n} - 1)} \right] dx$$

We recall the familiar limit $\lim_{n\to\infty} n\left[ A^{1/n} - 1 \right] = \log A$ and hence (formally) we are lead to the well-known integral expression for Euler's constant

(E.22b) $$\gamma = \lim_{n\to\infty} a_n = \int_0^\infty \left[ \frac{1}{e^x - 1} - \frac{e^{-x}}{x} \right] dx = \int_0^1 \left[ \frac{1}{1-y} + \frac{1}{\log y} \right] dy$$

A more rigorous proof is given in [110a].

We have from [110a]

$$H_n - \log n - \gamma = \int_0^\infty e^{-nx} \left[ \frac{1}{x} - \frac{1}{e^x - 1} \right] dx$$

and on summation we obtain

$$\sum_{n=1}^\infty \frac{H_n}{n^s} + \varsigma'(s) - \gamma\varsigma(s) = \sum_{n=1}^\infty \frac{1}{n^s} \int_0^\infty e^{-nx} \left[ \frac{1}{x} - \frac{1}{e^x - 1} \right] dx$$

In particular we have

$$\sum_{n=1}^\infty \frac{1}{n^2} \int_0^\infty e^{-nx} \left[ \frac{1}{x} - \frac{1}{e^x - 1} \right] dx = \sum_{n=1}^\infty \int_0^\infty e^{-nu} du \int_0^\infty e^{-nv} dv \int_0^\infty e^{-nx} \left[ \frac{1}{x} - \frac{1}{e^x - 1} \right] dx$$

$$= \sum_{n=1}^{\infty} \int_0^{\infty} du \int_0^{\infty} dv \int_0^{\infty} e^{-n(x+u+v)} \left[ \frac{1}{x} - \frac{1}{e^x - 1} \right] dx$$

$$= \int_0^{\infty} du \int_0^{\infty} dv \int_0^{\infty} \frac{1}{e^{x+u+v} - 1} \left[ \frac{1}{x} - \frac{1}{e^x - 1} \right] dx$$

We have

$$\int_0^{\infty} \frac{1}{e^{x+u+v} - 1} du = \int_0^{\infty} \frac{e^{-(x+u+v)}}{1 - e^{-(x+u+v)}} du = \log[1 - e^{-(x+u+v)}]\Big|_0^{\infty} = -\log[1 - e^{-(x+v)}]$$

and the Wolfram Integrator tells us that

$$\int \log[1 - e^{-(x+v)}] dv = \frac{1}{2} v^2 + v \log \frac{[1 - e^{-(x+v)}]}{[1 - e^{(x+v)}]} - Li_2[e^{(x+v)}]$$

With the substitution $y = e^{-v}$ we can easily find that

$$\int_0^{\infty} \log[1 - e^{-(x+v)}] dv = -Li_2[e^{-x}]$$

and hence we obtain

$$\sum_{n=1}^{\infty} \frac{H_n}{n^2} + \varsigma'(2) - \gamma\varsigma(2) = \int_0^{\infty} Li_2[e^{-x}] \left[ \frac{1}{x} - \frac{1}{e^x - 1} \right] dx$$

or equivalently

(E.22bi) $$2\varsigma(3) + \varsigma'(2) - \gamma\varsigma(2) = \int_0^{\infty} Li_2[e^{-x}] \left[ \frac{1}{x} - \frac{1}{e^x - 1} \right] dx$$

With the substitution $t = e^{-x}$ this may be written as

(E.22bii) $$2\varsigma(3) + \varsigma'(2) - \gamma\varsigma(2) = \int_0^{1} Li_2(t) \left[ \frac{1}{t-1} - \frac{1}{t \log t} \right] dt$$

We see from (E.22aa) that

$$\sum_{n=1}^{\infty} \frac{\log n}{n^s} = \sum_{n=1}^{\infty} \frac{1}{n^s} \int_0^{\infty} \frac{e^{-x} - e^{-nx}}{x} dx$$

which may be written as

$$-\varsigma'(s)=\int_0^\infty \frac{\varsigma(s)e^{-x}-Li_s[e^{-x}]}{x}dx$$

With $t=e^{-x}$ this becomes

(E.22biii) $$\varsigma'(s)=\int_0^1 \frac{\varsigma(s)t-Li_s(t)}{t\log t}dt$$

Does this hold for $s=0$ ? If so, we would have using (F.6)

$$\varsigma'(0)=-\frac{1}{2}\log(2\pi)=-\int_0^1 \frac{t/2+Li_0(t)}{t\log t}dt$$

$$=-\int_0^1 \frac{3-t}{2(1-t)\log t}dt$$

We might therefore have

$$\log(2\pi)=\int_0^1 \frac{3-t}{(1-t)\log t}dt$$

Alternatively, with $s=2$ we get

$$\varsigma'(2)=\int_0^1 \frac{\varsigma(2)t-Li_2(t)}{t\log t}dt$$

We note from (E.22bii) that

$$2\varsigma(3)+\varsigma'(2)-\gamma\varsigma(2)=\int_0^1 Li_2(t)\left[\frac{1}{t-1}-\frac{1}{t\log t}\right]dt$$

$$=\int_0^1 \left[\frac{Li_2(t)}{t-1}-\frac{Li_2(t)}{t\log t}\right]dt$$

$$=\int_0^1 \left[\frac{Li_2(t)}{t-1}-\frac{\varsigma(2)}{\log t}+\frac{\varsigma(2)}{\log t}-\frac{Li_2(t)}{t\log t}\right]dt$$

$$=\int_0^1 \left[\frac{Li_2(t)}{t-1}-\frac{\varsigma(2)}{\log t}\right]dt+\int_0^1 \left[\frac{\varsigma(2)}{\log t}-\frac{Li_2(t)}{t\log t}\right]dt$$

$$=\int_0^1 \left[\frac{Li_2(t)}{t-1}-\frac{\varsigma(2)}{\log t}\right]dt+\varsigma'(2)$$

We therefore see that

(E.22biv) $$2\varsigma(3)-\gamma\varsigma(2)=\int_0^1\left[\frac{Li_2(t)}{t-1}-\frac{\varsigma(2)}{\log t}\right]dt$$

which we may write as

$$2\varsigma(3)-\gamma\varsigma(2)=\lim_{x\to 1}\int_0^x\left[\frac{Li_2(t)}{t-1}-\frac{\varsigma(2)}{\log t}\right]dt$$

Using integration by parts we can easily determine that

$$\int\frac{Li_2(t)}{t-1}dt=Li_2(t)\log(1-t)+\log t\log^2(1-t)+2Li_2(1-t)\log(1-t)-2Li_3(1-t)$$

and we therefore have

$$\int_0^x\frac{Li_2(t)}{t-1}dt=Li_2(x)\log(1-x)+\log x\log^2(1-x)+2Li_2(1-x)\log(1-x)-2Li_3(1-x)+2\varsigma(3)$$

We have the limit

$$\lim_{x\to 1}\int_0^x\frac{Li_2(t)}{t-1}dt=\lim_{x\to 1}[Li_2(x)\log(1-x)+\log x\log^2(1-x)]+2\varsigma(3)$$

which then suggests that

$$-\gamma\varsigma(2)=\lim_{x\to 1}\left[Li_2(x)\log(1-x)+\log x\log^2(1-x)-\int_0^x\frac{\varsigma(2)}{\log t}dt\right]$$

We therefore have

(E.22bv) $$\gamma\varsigma(2)=\lim_{x\to 1}\left[-Li_2(x)\log(1-x)-\log x\log^2(1-x)+\varsigma(2)li(x)\right]$$

where $li(x)$ is the logarithmic integral. Nielsen [104a, p.3] has shown that

(E.22bvi) $$li(x)=\gamma+\log(-\log x)+\sum_{n=1}^{\infty}\frac{\log^n x}{n!n}$$

and hence we have

(E.22bvii) $$\gamma\varsigma(2)=\lim_{x\to 1}\left[-Li_2(x)\log(1-x)-\log x\log^2(1-x)+\varsigma(2)[\gamma+\log(-\log x)]\right]$$

which gives us

(E.22bviii) $$\lim_{x\to 1}\left[-Li_2(x)\log(1-x)-\log x\log^2(1-x)+\varsigma(2)\log(-\log x)\right]=0$$

Using Euler's identity (1.6c) we may write this as

$$\lim_{x\to 1}\left[Li_2(1-x)\log(1-x)+\varsigma(2)\log(1-x)+\varsigma(2)\log(-\log x)\right]=0$$

and hence we have

$$\lim_{x\to 1}\left[\log(1-x)+\log(-\log x)\right]=0$$

or the more obvious limit

$$\lim_{x\to 1}\left[\log[(x-1)\log x]\right]=0$$

The following was posed as a Quickie Problem Q974 in the October 2007 issue of Mathematics Magazine: prove that

$$\sum_{k=2}^{\infty}(-1)^k\frac{\varsigma(k)}{k+1}=1+\frac{1}{2}\gamma-\frac{1}{2}\log(2\pi)$$

We have Euler's formula for the gamma function

$$\Gamma(x)=\frac{1}{x}\prod_{n=1}^{\infty}\left[\left(1+\frac{1}{n}\right)^x\left(1+\frac{x}{n}\right)^{-1}\right]$$

and multiplying by $x$ and taking the logarithm of both sides gives us

$$\log\Gamma(x+1)=\sum_{n=1}^{\infty}\left[x\log\left(1+\frac{1}{n}\right)-\log\left(1+\frac{x}{n}\right)\right]$$

This may be written as

$$\log\Gamma(x+1)=\sum_{n=1}^{\infty}\left[x\log\left(1+\frac{1}{n}\right)+\sum_{k=1}^{\infty}\frac{(-1)^k x^k}{kn^k}\right]$$

$$=\sum_{n=1}^{\infty}\left[x\left[\log\left(1+\frac{1}{n}\right)-\frac{1}{n}\right]+\sum_{k=2}^{\infty}\frac{(-1)^k x^k}{kn^k}\right]$$

$$= x\sum_{n=1}^{\infty}\left[\log\left(1+\frac{1}{n}\right)-\frac{1}{n}\right]+\sum_{n=1}^{\infty}\sum_{k=2}^{\infty}\frac{(-1)^k x^k}{kn^k}$$

$$= x\sum_{n=1}^{\infty}\left[\log\left(1+\frac{1}{n}\right)-\frac{1}{n}\right]+\sum_{k=2}^{\infty}\frac{(-1)^k x^k}{k}\sum_{n=1}^{\infty}\frac{1}{n^k}$$

and hence we have the Maclaurin series for $\log\Gamma(x+1)$

$$\log\Gamma(x+1) = -\gamma x + \sum_{k=2}^{\infty}\frac{(-1)^k \varsigma(k)}{k}x^k$$

Differentiation results in

$$\psi(x+1) = -\gamma + \sum_{k=2}^{\infty}(-1)^k \varsigma(k) x^{k-1}$$

and we now multiply this equation by $x$ and integrate to obtain

$$\int_0^u x\psi(x+1)\,dx = -\frac{1}{2}\gamma u^2 + \sum_{k=2}^{\infty}(-1)^k \frac{\varsigma(k)}{k+1}u^{k+1}$$

Integration by parts gives us

$$\int_0^u x\psi(x+1)\,dx = x\log\Gamma(x+1)\Big|_0^u - \int_0^u \log\Gamma(x+1)\,dx$$

$$= u\log\Gamma(u+1) - \int_0^u \log\Gamma(x+1)\,dx$$

Hence we obtain

$$\int_0^u \log\Gamma(x+1)\,dx = u\log\Gamma(u+1) + \frac{1}{2}\gamma u^2 - \sum_{k=2}^{\infty}(-1)^k \frac{\varsigma(k)}{k+1}u^{k+1}$$

or alternatively

$$\int_0^u \log\Gamma(x)\,dx = u\log\Gamma(u) + u + \frac{1}{2}\gamma u^2 - \sum_{k=2}^{\infty}(-1)^k \frac{\varsigma(k)}{k+1}u^{k+1}$$

We have the well-known Raabe's integral

$$\int_0^1 \log\Gamma(x)\,dx = \frac{1}{2}\log(2\pi)$$

and hence with $u=1$ we have

$$\frac{1}{2}\log(2\pi) = 1 + \frac{1}{2}\gamma - \sum_{k=2}^{\infty} (-1)^k \frac{\varsigma(k)}{k+1}$$

The result may be generalised by reference to Alexeiewsky's theorem: in 1894 Alexeiewsky showed that

$$\int_0^u \log\Gamma(x+1)\,dx = \frac{1}{2}u(u-1) + \frac{1}{2}\log(2\pi) - \log G(u+1) + u\log\Gamma(u)$$

where $G(u)$ is the Barnes double gamma function defined by [126, p.25]

$$G(1+x) = (2\pi)^{x/2} \exp\left[-\frac{1}{2}(\gamma x^2 + x^2 + x)\right] \prod_{k=1}^{\infty} \left\{ \left(1+\frac{x}{k}\right)^k \exp\left(\frac{x^2}{2k} - x\right) \right\}$$

It is obvious from the definition that $G(1) = 1$.

We then obtain the following identity originally derived by Srivastava [126, p.210] in 1988

$$\sum_{k=2}^{\infty} (-1)^k \frac{\varsigma(k)}{k+1} u^{k+1} = \left[1 - \log(2\pi)\right]\frac{u}{2} + (1+\gamma)\frac{u^2}{2} + \log G(u+1)$$

□

From [75, p.351] we have

$$-\frac{\log(1-x)}{(1-x)^{m+1}} = \sum_{n=0}^{\infty} (H_{m+n} - H_m) \binom{m+n}{n} x^n$$

Reference to the Larcombe identity

$$m\binom{m+n}{n} \sum_{k=0}^{n} \binom{n}{k} \frac{(-1)^k}{(m+k)^2} = \sum_{k=m}^{m+n} \frac{1}{k} = H_{m+n} - H_m + \frac{1}{m}$$

then shows that

$$m\sum_{n=0}^{\infty} \binom{m+n}{n}^2 \sum_{k=0}^{n} \binom{n}{k} \frac{(-1)^k}{(m+k)^2} x^n - \frac{1}{m} \sum_{n=0}^{\infty} \binom{m+n}{n} x^n = -\frac{\log(1-x)}{(1-x)^{m+1}}$$

The following exposition is based on the interesting analysis given in Serret's classic text referred to in Appendix B.

From (E.10a) we have $\Gamma'(x) = \int_0^\infty t^{x-1} e^{-t} \log t \, dt$ and from (E.22aa) we have

$$\log t = \int_0^\infty \frac{e^{-y} - e^{-ty}}{y} dy$$

Therefore we may write

$$\Gamma'(x) = \int_0^\infty t^{x-1} e^{-t} \, dt \int_0^\infty \frac{e^{-y} - e^{-ty}}{y} dy$$

and interchanging the order of integration we get

$$\Gamma'(x) = \int_0^\infty \frac{dy}{y} \left[ e^{-y} \int_0^\infty t^{x-1} e^{-t} dt - \int_0^\infty t^{x-1} e^{-(1+y)t} dt \right]$$

$$= \Gamma(x) \int_0^\infty \left[ e^{-y} - \frac{1}{(1+y)^x} \right] \frac{dy}{y}$$

We therefore have

(E.22c) $$\psi(x) = \frac{d}{dx} \log \Gamma(x) = \frac{\Gamma'(x)}{\Gamma(x)} = \int_0^\infty \left[ e^{-y} - \frac{1}{(1+y)^x} \right] \frac{dy}{y}$$

and integrating the above over the interval $[1, x]$ we obtain

(E.22d) $$\log \Gamma(x) = \int_0^\infty \left[ (x-1) e^{-y} - \frac{(1+y)^{-1} - (1+y)^{-x}}{\log(1+y)} \right] \frac{dy}{y}$$

Since $\log \Gamma(2) = \log 1 = 0$ we have

(E.22e) $$0 = \int_0^\infty \left[ e^{-y} - \frac{y(1+y)^{-2}}{\log(1+y)} \right] \frac{dy}{y}$$

Multiplying (E.22e) by $x-1$ and using (E.22d) we find

(E.22f) $$\log \Gamma(x) = \int_0^\infty \left[ (x-1)(1+y)^{-2} - \frac{(1+y)^{-1} - (1+y)^{-x}}{y} \right] \frac{dy}{\log(1+y)}$$

Making the substitution $\log(1+y) = \alpha$ we obtain Malmstén's formula [126, p.16]

(E.22g) $$\log\Gamma(x) = \int_0^\infty \left[ (x-1)e^{-\alpha} - \frac{e^{-\alpha} - e^{-\alpha x}}{1-e^{-\alpha}} \right] \frac{d\alpha}{\alpha}$$

Later in this Appendix we show how (E.22g) is employed in the proof of Kummer's Fourier series expansion of $\log\Gamma(x)$.

Letting $x = 1$ in (E.22c) we get

$$\psi(1) = -\gamma = \int_0^\infty \left[ e^{-y} - \frac{1}{(1+y)} \right] \frac{dy}{y}$$

and subtracting this from (E.22c) we obtain

$$\psi(x) + \gamma = \int_0^\infty \left[ \frac{1}{1+y} - \frac{1}{(1+y)^x} \right] \frac{dy}{y}$$

Substituting $1 + y = 1/t$ we have the well-known result

(E.22gi) $$\psi(x) + \gamma = \int_0^1 \frac{1-t^{x-1}}{1-t} dt$$

and, where $x$ is a positive integer, this becomes

$$\psi(n) + \gamma = \int_0^1 \frac{1-t^{n-1}}{1-t} dt = H_{n-1}$$

Differentiating (E.22c) we get

$$\psi'(x) = \int_0^\infty \frac{\log(1+y)}{(1+y)^x} \frac{dy}{y}$$

Differentiating (E.22g) we get the following integral representation for $\psi(x)$ which is due to Gauss [126, p.15]

(E.22h) $$\psi(x) = \frac{d}{dx}\log\Gamma(x) = \int_0^\infty \left[ \frac{e^{-\alpha}}{\alpha} - \frac{e^{-\alpha x}}{1-e^{-\alpha}} \right] d\alpha$$

With $x = 1$ in (E.22h) we get the familiar result [25, p.177] for Euler's constant as reported in (E.22b)

(E.22i) $$\gamma = \int_0^\infty e^{-\alpha} \left[ \frac{1}{1-e^{-\alpha}} - \frac{1}{\alpha} \right] d\alpha$$

From (E.22h) we have formally

$$\psi(x) = \lim_{s\to1}\int_0^\infty \left[\frac{t^{s-1}e^{-t}}{t} - \frac{t^{s-1}e^{-xt}}{1-e^{-t}}\right]dt$$

$$= \lim_{s\to1}\left[\int_0^\infty t^{s-2}e^{-t}dt - \Gamma(s)\varsigma(s,x)\right]$$

$$= \lim_{s\to1}\left[\Gamma(s-1) - \Gamma(s)\varsigma(s,x)\right]$$

$$= \lim_{s\to1}\Gamma(s)\left[\frac{1}{s-1} - \varsigma(s,x)\right]$$

and we therefore see that

(E.22ii) $$\psi(x) = \lim_{s\to1}\left[\frac{1}{s-1} - \varsigma(s,x)\right]$$

which we previously saw in (4.3.202) in Volume II(b).

Differentiating (E.22ii) gives us

$$\psi'(x) = \frac{\partial}{\partial x}\lim_{s\to1}\left[\frac{1}{s-1} - \varsigma(s,x)\right] = -\frac{\partial}{\partial x}\lim_{s\to1}\varsigma(s,x) = \lim_{s\to1} s\varsigma(s+1,x) = \varsigma(2,x)$$

Differentiating (E.22h) we get

(E.22j) $$\psi'(x) = \frac{d^2}{dx^2}\log\Gamma(x) = \int_0^\infty \frac{\alpha e^{-\alpha x}}{1-e^{-\alpha}}d\alpha = \int_0^\infty \alpha e^{-\alpha x}\sum_{n=0}^\infty e^{-n\alpha}d\alpha = \sum_{n=0}^\infty \frac{1}{(x+n)^2}$$

Integrating (E.22j) over the interval $[1, x]$ we see that [126, p.14]

(E.22k) $$\psi(x) = \frac{d}{dx}\log\Gamma(x) = -\gamma - \frac{1}{x} + \sum_{n=1}^\infty \frac{x}{n(x+n)}$$

In (E.22j) let $x \to x+1$ to obtain

(E.22l) $$\frac{d^2}{dx^2}\log\Gamma(x+1) = \sum_{k=1}^\infty \frac{1}{(x+k)^2}$$

and differentiating (E.22l) $n-2$ times we get

(E.22m) $$\frac{1}{n!}\frac{d^n}{dx^n}\log\Gamma(x+1) = \frac{(-1)^n}{n}\sum_{k=1}^\infty \frac{1}{(x+k)^n}$$

Evaluating the above at $x = 0$ we have for $n > 1$

$$\frac{1}{n!}\frac{d^n}{dx^n}\log\Gamma(x+1)\Bigg|_{x=0} = \frac{(-1)^n}{n}\varsigma(n)$$

Accordingly we have the Maclaurin expansion

(E.22n) $$\log\Gamma(1+x) = -\gamma x + \sum_{n=2}^{\infty}(-1)^n \frac{\varsigma(n)}{n}x^n$$

and also see (E.6e) and (E.41).

Differentiating (E.22d) we obtain

$$\psi(x) = \int_0^\infty \left[e^{-y} - \frac{1}{(1+y)^x}\right]\frac{dy}{y}$$

$$\psi'(x) = \int_0^\infty \frac{\log(1+y)}{(1+y)^x}\frac{dy}{y}$$

$$\psi^{(p)}(x) = (-1)^{p+1}\int_0^\infty \frac{\log^p(1+y)}{(1+y)^x}\frac{dy}{y}$$

With the substitution $u = \log(1+y)$ this becomes

$$\psi^{(p)}(x) = (-1)^{p+1}\int_0^\infty \frac{u^p e^{-u(x-1)}}{e^u - 1}du$$

We note from [126, p.92] that

$$\varsigma(p+1,x)\Gamma(p+1) = \int_0^\infty \frac{u^p e^{-u(x-1)}}{e^u - 1}du$$

and hence we have [126, p.22]

$$\psi^{(p)}(x) = (-1)^{p+1}\Gamma(p+1)\varsigma(p+1,x)$$

We have from (A.2)

$$\frac{t}{e^t - 1} = 1 - \frac{t}{2} + \sum_{k=1}^{\infty} B_{2k}\frac{t^{2k}}{(2k)!} \qquad , (|t| < 2\pi)$$

With the inequality given by Póyla and Szegö [108a, p.212], namely

$$\sum_{k=1}^{2N} B_{2k} \frac{x^{2k-1}}{(2k)!} + \frac{1}{x} - \frac{1}{2} < \frac{1}{e^x - 1} < \sum_{k=1}^{2N+1} B_{2k} \frac{x^{2k-1}}{(2k)!} + \frac{1}{x} - \frac{1}{2}$$

and defining $a_n(x)$ by

$$a_n(x) = \frac{1}{e^x - 1} - \frac{e^{-x}}{n(e^{x/n} - 1)}$$

we have the inequality

$$\frac{1}{e^x - 1} - \frac{e^{-x}}{n}\left[ -\frac{1}{2} + \frac{n}{x} + \sum_{k=1}^{2N+1} B_{2k} \frac{x^{2k-1}}{n^{2k-1}(2k)!} \right] < a_n(x)$$

$$< \frac{1}{e^x - 1} - \frac{e^{-x}}{n}\left[ -\frac{1}{2} + \frac{n}{x} + \sum_{k=1}^{2N} B_{2k} \frac{x^{2k-1}}{n^{2k-1}(2k)!} \right]$$

With a little rearrangement this easily becomes

$$\frac{1}{e^x - 1} - \frac{e^{-x}}{x} + \frac{e^{-x}}{2n} - \frac{e^{-x}}{n}\left[ \sum_{k=1}^{2N+1} B_{2k} \frac{x^{2k-1}}{n^{2k-1}(2k)!} \right] < a_n(x)$$

$$< \frac{1}{e^x - 1} - \frac{e^{-x}}{x} + \frac{e^{-x}}{2n} - \frac{e^{-x}}{n}\left[ \sum_{k=1}^{2N} B_{2k} \frac{x^{2k-1}}{n^{2k-1}(2k)!} \right]$$

With integration we obtain

$$\int_0^\infty \frac{1}{e^x - 1} - \frac{e^{-x}}{x} + \frac{e^{-x}}{2n} - \frac{e^{-x}}{n}\left[ \sum_{k=1}^{2N+1} B_{2k} \frac{x^{2k-1}}{n^{2k-1}(2k)!} \right] dx < a_n$$

$$< \int_0^\infty \frac{1}{e^x - 1} - \frac{e^{-x}}{x} + \frac{e^{-x}}{2n} - \frac{e^{-x}}{n}\left[ \sum_{k=1}^{2N} B_{2k} \frac{x^{2k-1}}{n^{2k-1}(2k)!} \right] dx$$

because $a_n = \int_0^\infty a_n(x)\, dx$.

Since $\gamma = \int_0^\infty \left[ \frac{1}{e^x - 1} - \frac{e^{-x}}{x} \right] dx$ we have

$$\gamma + \frac{1}{2n} - \int_0^\infty \frac{e^{-x}}{n}\left[ \sum_{k=1}^{2N+1} B_{2k} \frac{x^{2k-1}}{n^{2k-1}(2k)!} \right] dx < a_n < \gamma + \frac{1}{2n} - \int_0^\infty \frac{e^{-x}}{n}\left[ \sum_{k=1}^{2N} B_{2k} \frac{x^{2k-1}}{n^{2k-1}(2k)!} \right] dx$$

Since by definition $\Gamma(\nu)=\int_0^\infty x^{\nu-1}e^{-x}dx$, we have $\int_0^\infty x^{2k-1}e^{-x}dx=\Gamma(2k)=(2k-1)!$ and hence

$$\int_0^\infty \frac{e^{-x}}{n}\left[\sum_{k=1}^{2N+1} B_{2k}\frac{x^{2k-1}}{n^{2k-1}(2k)!}\right]dx=\sum_{k=1}^{2N+1}\frac{B_{2k}}{2kn^{2k}}$$

We accordingly obtain

$$\gamma+\frac{1}{2n}-\sum_{k=1}^{2N+1}\frac{B_{2k}}{2kn^{2k}} < a_n < \gamma+\frac{1}{2n}-\sum_{k=1}^{2N}\frac{B_{2k}}{2kn^{2k}}$$

and therefore we have

(E.23) $$\gamma+\frac{1}{2n}-\sum_{k=1}^{2N+1}\frac{B_{2k}}{2kn^{2k}} < H_n-\log n < \gamma+\frac{1}{2n}-\sum_{k=1}^{2N}\frac{B_{2k}}{2kn^{2k}}$$

Hence we have the asymptotic expansion

(E.24) $$H_n-\log n\approx \gamma+\frac{1}{2n}-\sum_{k=1}^{2N+1}\frac{B_{2k}}{2kn^{2k}}$$

This formula is not quite new: an equivalent identity was obtained by Euler in 1736 by using the Euler-Maclaurin expansion which he had previously discovered in 1732 (and which was independently discovered by Maclaurin and published by him in 1742). A brief summary of the history of the Euler-Maclaurin expansion is given in [135, p.127] which also contains several other historical references. Proofs of the Euler-Maclaurin expansion are given, for example, in [14a], [135, p.127] and [75, p.469]. In 1962, using what we would now classify as a rudimentary computer, Knuth [90a] calculated $\gamma$ to 1,271 decimal places using the above formula (Euler having previously computed 16 digits in 1736 with $n=10$ and $N=3$).

From (E.24) it is clear that

(E.25) $$\lim_{n\to\infty} n\left[H_n-\log n-\gamma\right]=\frac{1}{2}$$

Since $\lim_{n\to\infty}\frac{\log n}{n^p}=0$ we also have from (E.24)

(E.26) $$\lim_{n\to\infty}\log n\left[H_n-\log n-\gamma\right]=0$$

We recall Adamchik's formula (4.1.14)

$$\sum_{k=1}^{n}\frac{H_k}{k}=\frac{1}{2}\left(H_n^{(1)}\right)^2+\frac{1}{2}H_n^{(2)}$$

Hence we obtain

$$\left[\sum_{k=1}^{n}\frac{H_k}{k}-\frac{1}{2}\left(H_n^{(1)}\right)^2\right]=\frac{1}{2}H_n^{(2)}$$

and substituting $\left(H_n^{(1)}\right)^2=\left(H_n^{(1)}-\log n-\gamma+\log n+\gamma\right)^2$ we have

(E.26a)

$$\left[\sum_{k=1}^{n}\frac{H_k}{k}-\frac{1}{2}\left(H_n^{(1)}-\log n-\gamma\right)^2-\left(H_n^{(1)}-\log n-\gamma\right)\left(\log n+\gamma\right)-\frac{1}{2}\left(\log n+\gamma\right)^2\right]=\frac{1}{2}H_n^{(2)}$$

Therefore we get

(E.27) $$\lim_{n\to\infty}\left[\sum_{k=1}^{n}\frac{H_k}{k}-\frac{1}{2}\left(\log n+\gamma\right)^2\right]=\frac{1}{2}\varsigma(2)$$

which implies that

(E.28) $$\lim_{n\to\infty}\left[\sum_{k=1}^{n}\frac{H_k}{k}-\gamma\log n-\frac{1}{2}\log^2 n\right]=\frac{1}{2}\left[\varsigma(2)+\gamma^2\right]=\frac{1}{2}\Gamma''(1)$$

I subsequently discovered that this formula had been previously discovered by Kanemitsu et al [82c] in 2004. They used this to show that

$$\sum_{n=1}^{\infty}H_n\left(\log\frac{n+1}{n}-\frac{1}{n}\right)=-\frac{1}{2}\left[\varsigma(2)+\gamma^2+2\gamma_1\right]$$

We may also write (E.28) as in (E.6d)

(E.29) $$\lim_{n\to\infty}\left[\frac{1}{2}\left(H_n^{(1)}\right)^2-\gamma\log n-\frac{1}{2}\log^2 n\right]=\frac{1}{2}\gamma^2$$

Equation (E.28) concurs with the asymptotic formula obtained by Flajolet and Sedgewick [68] (see also (3.13) of Volume I for more background information)

(E.30) $$-S_n(2)=\sum_{k=1}^{n}\frac{H_k}{k}=\frac{1}{2}\log^2 n+\gamma\log n+\frac{1}{2}\left[\varsigma(2)+\gamma^2\right]+O\left(\frac{\log n}{n}\right)$$

With reference to equation (E.29) we have

$$\left[\frac{1}{2}\left(H_n^{(1)}\right)^2-\gamma\log n-\frac{1}{2}\log^2 n\right]+\frac{1}{2}\left(H_n^{(1)}-\log n\right)^2=\left(H_n^{(1)}\right)^2-H_n^{(1)}\log n-\gamma\log n$$

and hence we see that

$$\lim_{n\to\infty}\left[\left(H_n^{(1)}\right)^2 - H_n^{(1)}\log n - \gamma\log n\right] = \gamma^2$$

Similarly, subtraction results in (E.26).

By extending the above procedure to $S_n(3)$, defined by (3.16c), and using (3.20c) we get

(E.30a)

$$-S_n(3) = -\sum_{k=1}^{n}\binom{n}{k}\frac{(-1)^k}{k^3} = \frac{1}{6}\left(H_n^{(1)}\right)^3 + \frac{1}{2}H_n^{(1)}H_n^{(2)} + \frac{1}{3}H_n^{(3)} = \frac{1}{2}\left\{\sum_{k=1}^{n}\frac{\left(H_k^{(1)}\right)^2}{k} + \sum_{k=1}^{n}\frac{H_k^{(2)}}{k}\right\}$$

The Flajolet and Sedgewick general asymptotic formula is

(E.30b)

$$-S_n(m) \approx \sum_{1m_1+2m_2+3m_3\ldots=m}\frac{1}{m_1!\,m_2!\,m_3!\ldots}\left(\log n+\gamma\right)^{m_1}\left(\frac{\varsigma(2)}{2}\right)^{m_2}\left(\frac{\varsigma(3)}{3}\right)^{m_3}\left(\frac{\varsigma(4)}{4}\right)^{m_4}\ldots$$

and, in particular, they obtained

(E.31)

$$-S_n(3) = \frac{1}{6}\log^3 n + \frac{\gamma}{2}\log^2 n + \frac{1}{2}\left[\varsigma(2)+\gamma^2\right]\log n + \frac{1}{2}\left[\varsigma(2)+\frac{1}{3}\gamma^2\right]\gamma + \frac{1}{3}\varsigma(3) + O\left(\frac{\log^2 n}{n}\right)$$

Using L'Hôpital's rule it is easily seen that $\lim_{n\to\infty}\left(\frac{\log^2 n}{n}\right) = 0$.

We can easily derive another one of the asymptotic formulae obtained by Flajolet and Sedgewick. By simple rearrangement of (E.30a) we have

$$\frac{1}{2}\left[\sum_{k=1}^{n}\frac{\left(H_k^{(1)}\right)^2}{k} + \sum_{k=1}^{n}\frac{H_k^{(2)}}{k}\right] - \frac{1}{6}\left(H_n^{(1)}\right)^3 - \frac{1}{2}H_n^{(1)}H_n^{(2)} = \frac{1}{3}H_n^{(3)}$$

and hence we get

(E.32)
$$\lim_{n\to\infty}\left[\frac{1}{2}\sum_{k=1}^{n}\frac{\left(H_k^{(1)}\right)^2}{k} + \frac{1}{2}\sum_{k=1}^{n}\frac{H_k^{(2)}}{k} - \frac{1}{6}H_n^{(3)} - \frac{1}{2}H_n^{(1)}H_n^{(2)}\right] = \frac{1}{3}\varsigma(3)$$

(E.32a) $$\lim_{n\to\infty}\left[\sum_{k=1}^{n}\frac{\left(H_k^{(1)}\right)^2}{k}+\sum_{k=1}^{n}\frac{H_k^{(2)}}{k}-H_n^{(1)}H_n^{(2)}\right]=\frac{4}{3}\varsigma(3)$$

We now use the following representations

$$\left(H_n^{(1)}\right)^3=\left(H_n^{(1)}-\log n-\gamma+\log n+\gamma\right)^3$$

$$=\left(H_n^{(1)}-\log n-\gamma\right)^3+3\left(H_n^{(1)}-\log n-\gamma\right)^2(\log n+\gamma)+3\left(H_n^{(1)}-\log n-\gamma\right)\left(\log n+\gamma\right)^2+\left(\log n+\gamma\right)^3$$

and

$$-\frac{1}{2}H_n^{(2)}=\frac{1}{2}\left(H_n^{(1)}\right)^2-\sum_{k=1}^{n}\frac{H_k}{k}$$

$$=\frac{1}{2}\left(H_n^{(1)}-\log n-\gamma+\log n+\gamma\right)^2-\sum_{k=1}^{n}\frac{H_k}{k}$$

$$=\frac{1}{2}\left(H_n^{(1)}-\log n-\gamma\right)^2+\left(H_n^{(1)}-\log n-\gamma\right)\left(\log n+\gamma\right)+\frac{1}{2}\left(\log n+\gamma\right)^2-\sum_{k=1}^{n}\frac{H_k}{k}$$

$$=\frac{1}{2}\left(H_n^{(1)}-\log n-\gamma\right)^2+\left(H_n^{(1)}-\log n-\gamma\right)\left(\log n+\gamma\right)+\frac{1}{2}\left(\log n+\gamma\right)^2$$

$$-\left[\sum_{k=1}^{n}\frac{H_k}{k}-\gamma\log n-\frac{1}{2}\log^2 n\right]-\gamma\log n-\frac{1}{2}\log^2 n$$

Hence, with some algebra, we obtain

$$-\frac{1}{2}H_n^{(1)}H_n^{(2)}=$$

$$\frac{1}{2}H_n^{(1)}\left(H_n^{(1)}-\log n-\gamma\right)^2+H_n^{(1)}\left(H_n^{(1)}-\log n-\gamma\right)\left(\log n+\gamma\right)$$

$$+\frac{1}{2}\left(H_n^{(1)}-\log n-\gamma\right)\left(\log n+\gamma\right)^2+\frac{1}{2}\left(\log n+\gamma\right)^3$$

$$-\left(H_n^{(1)}-\log n-\gamma\right)\left[\sum_{k=1}^{n}\frac{H_k}{k}-\gamma\log n-\frac{1}{2}\log^2 n\right]$$

$$-\left(\log n+\gamma\right)\left[\sum_{k=1}^{n}\frac{H_k}{k}-\gamma\log n-\frac{1}{2}\log^2 n-\frac{1}{2}\left(\varsigma(2)+\gamma^2\right)\right]$$

$$-\left(H_n^{(1)}-\log n-\gamma\right)\left(\gamma\log n+\frac{1}{2}\log^2 n\right)$$

$$-\left(\log n+\gamma\right)\left(\gamma\log n+\frac{1}{2}\log^2 n\right)-\frac{1}{2}\left(\log n+\gamma\right)\left(\varsigma(2)+\gamma^2\right)$$

Therefore, in the limit as $n\to\infty$ we have

$$\lim_{n\to\infty}\left[\begin{array}{l}\frac{1}{2}\sum_{k=1}^{n}\frac{\left(H_k^{(1)}\right)^2}{k}+\frac{1}{2}\sum_{k=1}^{n}\frac{H_k^{(2)}}{k}-\frac{1}{6}\left(\log n+\gamma\right)^3+\frac{1}{2}\left(\log n+\gamma\right)^3\\ -\left(\log n+\gamma\right)\left(\gamma\log n+\frac{1}{2}\log^2 n\right)-\frac{1}{2}\left(\log n+\gamma\right)\left(\varsigma(2)+\gamma^2\right)\end{array}\right]=\frac{1}{3}\varsigma(3)$$

This may be easily simplified and written as

(E.33)

$$\lim_{n\to\infty}\left[-S_n(3)-\frac{1}{6}\log^3 n-\frac{\gamma}{2}\log^2 n-\frac{1}{2}\left(\varsigma(2)+\gamma^2\right)\log n-\frac{1}{2}\varsigma(2)\gamma-\frac{1}{6}\gamma^3\right]=\frac{1}{3}\varsigma(3)$$

and this is equivalent to the result obtained by Flajolet and Sedgewick [68]. This may also be recorded as

(E.33a)

$$\lim_{n\to\infty}\left[-S_n(3)-\frac{1}{6}\log^3 n-\frac{\gamma}{2}\log^2 n-\frac{1}{2}\left(\varsigma(2)+\gamma^2\right)\log n\right]=\frac{1}{3}\varsigma(3)+\frac{1}{2}\varsigma(2)\gamma+\frac{1}{6}\gamma^3$$

or alternatively

(E.33b) $$\lim_{n\to\infty}\left[-S_n(3)-\frac{1}{6}\log^3 n-\frac{\gamma}{2}\log^2 n-\frac{1}{2}\left(\varsigma(2)+\gamma^2\right)\log n\right]=\frac{1}{6}\Gamma^{(3)}(1)$$

where in the last equation we have employed the identity in (E.16e) [25, p.213] for the derivative of the gamma function. Can (E.33b) be extended to give a result of the form $\lim_{n\to\infty}\left[-S_n(p)-\text{other factors }\right]=\frac{1}{p!}\Gamma^{(p)}(1)$ ? Wilf's result, which is referred to below, appears to be germane.

Using (E.30a), equation (E.33a) may also be written as

(E.33c)

$$\lim_{n\to\infty}\left[\frac{1}{6}\left(H_n^{(1)}\right)^3+\frac{1}{2}H_n^{(1)}H_n^{(2)}-\frac{1}{6}\log^3 n-\frac{\gamma}{2}\log^2 n-\frac{1}{2}\left(\varsigma(2)+\gamma^2\right)\log n\right]=\frac{1}{2}\varsigma(2)\gamma+\frac{1}{6}\gamma^3$$

From (E.29) we see that

$$\lim_{n\to\infty}\left[\frac{1}{2}\left(H_n^{(1)}\right)^2-\gamma\log n-\sum_{k=1}^{n}\frac{\log k}{k}+\sum_{k=1}^{n}\frac{\log k}{k}-\frac{1}{2}\log^2 n\right]=\frac{1}{2}\gamma^2$$

and referring to (4.3.224) in Volume II(b) for the Stieltjes constant $\gamma_1$ we get

(E.33ci) $$\lim_{n\to\infty}\left[\frac{1}{2}\left(H_n^{(1)}\right)^2-\gamma\log n-\sum_{k=1}^{n}\frac{\log k}{k}\right]=\frac{1}{2}\gamma^2-\gamma_1$$

From (E.28) we have

$$\lim_{n\to\infty}\left[\sum_{k=1}^{n}\frac{H_k}{k}-\gamma\log n-\frac{1}{2}\log^2 n\right]=\frac{1}{2}\left[\varsigma(2)+\gamma^2\right]$$

which may be written as

$$\lim_{n\to\infty}\left[\sum_{k=1}^{n}\frac{H_k}{k}-\sum_{k=1}^{n}\frac{\log k}{k}-\gamma\log n+\sum_{k=1}^{n}\frac{\log k}{k}-\frac{1}{2}\log^2 n\right]=\frac{1}{2}\left[\varsigma(2)+\gamma^2\right]$$

Then using the definition of Stieltjes constant $\gamma_1$ we have

$$\lim_{n\to\infty}\left[\sum_{k=1}^{n}\frac{H_k}{k}-\sum_{k=1}^{n}\frac{\log k}{k}-\gamma\log n\right]=\frac{1}{2}\left[\varsigma(2)+\gamma^2\right]-\gamma_1$$

or equivalently

$$\lim_{n\to\infty}\left[\sum_{k=1}^{n}\frac{H_k-\log k}{k}-\gamma\log n\right]=\frac{1}{2}\left[\varsigma(2)+\gamma^2\right]-\gamma_1$$

This may also be expressed as

$$\lim_{n\to\infty}\left[\sum_{k=1}^{n}\frac{H_k-\gamma-\log k}{k}+\sum_{k=1}^{n}\frac{\gamma}{k}-\gamma\log n\right]=\frac{1}{2}\left[\varsigma(2)+\gamma^2\right]-\gamma_1$$

or

$$\lim_{n\to\infty}\left[\sum_{k=1}^{n}\frac{H_k-\gamma-\log k}{k}+\gamma(H_n-\log n)\right]=\frac{1}{2}\left[\varsigma(2)+\gamma^2\right]-\gamma_1$$

and since $\lim_{n\to\infty}\gamma(H_n - \log n) = \gamma^2$ we therefore conclude that

$$\sum_{k=1}^{\infty}\frac{H_k - \gamma - \log k}{k} = \frac{1}{2}\left[\varsigma(2) - \gamma^2\right] - \gamma_1$$

This formula was recently proposed as a problem by Furdui [65a].

From (E.23) we see that

$$\frac{1}{2}\varsigma(2) - \sum_{n=1}^{\infty}\frac{1}{n}\sum_{k=1}^{2N+1}\frac{B_{2k}}{2kn^{2k}} < \sum_{n=1}^{\infty}\frac{H_n - \log n - \gamma}{n} < \frac{1}{2}\varsigma(2) - \sum_{n=1}^{\infty}\frac{1}{n}\sum_{k=1}^{2N}\frac{B_{2k}}{2kn^{2k}}$$

and hence we get

$$\sum_{n=1}^{\infty}\frac{1}{n}\sum_{k=1}^{2N+1}\frac{B_{2k}}{2kn^{2k}} > \frac{1}{2}\gamma^2 + \gamma_1 > \sum_{n=1}^{\infty}\frac{1}{n}\sum_{k=1}^{2N}\frac{B_{2k}}{2kn^{2k}}$$

The analysis may be extended as follows; from (E.33) we see that

$$\lim_{n\to\infty}\left[\sum_{k=1}^{n}\frac{\left(H_k^{(1)}\right)^2}{k} + \sum_{k=1}^{n}\frac{H_k^{(2)}}{k} - \frac{1}{3}\log^3 n - \gamma\log^2 n - \left(\varsigma(2)+\gamma^2\right)\log n\right] = \varsigma(2)\gamma + \frac{1}{3}\gamma^3 + \frac{2}{3}\varsigma(3)$$

The part in parenthesis may be written as

$$\sum_{k=1}^{n}\frac{\left(H_k^{(1)}\right)^2}{k} + \sum_{k=1}^{n}\frac{H_k^{(2)}}{k} - \sum_{k=1}^{n}\frac{\log^2 k}{k} + \left(\sum_{k=1}^{n}\frac{\log^2 k}{k} - \frac{1}{3}\log^3 n\right) - 2\gamma\sum_{k=1}^{n}\frac{\log k}{k}$$

$$+2\gamma\left(\sum_{k=1}^{n}\frac{\log k}{k} - \frac{1}{2}\log^2 n\right) + \left(\varsigma(2)+\gamma^2\right)H_n^{(1)} - \left(\varsigma(2)+\gamma^2\right)\left(H_n^{(1)} - \log n\right)$$

or equivalently

$$\sum_{k=1}^{n}\frac{\left(H_k^{(1)}\right)^2 + H_k^{(2)} - \log^2 k - 2\gamma\log k + \left(\varsigma(2)+\gamma^2\right)}{k} + \left(\sum_{k=1}^{n}\frac{\log^2 k}{k} - \frac{1}{3}\log^3 n\right)$$

$$+2\gamma\left(\sum_{k=1}^{n}\frac{\log k}{k} - \frac{1}{2}\log^2 n\right) - \left(\varsigma(2)+\gamma^2\right)\left(H_n^{(1)} - \log n\right)$$

and hence we obtain

$$\sum_{k=1}^{\infty} \frac{\left(H_k^{(1)}\right)^2 + H_k^{(2)} - \log^2 k - 2\gamma \log k + \left(\varsigma(2)+\gamma^2\right)}{k} = 2\varsigma(2)\gamma + \frac{4}{3}\gamma^3 + \frac{2}{3}\varsigma(3) - 2\gamma\gamma_1 - \gamma_2$$

or

$$\sum_{k=1}^{\infty} \frac{\left(H_k^{(1)}\right)^2 + H_k^{(2)} - \left(\log k + \gamma\right)^2 + \varsigma(2) + 2\gamma^2}{k} = 2\varsigma(2)\gamma + \frac{4}{3}\gamma^3 + \frac{2}{3}\varsigma(3) - 2\gamma\gamma_1 - \gamma_2$$

In his paper "The asymptotic behaviour of the Stirling numbers of the first kind" [138b], Wilf proved that if $\begin{bmatrix} n \\ k \end{bmatrix}$ is the (signless) Stirling number of the first kind, then for each fixed integer $k \ge 2$ we have

(E.33d) $$\frac{1}{(n-1)!}\begin{bmatrix} n \\ k \end{bmatrix} = \lambda_1 \frac{\log^{k-1} n}{(k-1)!} + \lambda_2 \frac{\log^{k-2} n}{(k-2)!} + \ldots + \lambda_k + O\left(\frac{\log^{k-2} n}{n}\right)$$

where $\lambda_j$ are the coefficients in the expansion

(E.33e) $$\frac{1}{\Gamma(z)} = \sum_{j=1}^{\infty} \lambda_j z^j$$

The signless Stirling number of the first kind is defined as the absolute value of $s(n,k)$.
In particular, using the recurrence (obtained from the logarithmic derivative of the infinite product for $\Gamma(z)$)

(E.33f) $$\lambda_{n+1} = \frac{1}{n}\left\{\gamma\lambda_n + \sum_{j=0}^{n-2} (-1)^{n-j-1} \varsigma(n-j)\lambda_{j+1}\right\} \qquad n \ge 1; \lambda_1 = 1$$

we obtain

$$\lambda_1 = 1$$

$$\lambda_2 = \gamma$$

$$\lambda_3 = \frac{1}{12}(6\gamma^2 - \pi^2)$$

$$\lambda_4 = \frac{1}{12}\left[2\gamma^3 - \gamma\pi^2 + 4\varsigma(3)\right]$$

$$\lambda_5 = \frac{1}{1440}\left[60\gamma^4 - 60\gamma^2\pi^2 + \pi^4 + 480\varsigma(3)\right]$$

From (3.105i) in Volume I we see that

$$s(n,0) = \delta_{n,0}$$

$$s(n,1) = (-1)^{n+1}(n-1)!$$

$$s(n,2) = (-1)^n (n-1)! H_{n-1}$$

$$s(n,3) = (-1)^{n+1}\frac{(n-1)!}{2}\left\{\left(H_{n-1}\right)^2 - H_{n-1}^{(2)}\right\}$$

$$s(n,4) = (-1)^n \frac{(n-1)!}{6}\left\{\left(H_{n-1}\right)^3 - 3H_{n-1}H_{n-1}^{(2)} + 2H_{n-1}^{(3)}\right\}$$

and we therefore get using Wilf's theorem

$$\frac{1}{(n-1)!}\begin{bmatrix} n \\ 2 \end{bmatrix} = H_{n-1} = \log n + \gamma + O\left(\frac{1}{n}\right)$$

For the next Stirling number we have

$$\frac{1}{(n-1)!}\begin{bmatrix} n \\ 3 \end{bmatrix} = \frac{1}{2}\left[\left(H_{n-1}\right)^2 - H_{n-1}^{(2)}\right] = \frac{1}{2}\log^2 n + \gamma\log n + \frac{1}{12}(6\gamma^2 - \pi^2) + O\left(\frac{\log n}{n}\right)$$

Since

$$\left(H_{n-1}\right)^2 - H_{n-1}^{(2)} = \left(H_n\right)^2 - 2\frac{H_n}{n} + 2\frac{1}{n^2} - H_n^{(2)}$$

we have

$$\left(H_n\right)^2 - 2\frac{H_n}{n} + 2\frac{1}{n^2} - H_n^{(2)} - \frac{1}{2}\log^2 n - \gamma\log n = \frac{1}{12}(6\gamma^2 - \pi^2) + O\left(\frac{\log n}{n}\right)$$

and taking the limit as $n \to \infty$ we see that

$$\lim_{n\to\infty}\left\{\frac{1}{2}\left(H_n\right)^2 - \frac{H_n}{n} + \frac{1}{n^2} - \frac{1}{2}H_n^{(2)} - \frac{1}{2}\log^2 n - \gamma\log n\right\} = \frac{1}{12}(6\gamma^2 - \pi^2)$$

We see from (E.24) that

$$\lim_{n\to\infty}\frac{H_n - \log n}{n} = 0$$

and hence it is clear that $\lim_{n\to\infty}\frac{H_n}{n}=0$. Accordingly we get

$$\lim_{n\to\infty}\left\{\frac{1}{2}\left(H_n\right)^2-\frac{1}{2}H_n^{(2)}-\frac{1}{2}\log^2 n-\gamma\log n\right\}=\frac{1}{12}(6\gamma^2-\pi^2)$$

and it will be noted that the above is equivalent to the previous result (E.29).

We also have

$$\frac{1}{(n-1)!}\begin{bmatrix} n \\ 4 \end{bmatrix}=\frac{1}{6}\left\{\left(H_{n-1}\right)^3-3H_{n-1}H_{n-1}^{(2)}+2H_{n-1}^{(3)}\right\}$$

$$=\frac{1}{3!}\log^3 n+\frac{1}{2!}\gamma\log^2 n+\frac{1}{12}(6\gamma^2-\pi^2)\log n+\frac{1}{12}\left[2\gamma^3-\gamma\pi^2+4\varsigma(3)\right]+O\left(\frac{\log^2 n}{n}\right)$$

Simple algebra shows that

$$\left(H_{n-1}\right)^3-3H_{n-1}H_{n-1}^{(2)}+2H_{n-1}^{(3)}=$$

$$\left(H_n^{(1)}\right)^3-3\frac{\left[H_n^{(1)}\right]^2}{n}+6\frac{H_n^{(1)}}{n^2}-6\frac{1}{n^3}-3H_n^{(1)}H_n^{(2)}+3\frac{H_n^{(2)}}{n}+2H_n^{(3)}$$

and we therefore have

$$\left\{\frac{1}{6}\left(H_n^{(1)}\right)^3-\frac{1}{2}\frac{\left[H_n^{(1)}\right]^2}{n}+\frac{H_n^{(1)}}{n^2}-\frac{1}{n^3}-\frac{1}{2}H_n^{(1)}H_n^{(2)}+\frac{1}{2}\frac{H_n^{(2)}}{n}+\frac{1}{3}H_n^{(3)}\right\}=$$

$$=\frac{1}{3!}\log^3 n+\frac{1}{2!}\gamma\log^2 n+\frac{1}{12}(6\gamma^2-\pi^2)\log n+\frac{1}{12}\left[2\gamma^3-\gamma\pi^2+4\varsigma(3)\right]+O\left(\frac{\log^2 n}{n}\right)$$

Noting that $\lim_{n\to\infty}\frac{H_n^{(2)}}{n}=\varsigma(2).\lim_{n\to\infty}\frac{1}{n}=0$, this may be written as

(E.33g)

$$\lim_{n\to\infty}\left\{\frac{1}{6}\left(H_n^{(1)}\right)^3-\frac{1}{2}\frac{\left[H_n^{(1)}\right]^2}{n}-\frac{1}{2}H_n^{(1)}H_n^{(2)}-\frac{1}{6}\log^3 n-\frac{1}{2}\gamma\log^2 n-\frac{1}{12}(6\gamma^2-\pi^2)\log n\right\}$$

$$=\frac{1}{12}\left[2\gamma^3-\gamma\pi^2\right]$$

and it will be noted that this is similar to (E.33c) above. Subtracting (E.33g) from (E.33c) we obtain

(E.33h) $$\lim_{n\to\infty}\left\{H_n^{(1)}H_n^{(2)}+\frac{1}{2}\frac{\left[H_n^{(1)}\right]^2}{n}-\varsigma(2)\log n\right\}=\gamma\varsigma(2)$$

This may be written in a more obvious form

(E.33i) $$\lim_{n\to\infty}\left\{H_n^{(2)}\left[H_n^{(1)}-\log n\right]+\log n\left[H_n^{(2)}-\varsigma(2)\right]+\frac{1}{2}\frac{\left[H_n^{(1)}\right]^2}{n}\right\}=\gamma\varsigma(2)$$

We showed in (4.4.233x) in Volume IV that

$$\lim_{n\to\infty}\log n\left[H_n^{(2)}-\varsigma(2)\right]=0$$

and later in (E.58a) we see that $\lim_{n\to\infty}\dfrac{\left(H_n^{(1)}\right)^2}{n+1}=0$ (and the result therefore becomes clear).

□

The following is extracted from a very interesting series of papers written by Snowden [120aa, p.68] whilst he was an undergraduate in 2003.

Let us consider the function $f(x)$ with the following Maclaurin expansion

(E.33j) $$\log f(x)=b_0+\sum_{n=1}^{\infty}\frac{b_n}{n}x^n$$

and we wish to determine the coefficients $a_n$ such that

$$f(x)=\sum_{n=0}^{\infty}a_nx^n$$

By differentiating (E.33j) and multiplying the two power series, we get

$$na_n=\sum_{k=1}^{n}b_ka_{n-k}$$

Upon examination of this recurrence relation it is easy to see that

$$n!a_n=a_0[b_1,-b_2,b_3,...,(-1)^{n+1}b_n]$$

where the symbol $[a_1, a_2, a_3, ..., a_n]$ is defined as the $n \times n$ determinant

$$\begin{vmatrix} a_1 & a_2 & a_3 & a_4 & . & . & . & a_n \\ (n-1) & a_1 & a_2 & a_3 & . & . & . & a_{n-1} \\ 0 & (n-2) & a_1 & a_2 & . & . & . & a_{n-2} \\ 0 & 0 & (n-3) & a_1 & . & . & . & a_{n-3} \\ \vdots & \vdots & \vdots & \vdots & \vdots & \vdots & \vdots & \vdots \\ 0 & 0 & 0 & 0 & 0 & 0 & 1 & a_1 \end{vmatrix}$$

Since $\log f(0) = \log a_0 = b_0$ we have

$$f(x) = e^{b_0} \left[ 1 + \sum_{n=1}^{\infty} [b_1, -b_2, b_3, ..., (-1)^{n+1} b_n] \frac{x^n}{n!} \right]$$

Multiplying (E.33j) by $\alpha$ it is easily seen that

(E.33k) $$f^{\alpha}(x) = e^{-b_0} \left[ 1 + \sum_{n=1}^{\infty} [\alpha b_1, -\alpha b_2, \alpha b_3, ..., (-1)^{n+1} \alpha b_n] \frac{x^n}{n!} \right]$$

and, in particular, with $\alpha = -1$ we obtain

(E.33l) $$\frac{1}{f(x)} = e^{-b_0} \left[ 1 + \sum_{n=1}^{\infty} [-b_1, b_2, -b_3, ..., (-1)^n b_n] \frac{x^n}{n!} \right]$$

Differentiating (E.33k) with respect to $\alpha$ would give us an expression for $f^{\alpha}(x) \log f(x)$. Snowden's analysis may well turn out to be useful in proving the Riemann Hypothesis.

We note from (E.22n) that

$$\log \Gamma(1+x) = -\gamma x + \sum_{n=2}^{\infty} (-1)^n \frac{\varsigma(n)}{n} x^n$$

and hence we have

(E.33m) $$\Gamma(1+x) = 1 + \sum_{n=1}^{\infty} [-\gamma, -\varsigma(2), -\varsigma(3), ..., -\varsigma(n)] \frac{x^n}{n!}$$

and

(E.33n) $$\frac{1}{\Gamma(1+x)} = 1 + \sum_{n=1}^{\infty} [\gamma, \varsigma(2), \varsigma(3), ..., \varsigma(n)] \frac{x^n}{n!}$$

where $\varsigma(1)$ is defined as equal to $\gamma$. Since $\dfrac{1}{\Gamma(1+x)} = \dfrac{1}{x\Gamma(x)}$ equation (E.33n) may be compared with (E.33e)

$$\frac{1}{\Gamma(x)} = \sum_{n=1}^{\infty} \lambda_n x^n$$

Shen [120] showed that

(E.33o) $$\frac{f'(z)}{f(z)} = 2\sum_{n=1}^{\infty}(2^n - 1)\varsigma(n+1)z^n$$

where $f(z)$ is defined as

(E.33p) $$f(z) = \frac{2^{-2z}\Gamma\left(\frac{1}{2} - z\right)}{\sqrt{\pi}\Gamma\left(1-z\right)}$$

Integrating (E.33o) and noting that $f(0) = 1$ we see that

$$\log f(z) = 2\sum_{n=1}^{\infty}(2^n - 1)\varsigma(n+1)\frac{z^{n+1}}{n+1}$$

$$= 2\sum_{n=2}^{\infty}(2^{n-1} - 1)\varsigma(n)\frac{z^n}{n}$$

Therefore, with regard to (E.33j) we have $b_0 = b_1 = 0$ and $b_n = (2^{n-1} - 1)\varsigma(n)$.

Differentiating (E.33p) logarithmically we see that

$$\frac{f'(z)}{f(z)} = -2\log 2 - \psi\left(\frac{1}{2} - z\right) + \psi\left(1-z\right)$$

Differentiating again results in

$$\frac{f(z)f''(z) - \left[f'(z)\right]^2}{f^2(z)} = \psi'\left(\frac{1}{2} - z\right) - \psi'\left(1-z\right)$$

We have [126, p.20]

$$\psi\left(\frac{1}{2}\right) = -\gamma - 2\log 2 \qquad \psi\left(1\right) = -\gamma$$

$$\psi^{(n)}\left(\frac{1}{2}\right) = (-1)^{n+1} n!(2^{n+1}-1)\varsigma(n+1)$$

$$\psi^{(n)}\left(1\right) = (-1)^{n+1} n!\varsigma(n+1)$$

and hence we see that

$$f'(0) = -2\log 2 - \psi\left(\frac{1}{2}\right) + \psi\left(1\right) = 0$$

$$f''(0) = \psi'\left(\frac{1}{2}\right) - \psi'\left(1\right) = 2\varsigma(2)$$

The resulting Maclaurin expansion agrees with Shen's result.

Differentiating (E.33m) we get

$$\Gamma^{(p)}(1+x) = \sum_{n=1}^{\infty} [-\gamma, -\varsigma(2), -\varsigma(3), ..., -\varsigma(n)] \frac{n(n-1)...(n-p+1)x^{n-p}}{n!}$$

and hence we have the $p$ th derivativative of the gamma function (see [126, p.264]) in the form of a determinant

(E.33pi) $$\Gamma^{(p)}(1) = [-\varsigma(1), -\varsigma(2), -\varsigma(3), ..., -\varsigma(p)]$$

where $\varsigma(1) = \gamma$. See also the paper by Kölbig and Strampp [91ab] where they employ the Bell polynomials.

Using the integral identity (4.4.246)

$$\frac{n}{2}\int_0^1 x^{n-1}\log^2(1-x)\,dx = \sum_{k=1}^{n} \frac{H_k}{k}$$

we may write the limit as

$$\lim_{n\to\infty}\left[\frac{n}{2}\int_0^1 x^{n-1}\log^2(1-x)\,dx - \gamma\log n - \frac{n}{2}\int_0^1 x^{n-1}\log^2 n\,dx\right] = \frac{1}{2}\left[\varsigma(2)+\gamma^2\right]$$

We have

$$\frac{n}{2}\int_0^1 x^{n-1}\log^2(1-x)\,dx - \gamma\log n - \frac{n}{2}\int_0^1 x^{n-1}\log^2 n\,dx =$$

$$\frac{n}{2}\left(\int_0^1 x^{n-1}\log^2(1-x)\,dx-\gamma\int_0^1 x^{n-1}\log^2 n\,dx-\int_0^1 x^{n-1}\log^2 n\,dx\right)$$

$$=\frac{n}{2}\left(\int_0^1 x^{n-1}\left[\log^2(1-x)-\gamma-\log^2 n\right]dx\right)$$

Using the Maclaurin expansion of $\log(1+t)$ we obtain

(E.33q) $$\log\left(1+\frac{x}{k}\right)=\sum_{n=1}^{\infty}(-1)^{n+1}\frac{x^n}{nk^n}=\frac{x}{k}-\frac{x^2}{2k^2}+\frac{x^3}{3k^3}-...$$

and thereby construct the finite sum

(E.33r) $$\sum_{k=1}^{N}\log\left(1+\frac{x}{k}\right)=\sum_{k=1}^{N}\frac{x}{k}-\frac{1}{2}\sum_{k=1}^{N}\frac{x^2}{k^2}+\frac{1}{3}\sum_{k=1}^{N}\frac{x^3}{k^3}-...$$

$$S_N=\sum_{k=1}^{N}\log\left(1+\frac{x}{k}\right)-\sum_{k=1}^{N}\frac{x}{k}=-\frac{1}{2}\sum_{k=1}^{N}\frac{x^2}{k^2}+\frac{1}{3}\sum_{k=1}^{N}\frac{x^3}{k^3}-...$$

Provided $-1<x\le 1$ we note that $S_N$ is an alternating series whose terms $\to 0$ as $N\to\infty$ and hence $S_N$ is convergent. Therefore we have

(E.34) $$\sum_{k=1}^{\infty}\left[\log\left(1+\frac{x}{k}\right)-\frac{x}{k}\right]=\sum_{k=2}^{\infty}\frac{(-1)^{k+1}}{k}x^k\varsigma(k)$$

and using (E.13a) we see that for $|x|\le 1$

(E.34a) $$\log\Gamma(x)+\log x+\gamma x=\sum_{k=2}^{\infty}\frac{(-1)^k}{k}x^k\varsigma(k)$$

In particular we have

(E.34b) $$\gamma=\sum_{k=2}^{\infty}\frac{(-1)^k}{k}\varsigma(k)$$

(E.34c) $$\frac{1}{2}\log\pi-\log 2+\frac{1}{2}\gamma=\sum_{k=2}^{\infty}\frac{(-1)^k}{k2^k}\varsigma(k)$$

(E.34ci) $$\frac{1}{2}\log\pi-\frac{1}{2}\gamma=\sum_{k=2}^{\infty}\frac{\varsigma(k)}{k2^k}$$

Alternatively, referring to (4.3.80a) in Volume II(a) we have

$$\log\left(1+\frac{x}{k}\right)=\sum_{j=1}^{n-1}\frac{(-1)^{j+1}}{j}\left(\frac{x}{k}\right)^j+\frac{(-1)^{n-1}}{k^{n-1}}\int_0^x\frac{t^{n-1}}{k+t}dt$$

and making the summation

$$\sum_{k=1}^{\infty}\left[\log\left(1+\frac{x}{k}\right)-\frac{x}{k}\right]=\sum_{k=1}^{\infty}\left[\sum_{j=1}^{n-1}\frac{(-1)^{j+1}}{j}\left(\frac{x}{k}\right)^j-\frac{x}{k}\right]+\sum_{k=1}^{\infty}\frac{(-1)^{n-1}}{k^{n-1}}\int_0^x\frac{t^{n-1}}{k+t}dt$$

We have as $n\to\infty$

$$=\sum_{k=1}^{\infty}\left[\sum_{j=1}^{\infty}\frac{(-1)^{j+1}}{j}\left(\frac{x}{k}\right)^j-\frac{x}{k}\right]$$

$$=\sum_{k=1}^{\infty}\left[\sum_{j=2}^{\infty}\frac{(-1)^{j+1}}{j}\left(\frac{x}{k}\right)^j\right]$$

$$=\sum_{k=2}^{\infty}\frac{(-1)^{k+1}}{k}x^k\varsigma(k)$$

Differentiating (E.34a) results in for $x<1$

(E.34d) $$\psi(x)+\frac{1}{x}+\gamma=\sum_{k=2}^{\infty}(-1)^k x^{k-1}\varsigma(k)$$

(which is reported by Sebah and Gourdon [119] and is the Maclaurin series for $\psi(x+1)+\gamma$ ) and with $x=1/2$ this gives us

(E.34e) $$2(1-\log 2)=\sum_{k=2}^{\infty}\frac{(-1)^k\varsigma(k)}{2^{k-1}}$$

In a similar vein, we note that Glaisher has shown that (see Lewin [100, p.269])

(E.43f) $$2\log 2-1=\sum_{k=1}^{\infty}\frac{\varsigma(2k+1)}{2^{2k}}$$

which may be compared with (E.34e).

As an alternative to (E.33q) we have

$$\sum_{k=1}^{N-1}\log\left(1+\frac{x}{k}\right)=\log\left(\frac{1+x}{1}\right)\left(\frac{2+x}{2}\right)\left(\frac{3+x}{3}\right)...\left(\frac{N-1+x}{N-1}\right)$$

$$= \log \frac{\Gamma(N+x)}{\Gamma(x+1)\Gamma(N)} = -\log \frac{\Gamma(x+1)\Gamma(N)}{\Gamma(N+x)}$$

$$= -\log \frac{x\Gamma(x)\Gamma(N)}{\Gamma(N+x)} = -\log x - \log B(x,N)$$

It is an exercise in Whittaker & Watson [135, p.262] to show that

(E.35) $$\log B(p,q) = \log\left(\frac{p+q}{pq}\right) + \int_0^1 \frac{(1-v^p)(1-v^q)}{(1-v)\log v} dv \qquad p,q > 0$$

and this formula is attributed to Euler (see also [104b, p.187]). We therefore have

$$\sum_{k=1}^{N-1} \log\left(1+\frac{x}{k}\right) = -\log x - \log\left(\frac{x+N}{xN}\right) - \int_0^1 \frac{(1-v^x)(1-v^N)}{(1-v)\log v} dv$$

and this is equivalent to

(E.35a) $$\sum_{k=1}^{N} \log\left(1+\frac{x}{k}\right) = -\int_0^1 \frac{(1-v^x)(1-v^N)}{(1-v)\log v} dv$$

Letting $x=1$ we get

$$\sum_{k=1}^{N} \log\left(1+\frac{1}{k}\right) = \log(N+1) = -\int_0^1 \frac{1-v^N}{\log v} dv$$

and as per (E.35b) we have $H_N = \int_0^1 \frac{(1-v^N)}{1-v} dv$. Therefore we see that

(E.35ai) $$H_N - \log(N+1) = \int_0^1 (1-v^N)\left[\frac{1}{1-v} + \frac{1}{\log v}\right] dv$$

and as $N \to \infty$ we see that

$$\gamma = \int_0^1 \left[\frac{1}{1-v} + \frac{1}{\log v}\right] dv$$

Letting $x = 1/2$ in (E.35a) we get

$$\sum_{k=1}^{N} \log\left(1+\frac{1}{2k}\right) = -\int_0^1 \frac{1-v^N}{(1+\sqrt{v})\log v} dv$$

We obtain from (E.35ai)

(E.35aii) $$\sum_{n=1}^{\infty}\frac{H_n}{n^s}-\sum_{n=1}^{\infty}\frac{\log(n+1)}{n^s}=\gamma\varsigma(s)-\int_0^1 Li_s(v)\left[\frac{1}{1-v}+\frac{1}{\log v}\right]dv$$

and with $s=2$ this becomes

$$2\varsigma(3)-\sum_{n=1}^{\infty}\frac{\log(n+1)}{n^2}-\gamma\varsigma(2)=\int_0^1 Li_2(v)\left[\frac{1}{v-1}-\frac{1}{\log v}\right]dv$$

which may be compared with (E.22bii)

$$2\varsigma(3)+\varsigma'(2)-\gamma\varsigma(2)=\int_0^1 Li_2(v)\left[\frac{1}{v-1}-\frac{1}{v\log v}\right]dv$$

Subtraction results in

$$\sum_{n=1}^{\infty}\frac{\log(n+1)}{n^2}-\sum_{n=1}^{\infty}\frac{\log n}{n^2}=\int_0^1\frac{Li_2(v)[v-1]}{v\log v}dv$$

or

(E.35aiii) $$\sum_{n=1}^{\infty}\frac{\log\left(1+\frac{1}{n}\right)}{n^2}=\int_0^1\frac{Li_2(v)[v-1]}{v\log v}dv$$

Integrating (E.35a) we obtain

$$\sum_{k=1}^{N}(x+k)\log\left(1+\frac{x}{k}\right)-Nx=-\int_0^1\frac{x(1-v^N)}{(1-v)\log v}dv-\int_0^1\frac{(1-v^x)(1-v^N)}{(1-v)\log^2 v}dv$$

$$=-\int_0^1\frac{(1-v^x+x\log v)(1-v^N)}{(1-v)\log^2 v}dv$$

Differentiating (E.35a) results in

(E.35b) $$\sum_{k=1}^{N}\frac{1}{k+x}=\int_0^1\frac{(1-v^N)v^x}{1-v}dv \quad x=0\Rightarrow \quad \sum_{k=1}^{N}\frac{1}{k}=\int_0^1\frac{(1-v^N)}{1-v}dv$$

and this may be easily verified by using the geometric series. Further differentiations give us for $p\geq 0$

(E.35c) $$(-1)^p p!\sum_{k=1}^{N}\frac{1}{(k+x)^{p+1}}=\int_0^1\frac{(1-v^N)v^x\log^p v}{1-v}dv$$

and with $x = 0$ we get

(E.35d) $$(-1)^p p! H_N^{(p+1)} = \int_0^1 \frac{(1-v^N)\log^p v}{1-v} dv$$

In the limit as $N \to \infty$ we have $v^N \to 0$ and it therefore seems that (please note that I have not directly confirmed the validity of this operation)

(E.36) $$(-1)^p p!\left[\varsigma(p+1,x) - \frac{1}{x^{p+1}}\right] = \int_0^1 \frac{v^x \log^p v}{1-v} dv$$

With the substitution $v = e^{-t}$ in (E.36) we obtain

$$\int_0^1 \frac{v^x \log^p v}{1-v} dv = (-1)^p \int_0^\infty \frac{t^p e^{-xt}}{e^t - 1} dt$$

and using the well-known integral representation of the Hurwitz zeta function this becomes

$$= (-1)^p \Gamma(p+1)\varsigma(p+1, x+1)$$

$$= (-1)^p p!\left[\varsigma(p+1,x) - \frac{1}{x^{p+1}}\right]$$

and thus we have confirmed the validity of (E.36).

Making the substitution $v = e^{-t}$ in (E.35c) we obtain

(E.36a) $$\sum_{k=1}^N \frac{1}{(k+x)^{p+1}} = \frac{1}{p!}\int_0^\infty \frac{t^p e^{-xt}(1-e^{-Nt})}{e^t - 1} dt$$

and with $x = 0$ we get an integral representation for the generalised harmonic numbers

(E.36b) $$H_N^{(p+1)} = \frac{1}{p!}\int_0^\infty \frac{t^p e^{-xt}(1-e^{-Nt})}{e^t - 1} dt$$

In fact this is quite evident by noting that

$$= \frac{1}{p!}\int_0^\infty \frac{t^p e^{-xt}}{e^t - 1} dt - \frac{1}{p!}\int_0^\infty \frac{t^p e^{-(x+N)t}}{e^t - 1} dt$$

$$= \varsigma(p+1,x) - \varsigma(p+1,x+N)$$

In (E.36a) can we treat $p$ as a variable and differentiate accordingly?

From (E.35a) and (E.35b) we see that

$$\sum_{k=1}^{N}\left[\log\left(1+\frac{x}{k}\right)-\frac{x}{k}\right]=-\int_0^1\frac{(1-v^x)(1-v^N)}{(1-v)\log v}dv-x\int_0^1\frac{(1-v^N)}{1-v}dv$$

$$=-\int_0^1\frac{(1-v^x+x\log v)(1-v^N)}{(1-v)\log v}dv$$

and this formally suggests that as $N\to\infty$ we have

(E.37) $$\log\Gamma(x)+\log x+\gamma x=\int_0^1\frac{1-v^x+x\log v}{(1-v)\log v}dv=\sum_{k=2}^{\infty}\frac{(-1)^k}{k}x^k\varsigma(k)$$

With $x=1$ we get the well-known integral

(E.37i) $$\gamma=\int_0^1\frac{1-v+\log v}{(1-v)\log v}dv$$

Differentiating (E.37) gives us

(E.38) $$\psi(x)+\frac{1}{x}+\gamma=\int_0^1\frac{1-v^x}{1-v}dv$$

which is in agreement with (E.50).

Since the integral in (E.35a) is symmetrical in $x$ and $N$ we see that

$$\sum_{k=1}^{N}\log\left(1+\frac{M}{k}\right)=\sum_{k=1}^{M}\log\left(1+\frac{N}{k}\right)$$

and with $M=N+p$ we have

$$\sum_{k=1}^{N}\log\left(1+\frac{N+p}{k}\right)=\sum_{k=1}^{N+p}\log\left(1+\frac{N}{k}\right)$$

Referring to (E.37) we have

$$I(x)=\int_0^1\frac{1-v^x+x\log v}{(1-v)\log v}dv$$

and letting $v=t^2$ we obtain

$$I(x)=\int_0^1 \frac{(1-t^{2x}+2x)t}{(1-t^2)\log t}dt$$

With $x=1/2$ we have

$$I(1/2)=\int_0^1 \left[\frac{t}{(1+t)\log t}+\frac{2t}{1-t^2}\right]dt$$

The latter integral has similarities to (4.4.102) which is reproduced below

$$\log\frac{\pi}{2}=\int_0^1 \frac{t-1}{(1+t)\log t}dt$$

and hence we obtain

$$I(1/2)-\log\frac{\pi}{2}=\int_0^1 \left[\frac{1}{(1+t)\log t}+\frac{2t}{1-t^2}\right]dt$$

Since $I(1/2)=\frac{1}{2}\log\pi-\log 2+\frac{1}{2}\gamma$ we see that

(E.39) $$\frac{1}{2}[\gamma-\log\pi]=\int_0^1 \left[\frac{1}{(1+t)\log t}+\frac{2t}{1-t^2}\right]dt$$

See also (4.3.66f) in Volume II(a) where we showed that

$$\sum_{k=1}^{\infty}\frac{(-1)^{k+1}}{k}H_{n-1}^{(k)}x^k=\log\left(1+\frac{x}{n}\right)+\int_0^1 \frac{(1-v^x)(1-v^n)}{(1-v)\log v}dv$$

Then using (E.35a) we obtain

$$\sum_{k=1}^{\infty}\frac{(-1)^{k+1}}{k}H_{n-1}^{(k)}x^k=-\sum_{k=1}^{n-1}\log\left(1+\frac{x}{k}\right)$$

which is equivalent to (E.33r).

As Glicksman [70a] showed in 1943, we may retrace our steps to obtain a familiar expression for the case where $x=1$. In that case we have

$$\log\left(1+\frac{1}{k}\right)=\sum_{n=1}^{\infty}(-1)^{n+1}\frac{1}{nk^n}=\frac{1}{k}-\frac{1}{2k^2}+\frac{1}{3k^3}-...$$

and, as before, construct the finite sum

$$\sum_{k=1}^{N-1}\log\left(1+\frac{1}{k}\right)=\sum_{k=1}^{N-1}\frac{1}{k}-\frac{1}{2}\sum_{k=1}^{N-1}\frac{1}{k^2}+\frac{1}{3}\sum_{k=1}^{N-1}\frac{1}{k^3}-...$$

Alternatively we have

$$\sum_{k=1}^{N-1}\log\left(1+\frac{1}{k}\right)=\log\left(\frac{2}{1}\right)\left(\frac{3}{2}\right)\left(\frac{4}{3}\right)...\left(\frac{N}{N-1}\right)=\log N$$

Therefore we get

$$\frac{1}{2}\sum_{k=1}^{N-1}\frac{1}{k^2}-\frac{1}{3}\sum_{k=1}^{N-1}\frac{1}{k^3}+...\pm\frac{1}{N-1}\sum_{k=1}^{N-1}\frac{1}{k^{N-1}}=\sum_{k=1}^{N-1}\frac{1}{k}-\log N$$

$$=-\frac{1}{N}+\left[\sum_{k=1}^{N}\frac{1}{k}-\log N\right]$$

In the limit as $N\to\infty$ we have the well-known result [25, p.202]

(E.40) $$\gamma=\sum_{k=2}^{\infty}(-1)^k\frac{\varsigma(k)}{k}$$

which may also be obtained by letting $x=1$ in (E.22n).

Alternatively, we could slightly generalise the analysis given by Hwang Chien-Lih in [43b]. We have

$$\frac{x}{1+x}=\sum_{n=2}^{\infty}(-1)^n x^{n-1} \quad \text{, for } |x|<1$$

Integration gives us

$$I_k=\int_0^{u/k}\frac{x}{1+x}\,dx=\int_0^{u/k}\sum_{n=2}^{\infty}(-1)^n x^{n-1}\,dx$$

$$=\sum_{n=2}^{\infty}\frac{(-1)^n}{n}\frac{u^n}{k^n}$$

We then have

$$\sum_{k=1}^{\infty}I_k=\sum_{k=1}^{\infty}\sum_{n=2}^{\infty}\frac{(-1)^n}{n}\frac{u^n}{k^n}$$

$$=\sum_{n=2}^{\infty}\frac{(-1)^n u^n}{n}\sum_{k=1}^{\infty}\frac{1}{k^n}$$

$$= \sum_{n=2}^{\infty} \frac{(-1)^n}{n} \varsigma(n) u^n$$

On the other hand we have

$$I_k = \int_0^{u/k} \frac{x}{1+x} dx = \int_0^{u/k} \left(1 - \frac{1}{1+x}\right) dx$$

$$= \frac{u}{k} - \log\left(1 + \frac{u}{k}\right)$$

Upon summation we have

$$\sum_{k=1}^{\infty} \int_0^{u/k} \frac{x}{1+x} dx = \lim_{N\to\infty} \sum_{k=1}^{N-1} \left[\frac{u}{k} - \log\left(1 + \frac{u}{k}\right)\right]$$

and we see from (E.13a) that

$$\lim_{N\to\infty} \sum_{k=1}^{N-1} \left[\frac{u}{k} - \log\left(1 + \frac{u}{k}\right)\right] = \log u + \gamma u + \log \Gamma(u)$$

Hence we obtain another proof of (E.34a) for $u \in [0,1]$

(E.41) $$\log u + \gamma u + \log \Gamma(u) = \sum_{n=2}^{\infty} \frac{(-1)^n}{n} u^n \varsigma(n)$$

and this concurs with (E.40) when $u = 1$.

An integration of (E.41) gives us for $0 \le x \le 1$

(E.42) $$\frac{1}{2}\gamma x^2 + x \log x - x + \int_0^x \log \Gamma(u) du = \sum_{n=2}^{\infty} \frac{(-1)^n}{n(n+1)} \varsigma(n) x^{n+1}$$

and with $x = 1$ we get [126, p.223]

(E.42a) $$\sum_{n=2}^{\infty} \frac{(-1)^n}{n(n+1)} \varsigma(n) = \frac{1}{2}\gamma - 1 + \frac{1}{2}\log(2\pi)$$

Dividing (E.41) by $u$ and integrating results in for $0 < x \le 1$

(E.43) $$\frac{1}{2}\log^2 x + \frac{1}{2}\gamma(x^2 - 1) + \int_1^x \frac{\log \Gamma(u)}{u} du = \sum_{n=2}^{\infty} \frac{(-1)^n}{n^2} [x^n - 1] \varsigma(n)$$

We also have

(E.43a) $$\frac{1}{2}\gamma x^2+\int_0^x \frac{\log\Gamma(1+u)}{u}du=\sum_{n=2}^{\infty}\frac{(-1)^n}{n^2}\varsigma(n)x^n$$

From (E.34a) we saw that

$$\log\Gamma(x)+\log x+\gamma x=\sum_{k=2}^{\infty}\frac{(-1)^k}{k}\varsigma(k)x^k$$

which may be written as

$$\log\Gamma(1+x)+\gamma x=\sum_{k=2}^{\infty}\frac{(-1)^k}{k}\varsigma(k)\,x^k$$

and letting $x\to -x$ we get

(E.43ai) $$\log\Gamma(1-x)-\gamma x=\sum_{k=2}^{\infty}\frac{x^k}{k}\varsigma(k)$$

Hence we obtain upon addition

$$\log\Gamma(x)+\log\Gamma(1-x)+\log x=\sum_{k=2}^{\infty}\frac{[1+(-1)^k]}{k}\varsigma(k)x^k$$

and Euler's reflection formula then gives us (see also (6.138) in Volume V)

$$\log\frac{\pi x}{\sin\pi x}=\sum_{k=1}^{\infty}\frac{\varsigma(2k)}{k}x^{2k}$$

Kanemitsu et al. [82b] have shown that for $|x|<\operatorname{Re}(a)$

(E.43aii) $$\sum_{k=2}^{\infty}\frac{x^k}{k}\varsigma(k,a)=\varsigma'(0,a-x)-\varsigma'(0,a)+\psi(a)x$$

and with $a=1$ we see that for $|x|<1$

(E.43aiii) $$\sum_{k=2}^{\infty}\frac{x^k}{k}\varsigma(k)=\varsigma'(0,1-x)-\varsigma'(0)-\gamma x$$

With Lerch's identity $\log\Gamma(1-x)=\varsigma'(0,1-x)-\varsigma'(0)$ we then recover (E.43ai).

More generally we see that (see also (6.99a) in Volume V)

(E.43aiv) $$\sum_{k=2}^{\infty} \frac{x^k}{k} \varsigma(k,a) = \log \Gamma(a-x) - \log \Gamma(a) + \psi(a)x$$

and with $x = -1, a = 1$ we obtain

(E.43av) $$\sum_{k=2}^{\infty} \frac{(-1)^k}{k} \varsigma(k) = \gamma$$

Further identities may be obtained by differentiating (E.43aiv) with respect to both $x$ and $a$.

THE GENERALISED EULER CONSTANT FUNCTION $\gamma(x)$

In 2006 Sondow and Hadjicostas [123ab] considered the generalised Euler constant function $\gamma(x)$ which they defined as

(E.43b) $$\gamma(x) = \sum_{n=1}^{\infty} x^{n-1} \left[ \frac{1}{n} - \log\left(1 + \frac{1}{n}\right) \right]$$

and hence $\gamma(1) = \gamma$ and using (E.6i) we have $\gamma(-1) = \log(\pi/4)$.

They showed that

$$x\gamma(x) = x(1 - \log 2) + \sum_{n=2}^{\infty} x^n \left[ \frac{1}{n} - \log\left(1 + \frac{1}{n}\right) \right]$$

and using

$$\sum_{n=2}^{\infty} x^n \left[ \frac{1}{n} - \log\left(1 + \frac{1}{n}\right) \right] = \sum_{n=2}^{\infty} x^n \left[ \frac{1}{n} + \sum_{k=1}^{\infty} \frac{(-1)^k}{kn^k} \right]$$

$$= \sum_{n=2}^{\infty} x^n \sum_{k=2}^{\infty} \frac{(-1)^k}{kn^k}$$

they obtained

(E.43c) $$x\gamma(x) = \sum_{n=2}^{\infty} \frac{(-1)^n}{n} Li_n(x)$$

With $x = 1$ we obtain (E.43av) above.

Coincidentally, I carried out the following analysis a few months before [123ab] was published. We have from (4.4.38e) in Volume II(b)

$$Li_n(u) = \frac{(-1)^{n-1}}{(n-2)!}\int_0^1 \frac{\log^{n-2} y \log(1-uy)}{y} dy$$

With reference to (E.43c) and, using (4.4.38e), we have

$$\sum_{n=2}^{\infty} \frac{(-1)^n}{n} Li_n(u) = -\sum_{n=2}^{\infty} \int_0^1 \frac{\log^{n-2} y}{n(n-2)!} \frac{\log(1-uy)}{y} dy$$

$$= -\int_0^1 \sum_{n=2}^{\infty} \frac{\log^{n-2} y}{n(n-2)!} \frac{\log(1-uy)}{y} dy$$

The exponential series gives us

$$f(x) = \frac{e^x}{x} = \frac{1}{x} + 1 + \sum_{n=2}^{\infty} \frac{x^{n-1}}{n!}$$

and differentiation results in

$$f'(x) = -\frac{1}{x^2} + \sum_{n=2}^{\infty} \frac{x^{n-2}}{n(n-2)!}$$

We therefore have

(E.43ci) $$\sum_{n=2}^{\infty} \frac{x^{n-2}}{n(n-2)!} = \frac{e^x(x-1)}{x^2} + \frac{1}{x^2}$$

Hence we obtain

$$\sum_{n=2}^{\infty} \frac{\log^{n-2} y}{n(n-2)!} = \frac{y(\log y - 1)}{\log^2 y} + \frac{1}{\log^2 y}$$

Accordingly we have

(E.43d) $$\sum_{n=2}^{\infty} \frac{(-1)^n}{n} Li_n(u) = -\int_0^1 \frac{[y(\log y - 1) + 1]\log(1-uy)}{y \log^2 y} dy$$

and hence we get

(E.43e) $$u\gamma(u) = -\int_0^1 \frac{[y(\log y - 1) + 1]\log(1-uy)}{y \log^2 y} dy$$

Letting $u = 1$ results in

(E.43f) $$\gamma(1)=\gamma=-\int_0^1 \frac{[y(\log y-1)+1]\log(1-y)}{y\log^2 y}dy$$

and letting $u=-1$ gives us

(E.43g) $$\gamma(-1)=\log\frac{4}{\pi}=\int_0^1 \frac{[y(\log y-1)+1]\log(1+y)}{y\log^2 y}dy$$

Alternatively, with reference to (E.43c) and, using (4.4.38b) from Volume II(b), we have

$$\sum_{n=2}^{\infty}\frac{(-1)^n}{n}Li_n(u)=-u\sum_{n=2}^{\infty}\int_0^1 \frac{\log^n y}{n!(1-uy)\log y}dy$$

$$=u\int_0^1 \frac{1-y+\log y}{(1-uy)\log y}dy$$

$$=u\int_0^1 \frac{1-y}{(1-uy)\log y}dy-\log(1-u)$$

We therefore obtain from (E.43c)

(E.43h) $$u\gamma(u)=u\int_0^1 \frac{1-y+\log y}{(1-uy)\log y}dy$$

which was also derived by Sondow and Hadjicostas [123ab, equation (11)] by a different method. See also (4.4.112b) in Volume III.

Letting $u=1$ results in

(E.43i) $$\sum_{n=2}^{\infty}\frac{(-1)^n}{n}\varsigma(n)=\int_0^1 \frac{1-y+\log y}{(1-y)\log y}dy=\gamma$$

and letting $u=-1$ results in

(E.43j) $$\log\frac{4}{\pi}=\int_0^1 \frac{1-y+\log y}{(1+y)\log y}dy$$

$$=\int_0^1 \frac{1-y}{(1+y)\log y}dy+\log 2$$

Hence we obtain (4.4.100u) in Volume III

$$\log\frac{\pi}{2}=\int_0^1 \frac{y-1}{(1+y)\log y}dy$$

We also note that

(E.43k) $$\int_a^b \frac{1-y+\log y}{(1-y)\log y}dy=\int_a^b \frac{1}{\log y}dy+\int_a^b \frac{1}{(1-y)}dy$$

$$=li(b)-li(a)+\log\frac{1-a}{1-b}$$

where $li(x)=\int_0^x \frac{dt}{\log t}$ for $x\in[0,1)$ is the logarithmic integral. This then indicates that

$$\lim_{b\to 1}[li(b)-\log(1-b)]=\gamma$$

Having regard to (E.22bvi)

$$li(x)=\gamma+\log(-\log x)+\sum_{n=1}^{\infty}\frac{\log^n x}{n!n}$$

we see that

$$\lim_{b\to 1}[li(b)-\log(1-b)]=\lim_{b\to 1}[\gamma+\log(-\log b)-\log(1-b)]=\gamma$$

since by L'Hopital's rule we have $\lim_{b\to 1}\log\left[\frac{-\log b}{1-b}\right]=1$.

We now consider the integral

$$\int_0^b \frac{1-y+\log y}{(1+y)\log y}dy=\int_0^b\left[\frac{2}{1+y}-1\right]\frac{dy}{\log y}+\int_0^b \frac{1}{1+y}dy$$

and as $b\to 1$ we have

$$\log\frac{4}{\pi}=\lim_{b\to 1}\left[\int_0^b \frac{2dy}{(1+y)\log y}-li(b)\right]+\log 2$$

This may be written as

$$\log\frac{2}{\pi}=\lim_{b\to 1}\left[\int_0^b \frac{2dy}{(1+y)\log y}-\gamma-\log(-\log b)\right]$$

or as

$$\log\frac{2}{\pi}=\lim_{b\to 1}\left[\int_0^b \frac{2}{(1+y)\log y}-\frac{1}{\log y}\right]dy$$

$$=-\int_0^1 \frac{y-1}{(1+y)\log y}dy$$

We have seen that

$$\log\frac{2}{\pi}=\lim_{b\to 1}\left[\int_0^b \frac{2dy}{(1+y)\log y}-\gamma-\log(-\log b)\right]$$

$$=\lim_{b\to 1}\left[\int_0^b \frac{2dy}{(1+y)\log y}+\int_0^b \log(-\log y)dy-\log(-\log b)\right]$$

$$=\lim_{b\to 1}\left[\int_0^b \left(\frac{2}{(1+y)\log y}+\log\frac{\log y}{\log b}\right)dy\right]$$

□

Equating (E.43h) and (E.43e) gives us

$$-\int_0^1 \frac{[y(\log y-1)+1]\log(1-uy)}{y\log^2 y}dy=u\int_0^1 \frac{1-y+\log y}{(1-uy)\log y}dy$$

Differentiating the above results in

$$\int_0^1 \frac{y(\log y-1)+1}{(1-uy)\log^2 y}dy=u\int_0^1 \frac{y[1-y+\log y]}{(1-uy)^2\log y}dy+\int_0^1 \frac{1-y+\log y}{(1-uy)\log y}dy$$

$$=\int_0^1 \frac{1-y+\log y}{(1-uy)^2\log y}dy$$

The above formula (E.43d) may be extended as follows: with reference to (E.43c) and, using (4.4.38e), we have

$$\sum_{n=2}^{\infty}\frac{(-1)^n}{n}t^{n-2}Li_n(u)=-\sum_{n=2}^{\infty}\int_0^1 \frac{t^{n-2}\log^{n-2}y}{n(n-2)!}\frac{\log(1-uy)}{y}dy$$

$$=-\int_0^1 \sum_{n=2}^{\infty}\frac{(t\log y)^{n-2}}{n(n-2)!}\frac{\log(1-uy)}{y}dy$$

and, continuing as before, we find

(E.43l) $$\sum_{n=2}^{\infty}\frac{(-1)^n}{n}t^n Li_n(u) = -\int_0^1 \frac{[y^t(t\log y-1)+1]\log(1-uy)}{y\log^2 y}dy$$

With $u=1$ this becomes

(E.43l) $$\sum_{n=2}^{\infty}\frac{(-1)^n}{n}\varsigma(n)t^n = -\int_0^1 \frac{[y^t(t\log y-1)+1]\log(1-y)}{y\log^2 y}dy$$

and from (E.22n) and [126, p.159] we have for $-1 < t \le 1$

$$\sum_{n=2}^{\infty}\frac{(-1)^n}{n}\varsigma(n)t^n = \log\Gamma(1+t)+\gamma t$$

Hence we have

(E.43m) $$-\int_0^1 \frac{[y^t(t\log y-1)+1]\log(1-y)}{y\log^2 y}dy = \log\Gamma(1+t)+\gamma t$$

Differentiating (E.43m) results in

(E.43n) $$-t\int_0^1 \frac{y^{t-1}\log(1-y)}{\log y}dy = \psi(1+t)+\gamma$$

which may also be obtained by differentiating the Beta function.
Multiplying (E.43ci) by $x$ and integrating results in

$$\sum_{n=2}^{\infty}\frac{v^{n-2}}{n^2(n-2)!} = \frac{1}{v^2}\int_0^v \left[\frac{e^x(x-1)}{x}+\frac{1}{x}\right]dx$$

and we see that

$$\sum_{n=2}^{\infty}\frac{(-1)^n}{n^2}t^{n-2}Li_n(u) = -\int_0^1 \sum_{n=2}^{\infty}\frac{(t\log y)^{n-2}}{n^2(n-2)!}\frac{\log(1-uy)}{y}dy$$

This gives us

(E.43ni $$\sum_{n=2}^{\infty}\frac{(-1)^n}{n^2}t^n Li_n(u) = -\int_0^1\int_0^{t\log y}\left[\frac{e^x(x-1)}{x}+\frac{1}{x}\right]\frac{\log(1-uy)}{y\log^2 y}dx\,dy$$

Integrating (E.43b) gives us

$$\int_0^u \gamma(x)dx = \sum_{n=1}^{\infty} \frac{u^n}{n}\left[\frac{1}{n} - \log\left(1+\frac{1}{n}\right)\right]$$

$$= Li_2(u) - \sum_{n=1}^{\infty} \frac{u^n}{n} \log\left(1+\frac{1}{n}\right)$$

Then using (E.43c) we obtain

$$\int_0^u \gamma(x)dx = \sum_{n=2}^{\infty} \frac{(-1)^n}{n} \int_0^u \frac{Li_n(x)}{x} dx$$

and thus we get

(E.43o) $$\int_0^u \gamma(x)dx = \sum_{n=2}^{\infty} \frac{(-1)^n}{n} Li_{n+1}(u)$$

Hence we have

(E.43p) $$\sum_{n=1}^{\infty} \frac{u^n}{n} \log\left(1+\frac{1}{n}\right) = -\sum_{n=1}^{\infty} \frac{(-1)^n}{n} Li_{n+1}(u)$$

Further integrations result in (after dividing both sides by $u$ )

(E.43q) $$\sum_{n=1}^{\infty} \frac{u^n}{n^p} \log\left(1+\frac{1}{n}\right) = -\sum_{n=1}^{\infty} \frac{(-1)^n}{n} Li_{n+p}(u)$$

and, with $u=1$, this becomes for $p \geq 2$

(E.43r) $$\varsigma'(p) + \sum_{n=1}^{\infty} \frac{\log(n+1)}{n^p} = -\sum_{n=1}^{\infty} \frac{(-1)^n}{n} \varsigma(n+p)$$

From (4.4.42i) in Volume II(b) we have

$$Li_{n+2}(u) = \frac{(-1)^{n+1}}{\Gamma(n)} \int_0^1 \frac{\log^{n-1} y \, Li_2(uy)}{y} dy$$

and therefore we see that

$$\sum_{n=1}^{\infty} \frac{u^n}{n^2} \log\left(1+\frac{1}{n}\right) = \sum_{n=1}^{\infty} \frac{1}{n!} \int_0^1 \frac{\log^{n-1} y \, Li_2(uy)}{y} dy$$

$$= \int_0^1 \sum_{n=1}^{\infty} \frac{\log^n y}{n!} \frac{Li_2(uy)}{y \log y} dy$$

$$= \int_0^1 (e^{\log y} - 1) \frac{Li_2(uy)}{y \log y} dy$$

Hence we have

$$\sum_{n=1}^{\infty} \frac{u^n}{n^2} \log\left(1 + \frac{1}{n}\right) = \int_0^1 \frac{(y-1) Li_2(uy)}{y \log y} dy$$

and with $u = 1$ we have

$$\varsigma'(2) + \sum_{n=1}^{\infty} \frac{\log(n+1)}{n^2} = \int_0^1 \frac{(y-1) Li_2(y)}{y \log y} dy$$

We have from (4.4.38e) in Volume II(a)

$$Li_{n+1}(u) = \frac{(-1)^n}{(n-1)!} \int_0^1 \frac{\log^{n-1} y \log(1-uy)}{y} dy$$

and therefore we get

$$\sum_{n=2}^{\infty} \frac{(-1)^n}{n} Li_{n+1}(u) = \sum_{n=2}^{\infty} \int_0^1 \frac{\log^n y}{n!} \frac{\log(1-uy)}{y \log y} dy$$

$$= \int_0^1 \sum_{n=2}^{\infty} \frac{\log^n y}{n!} \frac{\log(1-uy)}{y \log y} dy$$

$$= \int_0^1 \left[ e^{\log y} - 1 - \log y \right] \frac{\log(1-uy)}{y \log y} dy$$

$$= \int_0^1 \left[ y - 1 - \log y \right] \frac{\log(1-uy)}{y \log y} dy$$

$$= \int_0^1 \frac{(y-1) \log(1-uy)}{y \log y} dy - \int_0^1 \frac{\log(1-uy)}{y} dy$$

$$= \int_0^1 \frac{(y-1) \log(1-uy)}{y \log y} dy + Li_2(u)$$

Hence we see that

$$\sum_{n=1}^{\infty} \frac{(-1)^n}{n} Li_{n+1}(u) = \int_0^1 \frac{(y-1) \log(1-uy)}{y \log y} dy$$

It is possible that Adamchik's paper "A class of logarithmic integrals" [2a] may be germane to the evaluation of the above integral.

Integrating (E.43c) we get

$$\int_0^x u\gamma(u)du = \sum_{n=2}^{\infty} \frac{(-1)^n}{n} \int_0^x Li_n(u)du$$

By trial and error we see that for $n \geq 2$

$$\int_0^x Li_n(u)\,du = (-1)^n \left[(x-1)\log(1-x) - x\right] + (-1)^n x \sum_{j=2}^{n} (-1)^j Li_j(x)$$

and therefore we see that

$$\int_0^x u\gamma(u)du = \sum_{n=2}^{\infty} \frac{1}{n} \left( \left[(x-1)\log(1-x) - x\right] + x \sum_{j=2}^{n} (-1)^j Li_j(x) \right)$$

(the convergence of which appears to be questionable).

Since from (1.5b) in Volume I $Li_1(x) = -\log(1-x)$ we may write this as

$$\int_0^x u\gamma(u)du = \sum_{n=2}^{\infty} \frac{1}{n} \left( -\left[\log(1-x) + x\right] + x \sum_{j=1}^{n} (-1)^j Li_j(x) \right)$$

and differentiating this we get

$$x\gamma(x) = \sum_{n=2}^{\infty} \frac{1}{n} \left( \frac{1}{1-x} - 1 + \sum_{j=1}^{n} (-1)^j [Li_j(x) + Li_{j-1}(x)] \right)$$

This telescopes to

$$= \sum_{n=2}^{\infty} \frac{1}{n} \left( \frac{x}{1-x} - Li_0(x) + (-1)^n Li_n(x) \right)$$

which simplifies to (E.43c)

$$x\gamma(x) = \sum_{n=2}^{\infty} \frac{(-1)^n}{n} Li_n(x)$$

□

In (4.4.43) in Volume II(b) we saw that

$$\sum_{n=1}^{\infty} t^n \sum_{k=1}^{n} \binom{n}{k} \frac{x^k}{(k+y)^s} = \frac{1}{1-t} \sum_{n=0}^{\infty} \frac{1}{(n+y)^s} \left[ \frac{xt}{(1-t)} \right]^n - \frac{1}{(1-t)y^s}$$

$$= \frac{1}{1-t} \sum_{n=1}^{\infty} \frac{1}{(n+y)^s} \left[ \frac{xt}{(1-t)} \right]^n$$

and with $s = 1$ we get

$$\sum_{n=1}^{\infty} t^n \sum_{k=1}^{n} \binom{n}{k} \frac{x^k}{k+y} = \frac{1}{1-t} \sum_{n=1}^{\infty} \frac{1}{n+y} \left[ \frac{xt}{(1-t)} \right]^n$$

Integration results in

$$\sum_{n=1}^{\infty} t^n \sum_{k=1}^{n} \binom{n}{k} (-1)^k \log \frac{(k+v)}{k} = \frac{1}{1-t} \sum_{n=1}^{\infty} \log \frac{(n+v)}{n} \left[ \frac{-t}{1-t} \right]^n$$

Letting $v = 1$ gives us

$$\sum_{n=1}^{\infty} t^n \sum_{k=1}^{n} \binom{n}{k} (-1)^k \log \frac{(k+1)}{k} = \frac{1}{1-t} \sum_{n=1}^{\infty} \left[ \frac{-t}{1-t} \right]^n \log \left( 1 + \frac{1}{n} \right)$$

From (E.43c) we have

$$\gamma(x) = \sum_{n=1}^{\infty} x^{n-1} \left[ \frac{1}{n} - \log \left( 1 + \frac{1}{n} \right) \right]$$

and therefore we get

$$-x\gamma(-x) = \sum_{n=1}^{\infty} (-x)^n \left[ \frac{1}{n} - \log \left( 1 + \frac{1}{n} \right) \right]$$

$$= -\log(1+x) - \sum_{n=1}^{\infty} (-x)^n \log \left( 1 + \frac{1}{n} \right)$$

With $x = \dfrac{t}{1-t}$ we see that

$$-\frac{t}{1-t} \gamma \left( \frac{-t}{1-t} \right) = \log(1-t) - \sum_{n=1}^{\infty} \left[ \frac{-t}{1-t} \right]^n \log \left( 1 + \frac{1}{n} \right)$$

$$= \log(1-t) - (1-t) \sum_{n=1}^{\infty} t^n \sum_{k=1}^{n} \binom{n}{k} (-1)^k \log \frac{(k+1)}{k}$$

$$t\gamma\left(\frac{-t}{1-t}\right)=-(1-t)\log(1-t)+(1-t)^2\sum_{n=1}^{\infty}t^n\sum_{k=1}^{n}\binom{n}{k}(-1)^k\log\frac{(k+1)}{k}$$

But Sondow and Hadjicostas [123ab] have shown that

$$t\gamma\left(\frac{-t}{1-t}\right)=-(1-t)\log(1-t)+(1-t)\sum_{n=1}^{\infty}t^n\sum_{k=1}^{n}\binom{n}{k}(-1)^k\log(k+1)$$

and this therefore implies that

$$(1-t)^2\sum_{n=1}^{\infty}t^n\sum_{k=1}^{n}\binom{n}{k}(-1)^k\log\frac{(k+1)}{k}=(1-t)\sum_{n=1}^{\infty}t^n\sum_{k=1}^{n}\binom{n}{k}(-1)^k\log(k+1)$$

We therefore obtain

$$-t\sum_{n=1}^{\infty}t^n\sum_{k=1}^{n}\binom{n}{k}(-1)^k\log(k+1)=(1-t)\sum_{n=1}^{\infty}t^n\sum_{k=1}^{n}\binom{n}{k}(-1)^k\log k$$

and with $t=1/2$ we get

$$-\sum_{n=1}^{\infty}\frac{1}{2^{n+1}}\sum_{k=1}^{n}\binom{n}{k}(-1)^k\log(k+1)=\sum_{n=1}^{\infty}\frac{1}{2^{n+1}}\sum_{k=1}^{n}\binom{n}{k}(-1)^k\log k$$

□

The motivation for the following work came from Villarino's paper [133a] "Ramanujan's approximation to the $n$th partial sum of the harmonic series".

Since $\gamma=\lim_{n\to\infty}\left(H_n-\log n\right)$ and $\gamma=\lim_{n\to\infty}\left(H_n-\log[n+1]\right)$ it is clear that

$$\gamma=\lim_{n\to\infty}\left(H_n-\frac{1}{2}\log[n(n+1)]\right)$$

With $\varepsilon_n=H_n-\frac{1}{2}\log[n(n+1)]-\gamma$ we get

$$\varepsilon_{n-1}-\varepsilon_n=\int_0^1\frac{x^2}{n(n^2-x^2)}dx$$

and therefore

$$\varepsilon_n=(\varepsilon_n-\varepsilon_{n+1})+(\varepsilon_{n+1}-\varepsilon_{n+2})+...$$

$$= \sum_{k=n+1}^{\infty} \int_0^1 \frac{x^2}{k(k^2-x^2)} dx$$

We now consider the series (note that this $\varepsilon_n$ is different from the one used just before (E.7))

$$\sum_{n=1}^{\infty} \frac{\varepsilon_n}{n^s} = \sum_{n=1}^{\infty} \frac{1}{n^s} \left[ H_n - \frac{1}{2} \log[n(n+1)] - \gamma \right]$$

$$= \sum_{n=1}^{\infty} \frac{H_n}{n^s} - \frac{1}{2} \sum_{n=1}^{\infty} \frac{\log n}{n^s} - \frac{1}{2} \sum_{n=1}^{\infty} \frac{\log(n+1)}{n^s} - \gamma \sum_{n=1}^{\infty} \frac{1}{n^s}$$

$$= \sum_{n=1}^{\infty} \frac{H_n}{n^s} + \frac{1}{2} \varsigma'(s) - \frac{1}{2} \sum_{n=1}^{\infty} \frac{\log(n+1)}{n^s} - \gamma \varsigma(s)$$

Alternatively, we also have

$$\sum_{n=1}^{\infty} \frac{\varepsilon_n}{n^s} = \sum_{n=1}^{\infty} \frac{1}{n^s} \sum_{k=n+1}^{\infty} \int_0^1 \frac{x^2}{k(k^2-x^2)} dx$$

$$= \sum_{n=1}^{\infty} \frac{1}{n^s} \left[ \sum_{k=n}^{\infty} \int_0^1 \frac{x^2}{k(k^2-x^2)} dx - \int_0^1 \frac{x^2}{n(n^2-x^2)} dx \right]$$

$$= \sum_{n=1}^{\infty} \frac{1}{n^s} \sum_{k=n}^{\infty} \int_0^1 \frac{x^2}{k(k^2-x^2)} dx - \sum_{n=1}^{\infty} \frac{1}{n^s} \int_0^1 \frac{x^2}{n(n^2-x^2)} dx$$

$$= \sum_{n=1}^{\infty} \left[ \int_0^1 \frac{x^2}{n(n^2-x^2)} dx \right] \sum_{k=n}^{\infty} \frac{1}{k^s} - \sum_{n=1}^{\infty} \frac{1}{n^s} \int_0^1 \frac{x^2}{n(n^2-x^2)} dx$$

We have

$$\int_0^1 \frac{x^2}{n(n^2-x^2)} dx = -\frac{1}{n} + \tanh^{-1} \frac{1}{n} = -\frac{1}{n} + \frac{1}{2} \log \frac{1+1/n}{1-1/n}$$

and the series becomes

$$= \sum_{n=1}^{\infty} \left[ -\frac{1}{n} + \tanh^{-1} \frac{1}{n} \right] \sum_{k=n}^{\infty} \frac{1}{k^s} - \sum_{n=1}^{\infty} \frac{1}{n^s} \left[ -\frac{1}{n} + \tanh^{-1} \frac{1}{n} \right]$$

$$= \sum_{n=1}^{\infty} \left[ -\frac{1}{n} + \tanh^{-1} \frac{1}{n} \right] \left( \varsigma(s) - H_{n-1}^{(s)} \right) + \varsigma(s+1) - \sum_{n=1}^{\infty} \frac{1}{n^s} \left[ \tanh^{-1} \frac{1}{n} \right]$$

Therefore we get

$$\sum_{n=1}^{\infty}\frac{H_n}{n^s}+\frac{1}{2}\varsigma'(s)-\frac{1}{2}\sum_{n=1}^{\infty}\frac{\log(n+1)}{n^s}-\gamma\varsigma(s)=$$

$$\sum_{n=1}^{\infty}\left[-\frac{1}{n}+\tanh^{-1}\frac{1}{n}\right]\left(\varsigma(s)-H_{n-1}^{(s)}\right)+\varsigma(s+1)-\sum_{n=1}^{\infty}\frac{1}{n^s}\left[\tanh^{-1}\frac{1}{n}\right]$$

and hence we have

(E.43b) $$\sum_{n=1}^{\infty}\frac{H_n}{n^s}+\frac{1}{2}\varsigma'(s)-\frac{1}{2}\sum_{n=1}^{\infty}\frac{\log(n+1)}{n^s}-\gamma\varsigma(s)=\sum_{n=1}^{\infty}\left[-\frac{1}{n}+\tanh^{-1}\frac{1}{n}\right]\left(\varsigma(s)-H_n^{(s)}\right)$$

KUMMER'S FOURIER SERIES FOR $\log\Gamma(x)$

In 1985 Berndt [20] gave an elementary proof of the Fourier series expansion for $\log\Gamma(x)$ which, inter alia, we referred to in (4.4.210)

(E.44) $$\log\Gamma(x)=\frac{1}{2}\log\pi-\frac{1}{2}\log\sin\pi x+\sum_{n=1}^{\infty}\frac{(\gamma+\log 2\pi n)\sin 2\pi nx}{\pi n}\qquad(0<x<1)$$

(this formula was originally derived by Kummer in 1847 [94]). Reference to (7.8) confirms that (4.4.210) is properly described as a Fourier series expansion for $\log\Gamma(x)$. Using (7.5) we may also write Kummer's formula (4.4.210a) as (cf. Nielsen [104a, p.79])

(E.44a) $$\log\Gamma(x)=\frac{1}{2}\log\frac{\pi}{\sin\pi x}+\left(\frac{1}{2}-x\right)(\gamma+\log 2\pi)+\frac{1}{\pi}\sum_{n=1}^{\infty}\frac{\log n}{n}\sin 2\pi nx$$

Berndt's proof employed Hurwitz's formula for the Fourier series expansion of the Hurwitz zeta function $\varsigma(s,x)$, but the following proof is based on Kummer's method as outlined in [8a, p.29].

From (7.5) and (7.8) of Volume V we have

$$\frac{\pi}{2}(1-2x)=\sum_{n=1}^{\infty}\frac{\sin 2n\pi x}{n}$$

$$-\log\left[2\sin\pi x\right]=\sum_{n=1}^{\infty}\frac{\cos 2n\pi x}{n}$$

Since $\log\Gamma(x)$ is differentiable in (0,1), it has a Fourier expansion

$$\log\Gamma(x)=C_0+2\sum_{k=1}^{\infty}C_k\cos 2k\pi x+2\sum_{k=1}^{\infty}D_k\sin 2k\pi x$$

where $C_k = \int_0^1 \log\Gamma(x)\cos 2k\pi x\,dx$ and $D_k = \int_0^1 \log\Gamma(x)\sin 2k\pi x\,dx$

The $C_k$ coefficients are relatively easy to determine using Euler's reflection formula for the gamma function

$$\log\Gamma(x) + \log\Gamma(1-x) = \log 2\pi - \log\left[2\sin\pi x\right]$$

$$= \log 2\pi + \sum_{n=1}^{\infty} \frac{\cos 2n\pi x}{n}$$

The corresponding Fourier series is

$$\log\Gamma(x) + \log\Gamma(1-x) = 2C_0 + 4\sum_{k=1}^{\infty} C_k \cos 2k\pi x$$

and equating coefficients in the last two equations gives

$$C_0 = \frac{1}{2}\log 2\pi \quad \text{and} \quad C_k = \frac{1}{4k} \quad \text{for } k \ge 1$$

We now refer to Malmstén's formula [126, p.16] which was derived in (E.22g)

(E.45) $$\log\Gamma(x) = \int_0^{\infty} \left[ (x-1)e^{-\alpha} - \frac{e^{-\alpha} - e^{-\alpha x}}{1-e^{-\alpha}} \right] \frac{d\alpha}{\alpha}$$

and a change of variables $u = e^{-\alpha}$ gives us

(E.45a) $$\log\Gamma(x) = \int_0^1 \left[ \frac{1-u^{x-1}}{1-u} - x + 1 \right] \frac{du}{\log u}$$

We therefore obtain

$$D_k = \int_0^1 \log\Gamma(x)\sin 2k\pi x\,dx$$

$$= \int_0^1\int_0^1 \left[ \frac{1-u^{x-1}}{1-u} - x + 1 \right] \frac{\sin 2k\pi x\,du dx}{\log u}$$

We have

$$\int_0^1 \sin 2k\pi x\, dx = 0$$

$$\int_0^1 x \sin 2k\pi x\, dx = -\frac{1}{2k\pi}$$

$$\int_0^1 u^{x-1} \sin 2k\pi x\, dx = \frac{1}{u} \operatorname{Im} \int_0^1 \exp\left[x(\log u + 2k\pi i)\right] dx$$

$$= \frac{1}{u} \operatorname{Im} \frac{u-1}{\log u + 2k\pi i}$$

Therefore we get

$$D_k = \int_0^1 \left[ \frac{-2k\pi}{u\left(\log^2 u + 4k^2\pi^2\right)} + \frac{1}{2k\pi} \right] \frac{du}{\log u}$$

With $u = e^{-2k\pi t}$ we obtain

$$D_k = \frac{1}{2\pi k} \int_0^\infty \left[ \frac{1}{\left(1+t^2\right)} - e^{-2k\pi t} \right] \frac{dt}{t}$$

Taking $k = 1$ we have

$$D_1 = \frac{1}{2\pi} \int_0^\infty \left[ \frac{1}{\left(1+t^2\right)} - e^{-2\pi t} \right] \frac{dt}{t}$$

With $x = 1$ in Malmstén's formula (E.45) we get

$$-\frac{\gamma}{2\pi} = \frac{1}{2\pi} \int_0^\infty \left[ e^{-t} - \frac{1}{1+t} \right] \frac{dt}{t}$$

Hence we have

$$D_1 - \frac{\gamma}{2\pi} = \frac{1}{2\pi} \int_0^\infty \frac{e^{-t} - e^{-2\pi t}}{t} dt + \frac{1}{2\pi} \int_0^\infty \left[ \frac{1}{\left(1+t^2\right)} - \frac{1}{(1+t)} \right] \frac{dt}{t}$$

By (E.22a) the first integral is equal to $\log 2\pi$ and a change of variables $t \to 1/t$ shows that the second integral vanishes. Therefore we have

$$D_1 = \frac{\gamma}{2\pi} + \frac{1}{2\pi}\log 2\pi$$

$D_k$ is found by noting that

$$kD_k - D_1 = \frac{1}{2\pi}\int_0^\infty \frac{e^{-2\pi t} - e^{-2k\pi t}}{t}\,dt = \frac{1}{2\pi}\log k$$

where the integral is also evaluated with (E.22a). Thus we have

(E.46) $$D_k = \frac{1}{2\pi k}\left(\gamma + \log 2k\pi\right) = \int_0^1 \log\Gamma(x)\sin 2k\pi x\,dx$$

and Kummer's formula thereby follows.

Kummer's Fourier series expansion (4.4.210) may also be written as

$$\log\Gamma(x) = -\frac{1}{2}\log x - \frac{1}{2}\log\frac{\sin\pi x}{\pi x} + \sum_{n=1}^\infty \frac{(\gamma + \log 2\pi n)\sin 2\pi nx}{\pi n} \qquad (0 < x < 1)$$

and hence we have

$$\frac{1}{2}\log\frac{\Gamma(x)}{\Gamma(1-x)} = \sum_{n=1}^\infty \frac{(\gamma + \log 2\pi n)\sin 2\pi nx}{\pi n}$$

We have $D_k = \int_0^1 \log\Gamma(x)\sin 2k\pi x\,dx = \frac{1}{2\pi k}\left(\gamma + \log 2k\pi\right)$ for $k \ge 1$ and hence we get

$$\sum_{k=1}^\infty \frac{D_k}{k} = \sum_{k=1}^\infty \frac{1}{2\pi k^2}\left(\gamma + \log 2k\pi\right)$$

$$= \frac{\left(\gamma + \log 2\pi\right)}{2\pi}\varsigma(2) - \frac{1}{2\pi}\varsigma'(2)$$

We also have

$$\sum_{k=1}^\infty \frac{D_k}{k} = \sum_{k=1}^\infty \frac{1}{k}\int_0^1 \log\Gamma(x)\sin 2k\pi x\,dx = \int_0^1 \log\Gamma(x)\sum_{k=1}^\infty \frac{\sin 2k\pi x}{k}\,dx$$

Now using (7.5) we get which is valid for $x \in (0,1)$

$$\sum_{k=1}^\infty \frac{\sin 2k\pi x}{k} = \frac{\pi}{2}(1-2x)$$

and, accordingly, since $\log\Gamma(0) = \log\Gamma(1) = 0$

$$\sum_{k=1}^{\infty} \log\Gamma(x) \frac{\sin 2k\pi x}{k} = \frac{\pi}{2} \log\Gamma(x)(1-2x)$$

is valid for $x \in [0,1]$. It may be noted that

$$\lim_{x\to 0} \left[\sin 2k\pi x \log\Gamma(x)\right] = \lim_{x\to 0} \left[\sin 2k\pi x \log\Gamma(x+1) - \sin 2k\pi x \log x\right] = 0$$

and hence the Fourier series expansion (when multiplied by $\log\Gamma(x)$) may be validly used at $x = 0$.

Hence we have

$$\sum_{k=1}^{\infty} \frac{D_k}{k} = \frac{\pi}{2} \int_0^1 \log\Gamma(x)(1-2x)\, dx$$

We note from (C.43b) that $\int_0^1 \log\Gamma(x)\, dx = \frac{1}{2} \log 2\pi$ and we obtain

$$\frac{\pi}{4} \log 2\pi - \pi \int_0^1 x \log\Gamma(x)\, dx = \frac{(\gamma + \log 2\pi)}{2\pi} \varsigma(2) - \frac{1}{2\pi} \varsigma'(2)$$

Therefore we get a further proof of (4.4.213d), namely

(E.47) $$\int_0^1 x \log\Gamma(x)\, dx = \frac{1}{6} \log 2\pi - \frac{\gamma}{12} + \frac{1}{2\pi^2} \varsigma'(2)$$

With reference to series of the form $\sum_{k=1}^{\infty} \frac{D_k}{k^p}$ and $\sum_{k=1}^{\infty} (-1)^k \frac{D_k}{k^p}$ we may easily obtain further identities using the Fourier series reported, for example, in [130, p.148].

We have $C_k = \int_0^1 \log\Gamma(x) \cos 2k\pi x\, dx = \frac{1}{4k}$ for $k \geq 1$ and hence we get

$$\sum_{k=1}^{\infty} (-1)^{k+1} \frac{C_k}{k} = \frac{1}{4} \varsigma_a(2)$$

We also have

$$\sum_{k=1}^{\infty} (-1)^{k+1} \frac{C_k}{k} = \sum_{k=1}^{\infty} \frac{1}{k} \int_0^1 \log\Gamma(x)(-1)^{k+1} \cos 2k\pi x\, dx = \int_0^1 \log\Gamma(x) \sum_{k=1}^{\infty} (-1)^{k+1} \frac{\cos 2k\pi x}{k} dx$$

Now using [130, p.148] we get

$$\text{(E.48)} \qquad \sum_{k=1}^{\infty}(-1)^{k+1}\frac{\cos 2k\pi x}{k} = \log[2\cos \pi x]$$

and accordingly

$$\text{(E.49)} \qquad \sum_{k=1}^{\infty}(-1)^{k+1}\frac{C_k}{k} = \int_0^1 \log\Gamma(x)\log[2\cos \pi x]\,dx = \frac{1}{4}\varsigma_a(2)$$

We therefore have

$$\text{(E.49a)} \qquad \int_0^1 \log\Gamma(x)\log\cos \pi x\,dx = -\frac{1}{2}\log 2\log 2\pi + \frac{\pi^2}{48}$$

It is noted in [59, p.18] that

$$\text{(E.49b)} \qquad \int_0^1 \log\Gamma(x)\log\sin \pi x\,dx = -\frac{1}{2}\log 2\log 2\pi - \frac{\pi^2}{24}$$

and addition of the two integrals results in

$$\text{(E.49c)} \qquad \int_0^1 \log\Gamma(x)\log\sin 2\pi x\,dx = -\frac{1}{2}\log 2\log 2\pi - \frac{\pi^2}{48}$$

From (6.92a) we have

$$\int_0^1 \log\Gamma(x+1)\log\left[2\sin(\pi x)\right]dx = \frac{1}{2\pi}\sum_{n=1}^{\infty}\frac{si(2n\pi)}{n^2}$$

and we see that

$$\int_0^1 \log\Gamma(x+1)\log\left[2\sin(\pi x)\right]dx = \int_0^1 \left[\log x + \log\Gamma(x)\right]\left[\log 2 + \log\sin(\pi x)\right]dx$$

$$= \log 2\int_0^1 \log x\,dx + \int_0^1 \log x\log\sin(\pi x)dx + \log 2\int_0^1 \log\Gamma(x)dx + \int_0^1 \log\Gamma(x)\log\sin(\pi x)\,dx$$

$$= -\log 2 + \int_0^1 \log x\log\sin(\pi x)dx + \frac{1}{2}\log 2\log 2\pi - \frac{1}{2}\log 2\log 2\pi - \frac{\pi^2}{24}$$

Hence we have

(E.49c) $$\int_0^1 \log x \log\sin(\pi x)dx = \frac{1}{2\pi}\sum_{n=1}^{\infty}\frac{si(2n\pi)}{n^2} + \log 2 + \frac{\pi^2}{24}$$

and $\sum_{n=1}^{\infty}\frac{si(2n\pi)}{n^2}$ may be determined from (6.117j) in Volume V.

In [20] Berndt used Kummer's formula (4.4.210) to derive Euler's reflection formula in the following elementary manner. We have

$$\log\Gamma(x) = \frac{1}{2}\log\pi - \frac{1}{2}\log\sin\pi x + \sum_{n=1}^{\infty}\frac{(\gamma+\log 2\pi n)\sin 2\pi nx}{\pi n} \quad (0 < x < 1)$$

Letting $x \to 1-x$ we get

$$\log\Gamma(1-x) = \frac{1}{2}\log\pi - \frac{1}{2}\log\sin\pi(1-x) + \sum_{n=1}^{\infty}\frac{(\gamma+\log 2\pi n)\sin 2\pi n(1-x)}{\pi n}$$

$$= \frac{1}{2}\log\pi - \frac{1}{2}\log\sin\pi x - \sum_{n=1}^{\infty}\frac{(\gamma+\log 2\pi n)\sin 2\pi nx}{\pi n}$$

and hence

$$\log\Gamma(x) + \log\Gamma(1-x) = \log\frac{\pi}{\sin\pi x}$$

We then have

$$\Gamma(x)\Gamma(1-x) = \frac{\pi}{\sin\pi x}$$

And that's all there is to it!

Alternatively, we may subtract the above two equations to obtain the following formula which appears in Ramanujan's Notebooks [Berndt, Vol.1, p.199]

$$\log\frac{\Gamma(x)}{\Gamma(1-x)} + [\gamma + \log(2\pi)](2x-1) = \frac{2}{\pi}\sum_{n=1}^{\infty}\frac{\log n}{n}\sin 2\pi nx$$

where we have used (7.5) $\sum_{n=1}^{\infty}\frac{\sin 2\pi nx}{n\pi} = \frac{1}{2} - x$ which is valid for $0 < x < 1$. Integration then results in

$$\int_0^u \log\Gamma(x)dx - \int_0^u \log\Gamma(1-x)dx + [\gamma + \log(2\pi)](u^2 - u) =$$

$$\frac{1}{\pi^2}\sum_{n=1}^{\infty}\frac{\log n}{n^2}-\frac{1}{\pi^2}\sum_{n=1}^{\infty}\frac{\log n}{n^2}\cos 2\pi nu$$

We have

$$\int_0^u \log\Gamma(1-x)dx=-\int_0^{-u}\log\Gamma(1+x)dx$$

and we noted in (6.112) that

$$\int_0^u \log\Gamma\left(1+x\right)dx=\frac{1}{2}\left[\log(2\pi)-1\right]u-\frac{u^2}{2}+u\log\Gamma\left(1+u\right)-\log G(1+u)$$

Hence we have

$$\int_0^u \log\Gamma(x)dx-\int_0^u \log\Gamma(1-x)dx=\int_0^u \log\Gamma(x)dx+\int_0^{-u}\log\Gamma(1+x)dx$$

$$=\int_0^u \log\Gamma(1+x)dx+\int_0^{-u}\log\Gamma(1+x)dx-\int_0^u \log x\,dx$$

We see that

$$\int_0^u \log\Gamma(1+x)dx+\int_0^{-u}\log\Gamma(1+x)dx=$$

$$-u^2+u\log\Gamma\left(1+u\right)-u\log\Gamma\left(1-u\right)-\log G(1+u)-\log G(1-u)$$

and with $u=1/2$ we get

$$\int_0^{1/2}\log\Gamma(1+x)dx+\int_0^{-1/2}\log\Gamma(1+x)dx=$$

$$-\frac{1}{4}+\frac{1}{2}\log\Gamma\left(3/2\right)-\frac{1}{2}\log\Gamma\left(1/2\right)-\log G(3/2)-\log G(1/2)$$

$$=-\frac{1}{4}-\frac{1}{2}\log 2-\log\Gamma\left(1/2\right)-2\log G(1/2)$$

We note from (6.127) that

$$\log G(1/2)=-\frac{3}{2}\log A-\frac{1}{4}\log\pi+\frac{1}{8}+\frac{1}{24}\log 2$$

and hence we get

$$\int_0^{1/2} \log\Gamma(1+x)dx + \int_0^{-1/2} \log\Gamma(1+x)dx = 3\log A - \frac{7}{12}\log 2$$

Therefore we obtain

$$3\log A - \frac{1}{12}\log 2 - \frac{1}{4}[\gamma + \log(2\pi)] = \frac{1}{\pi^2}\sum_{n=1}^{\infty}\frac{\log n}{n^2} - \frac{1}{\pi^2}\sum_{n=1}^{\infty}(-1)^n\frac{\log n}{n^2}$$

$$= -\frac{1}{\pi^2}[\varsigma'(2) + \varsigma_a'(2)]$$

Since $\varsigma_a'(s) = \left(1-2^{1-s}\right)\varsigma'(s) + 2^{1-s}\log 2.\varsigma(s)$ we have

$$\varsigma_a'(2) = \frac{1}{2}\varsigma'(2) + \frac{1}{2}\log 2.\varsigma(2)$$

and this gives us

$$3\log A - \frac{1}{12}\log 2 - \frac{1}{4}[\gamma + \log(2\pi)] = -\frac{1}{\pi^2}\left[\frac{3}{2}\varsigma'(2) + \frac{1}{2}\log 2.\varsigma(2)\right]$$

$$= -\frac{3}{2\pi^2}\varsigma'(2) - \frac{1}{12}\log 2$$

We then obtain

$$\log A - \frac{1}{12}[\gamma + \log(2\pi)] = -\frac{1}{2\pi^2}\varsigma'(2)$$

which is in agreement with (4.4.228) in Volume IV.

We now consider the integral

$$\int_0^1 \frac{1-t^{z-1}}{1-t}dt = \int_0^1 (1-t^{z-1})\sum_{n=0}^{\infty} t^n dt$$

We have $\int_0^1 (1-t^{z-1})t^n dt = \frac{1}{n+1} - \frac{1}{z+n} = \frac{z-1}{(n+1)(z+n)}$ and therefore

$$\int_0^1 \frac{1-t^{z-1}}{1-t}dt = (z-1)\sum_{n=0}^{\infty}\frac{1}{(n+1)(z+n)}$$

Reference to (E.15) then shows that (see (E.22gi) and [126, p.15])

(E.50) $$\int_0^1 \frac{1-t^{z-1}}{1-t} dt = \psi(z) + \gamma$$

With $z = n+1$ we get using (4.1.7a)

$$\int_0^1 \frac{1-t^n}{1-t} dt = \psi(n+1) + \gamma = H_n^{(1)}$$

Letting $u = 1-t$ in the integral we have

(E.51) $$\psi(z) + \gamma = \int_0^1 \frac{1-(1-u)^{z-1}}{u} du$$

and upon differentiation we get

(E.52) $$\psi'(z) = -\int_0^1 \frac{(1-u)^{z-1} \log(1-u)}{u} du$$

From (E.22a) we have

(E.53) $$\psi'(z) = \frac{\Gamma(z)\Gamma''(z) - \left(\Gamma'(z)\right)^2}{\Gamma^2(z)} = \varsigma(2, z)$$

In particular we have

$$\psi'(1) = -\int_0^1 \frac{\log(1-u)}{u} du = \varsigma(2)$$

Further differentiation of (E.52) results in

(E.54) $$\psi^{(n)}(z) = -\int_0^1 \frac{(1-u)^{z-1} \log^n(1-u)}{u} du$$

and using (E.16a) we obtain

(E.55) $$\int_0^1 \frac{(1-u)^{z-1} \log^n(1-u)}{u} du = (-1)^{n+1}(n-1)!\varsigma(n, z)$$

and (see also (4.4.233c))

(E.56) $$\int_0^1 \frac{t^{z-1}\log^n t}{1-t}\,dt = (-1)^{n+1}(n-1)!\varsigma(n,z)$$

Integrating (E.51) with respect to $z$ we obtain

$$\int_a^x \psi(z)\,dz + \gamma(x-a) = \int_a^x dz \int_0^1 \frac{1-t^{z-1}}{1-t}\,dt = \int_0^1 \frac{dt}{1-t}\int_a^x (1-t^{z-1})\,dz$$

and hence

(E.57) $$\log\Gamma(x) - \log\Gamma(a) + \gamma(x-a) = \int_0^1 \left[\frac{x-a}{1-t} - \frac{t^{x-1}-t^{a-1}}{(1-t)\log t}\right]dt$$

It is then clear that

(E.57a) $$\log\Gamma(x) + \gamma x + c = \int_0^1 \left[\frac{x}{1-t} - \frac{t^{x-1}}{(1-t)\log t}\right]dt$$

and with $x = 1$ we see that

$$\gamma + c = \int_0^1 \left[\frac{1}{1-t} - \frac{1}{(1-t)\log t}\right]dt$$

More work needs to be carried out to determine $c$.

For a positive integer $n$ and $\text{Re}(z) > 0$, repeated integration by parts yields [58b, p.2]

$$\int_0^n \left(1-\frac{t}{n}\right)^n t^{z-1}dt = \frac{n!n^z}{z(z+1)...(z+n)}$$

so that by Tannery's theorem we have

$$\lim_{n\to\infty}\int_0^n \left(1-\frac{t}{n}\right)^n t^{z-1}dt = \int_0^\infty e^{-t}t^{z-1}dt = \lim_{n\to\infty}\frac{n!n^z}{z(z+1)...(z+n)} = \Gamma(z)$$

This gives an alternative derivation of (E.11).

The following analysis is based on Locke's solution to the problem posed by Levenson in [98] in 1938. We have the following inequality from [135, p.242]

(E.58) $$0 \le e^{-x} - (1-x/n)^n \le x^2 e^{-x}/n \qquad 0 \le x \le n$$

and hence

$$I=\int_0^\infty e^{-x}\log x\,dx=\lim_{n\to\infty}\int_0^n (1-x/n)^n \log x\,dx$$

$$=\lim_{n\to\infty}\left[\sum_{k=0}^n \binom{n}{k}\frac{(-1)^k}{n^k}\int_0^n x^k \log x\,dx\right]$$

We note from (3.237i) of Volume I that

$$\int_0^n x^k \log x\,dx = n^{k+1}\left(\frac{\log n}{k+1}-\frac{1}{(k+1)^2}\right)$$

and hence

$$I=\lim_{n\to\infty} n\left[\sum_{k=0}^n \binom{n}{k}(-1)^k\left(\frac{\log n}{k+1}-\frac{1}{(k+1)^2}\right)\right]$$

$$=\lim_{n\to\infty} n\left[\log n\sum_{k=0}^n \binom{n}{k}\frac{(-1)^k}{k+1}-\sum_{k=0}^n \binom{n}{k}\frac{(-1)^k}{(k+1)^2}\right]$$

and using (E.18a) and (E.18b) this becomes

$$=\lim_{n\to\infty}\frac{n}{n+1}\left[\log n-H_{n+1}^{(1)}\right]=\lim_{n\to\infty}\left[\log n-H_n^{(1)}-\frac{1}{n+1}\right]=-\gamma$$

We therefore get the familiar result

$$\int_0^\infty e^{-x}\log x\,dx=-\gamma$$

Continuing in the same vein we have

$$\int_0^\infty e^{-x}\log^2 x\,dx=\lim_{n\to\infty}\left[\sum_{k=0}^n \binom{n}{k}\frac{(-1)^k}{n^k}\int_0^n x^k \log^2 x\,dx\right]$$

and using (3.237ii) this becomes

$$=\lim_{n\to\infty} n\left[\log^2 n\sum_{k=0}^n \binom{n}{k}\frac{(-1)^k}{(k+1)}-2\log n\sum_{k=0}^n \binom{n}{k}\frac{(-1)^k}{(k+1)^2}+2\sum_{k=0}^n \binom{n}{k}\frac{(-1)^k}{(k+1)^3}\right]$$

Using (E.18a), (E.18b) and (E.18c) we obtain

$$=\lim_{n\to\infty}\frac{n}{n+1}\left[\log^2 n-2H_{n+1}^{(1)}\log n+(H_{n+1}^{(1)})^2+H_{n+1}^{(2)}\right]$$

$$= \lim_{n\to\infty}\left[\left(\log n - H_{n+1}^{(1)}\right)^2 + H_{n+1}^{(2)}\right] = \gamma^2 + \varsigma(2)$$

Therefore we see in accordance with (E.16d) that

$$\int_0^\infty e^{-x}\log^2 x\,dx = \gamma^2 + \varsigma(2)$$

An Editorial Note in a 1938 edition of The American Mathematical Monthly [98] states that Locke's method "may be used for other positive integral powers of $\log x$ in the integrand, but it appears that no such simple form for the result is known when the exponent is odd and greater than unity". However, as noted below, it will be seen that Locke's procedure can in fact be successfully employed for $\int_0^\infty e^{-x}\log^3 x\,dx$. As before we have

$$\int_0^\infty e^{-x}\log^3 x\,dx = \lim_{n\to\infty}\left[\sum_{k=0}^n \binom{n}{k}\frac{(-1)^k}{n^k}\int_0^n x^k \log^3 x\,dx\right]$$

and with (3.237iii) of Volume I this becomes

$$= \lim_{n\to\infty} n\left[\sum_{k=0}^n \binom{n}{k}(-1)^k\left(\frac{\log^3 n}{k+1} - \frac{3\log^2 n}{(k+1)^2} + \frac{6\log n}{(k+1)^3} - \frac{6}{(k+1)^4}\right)\right]$$

$$= \lim_{n\to\infty} n\left[\frac{1}{n+1}\log^3 n - 3\log^2 n\sum_{k=0}^n \binom{n}{k}\frac{(-1)^k}{(k+1)^2} + 6\log n\sum_{k=0}^n \binom{n}{k}\frac{(-1)^k}{(k+1)^3} - 6\sum_{k=0}^n \binom{n}{k}\frac{(-1)^k}{(k+1)^4}\right]$$

$$= \lim_{n\to\infty}\frac{n}{n+1}\left[\log^3 n - 3H_{n+1}^{(1)}\log^2 n + 3(H_{n+1}^{(1)})^2\log n + 3H_{n+1}^{(2)}\log n - 6(n+1)\sum_{k=0}^n \binom{n}{k}\frac{(-1)^k}{(k+1)^4}\right]$$

$$= \lim_{n\to\infty}\left[\left(\log n - H_{n+1}^{(1)}\right)^3 + (H_{n+1}^{(1)})^3 + 3H_{n+1}^{(2)}\log n - 6(n+1)\sum_{k=0}^n \binom{n}{k}\frac{(-1)^k}{(k+1)^4}\right]$$

$$= -\gamma^3 + \lim_{n\to\infty}\left[(H_{n+1}^{(1)})^3 + 3H_{n+1}^{(2)}\log n - 6(n+1)\sum_{k=0}^n \binom{n}{k}\frac{(-1)^k}{(k+1)^4}\right]$$

$$= -\gamma^3 + \lim_{n\to\infty}\left[(H_{n+1}^{(1)})^3 - 3H_{n+1}^{(2)}\left(H_{n+1}^{(1)} - \log n\right) + 3H_{n+1}^{(2)}H_{n+1}^{(1)} - 6(n+1)\sum_{k=0}^n \binom{n}{k}\frac{(-1)^k}{(k+1)^4}\right]$$

$$= -\gamma^3 - 3\gamma\varsigma(2) + \lim_{n\to\infty}\left[(H_{n+1}^{(1)})^3 + 3H_{n+1}^{(2)}H_{n+1}^{(1)} - 6(n+1)\sum_{k=0}^n \binom{n}{k}\frac{(-1)^k}{(k+1)^4}\right]$$

We have

$$(H_{n+1}^{(1)})^3 + 3H_{n+1}^{(2)}H_{n+1}^{(1)} = (H_n^{(1)})^3 + 3\frac{(H_n^{(1)})^2}{n+1} + 6\frac{H_n^{(1)}}{(n+1)^2} + \frac{4}{(n+1)^3} + 3H_n^{(2)}H_n^{(1)} + 3\frac{H_n^{(2)}}{n+1}$$

In (E.29) we showed that

$$\lim_{n\to\infty}\left[\frac{1}{2}\left(H_n^{(1)}\right)^2 - \gamma\log n - \frac{1}{2}\log^2 n\right] = \frac{1}{2}\gamma^2$$

and hence we have

$$\lim_{n\to\infty}\frac{1}{n+1}\left[\frac{1}{2}\left(H_n^{(1)}\right)^2 - \gamma\log n - \frac{1}{2}\log^2 n\right] = 0$$

By L'Hôpital's rule we have

$$\lim_{n\to\infty}\frac{\log n}{n+1} = 0 \quad \text{and} \quad \lim_{n\to\infty}\frac{\log^2 n}{n+1} = 0$$

and thus we see that

(E.58a) $$\lim_{n\to\infty}\frac{\left(H_n^{(1)}\right)^2}{n+1} = 0$$

(which indicates how slowly $H_n^{(1)}$ diverges).

We therefore we have the limit

(E.59) $$\lim_{n\to\infty}\left[(H_{n+1}^{(1)})^3 + 3H_{n+1}^{(2)}H_{n+1}^{(1)} - 6(n+1)\sum_{k=0}^{n}\binom{n}{k}\frac{(-1)^k}{(k+1)^4}\right]$$

$$= \lim_{n\to\infty}\left[(H_n^{(1)})^3 + 3H_n^{(2)}H_n^{(1)} - 6(n+1)\sum_{k=0}^{n}\binom{n}{k}\frac{(-1)^k}{(k+1)^4}\right]$$

We now note Adamchik's identity [2] in (3.19) from Volume I

$$3\sum_{k=1}^{n}\frac{\left(H_k^{(1)}\right)^2}{k} + 3\sum_{k=1}^{n}\frac{H_k^{(2)}}{k} = \left(H_n^{(1)}\right)^3 + 3H_n^{(1)}H_n^{(2)} + 2H_n^{(3)} = 6\sum_{k=1}^{n}\frac{1}{k}\sum_{j=1}^{k}\frac{H_j^{(1)}}{j}$$

and (E.59) then becomes

$$= \lim_{n\to\infty}\left[6\sum_{k=1}^{n}\frac{1}{k}\sum_{j=1}^{k}\frac{H_j^{(1)}}{j} - 6(n+1)\sum_{k=0}^{n}\binom{n}{k}\frac{(-1)^k}{(k+1)^4} - 2H_n^{(3)}\right]$$

Alternative proofs of (E.18a) and (E.18b) were also given by Locke in [94aa] a couple of years earlier in 1936 and that paper also contains a proof by Olds that

(E.60) $$n\sum_{k=1}^{n-1}\binom{n-1}{k}\frac{(-1)^k}{(k+1)^m}=\sum_{i_1=1}^{n}\frac{1}{i_1}\sum_{i_2=1}^{i_1}\frac{1}{i_2}\sum_{i_3=1}^{i_2}\frac{1}{i_3}\cdots\sum_{i_{m-1}=1}^{i_{m-2}}\frac{1}{i_{m-1}}$$

and we note the connection with Dilcher's formula (4.1.18a)

$$\sum_{k=1}^{n}\binom{n}{k}\frac{(-1)^{k+1}}{k^s}=\sum_{1\le i_1\le i_2\le\ldots\le i_s\le n}\frac{1}{i_1 i_2 \ldots i_s}$$

Fortunately for me, [98] contained a cross reference to the work in [94aa].

We therefore see that

(E.61) $$\sum_{k=1}^{n}\frac{1}{k}\sum_{j=1}^{k}\frac{H_j^{(1)}}{j}=(n+1)\sum_{k=0}^{n}\binom{n}{k}\frac{(-1)^k}{(k+1)^4}=(n+1)\sum_{k=1}^{n}\frac{1}{k}\sum_{j=1}^{k}\frac{1}{j}\sum_{l=1}^{j}\frac{1}{l}$$

and hence (E.59) ends up as $-2\varsigma(3)$. Accordingly we have shown that

(E.62) $$\int_0^\infty e^{-x}\log^3 x\,dx=-\gamma^3-3\gamma\varsigma(2)-2\varsigma(3)$$

in accordance with [126, p.265] and (E16e).

It may be see from (E.6d) that

$$\lim_{n\to\infty}\frac{1}{n}\left\{\left(H_n\right)^2-2H_n\log n+\log^2 n\right\}=0$$

and it is then clear that

$$\lim_{n\to\infty}\frac{H_n\log n}{n}=0\,.$$

We may also compare (E.61) with (4.1.18) from Volume II

$$\sum_{k=1}^{n}\binom{n}{k}\frac{(-1)^{k+1}}{k^3}=\sum_{k=1}^{n}\frac{1}{k}\sum_{l=1}^{k}\frac{1}{l}\sum_{m=1}^{l}\frac{1}{m}=\sum_{k=1}^{n}\frac{1}{k}\sum_{l=1}^{k}\frac{H_l}{l}$$

Van der Pol gave the following expression for the gamma function

(E.63) $$\Gamma(1+x)=\left(\frac{x}{e}\right)^x\prod_{k=0}^{\infty}\frac{e_{x+k}}{e_k}$$

where $e_u=\left(1+\frac{1}{u}\right)^u$ and $e_0=1$. The following proof is based on the short paper by Nanjundiah [103ai] in 1958. We have the relation (very easy to prove, but not particularly obvious)

(E.63a) $$\frac{1+\frac{1}{k}}{1+\frac{x}{k}}=\frac{1+\frac{1}{x+k}}{1+\frac{x}{k+1}}$$

and raising both sides to the power $x+k$ we obtain

$$\frac{\left(1+\frac{1}{k}\right)^x}{1+\frac{x}{k}}=\frac{\left(1+\frac{x}{k}\right)^{x+k-1}}{\left(1+\frac{x}{k+1}\right)^{x+k}}\frac{e_{x+k}}{e_k}$$

We then have from (E.12)

$$\Gamma(x)=\frac{1}{x}\prod_{n=1}^{\infty}\left[\left(1+\frac{1}{n}\right)^x\left(1+\frac{x}{n}\right)^{-1}\right]$$

and hence we obtain

(E.63b) $$\Gamma(1+x)=\left(\frac{x+1}{e}\right)^x\prod_{k=1}^{\infty}\frac{e_{x+k}}{e_k}=\frac{e_0}{e_x}\left(\frac{x+1}{e}\right)^x\prod_{k=0}^{\infty}\frac{e_{x+k}}{e_k}$$

which is the same as (E.63) since $e_0=1$ by definition.
From (E.63a) we see that

$$\log\left(1+\frac{1}{k}\right)-\log\left(1+\frac{x}{k}\right)=\log\left(1+\frac{1}{x+k}\right)-\log\left(1+\frac{x}{k+1}\right)$$

Let us now take the logarithm of (E.63) to get

$$\log\Gamma(1+x)=x\log x-x+\sum_{k=0}^{\infty}[\log e_{x+k}-\log e_k]$$

or

(E.64) $$\log\Gamma(1+x)=x\log x-x+\sum_{k=0}^{\infty}\left[(x+k)\log\left(1+\frac{1}{x+k}\right)-k\log\left(1+\frac{1}{k}\right)\right]$$

$$= x\log x - x + \sum_{k=0}^{\infty}\left[(x+k)\log\left(x+k+1\right) - (x+k)\log(x+k) - k\log\left(1+\frac{1}{k}\right)\right]$$

Letting $x=1$ in (E.64) gives us

$$1 = \sum_{k=0}^{\infty}\left[(1+k)\log\left(1+\frac{1}{1+k}\right) - k\log\left(1+\frac{1}{k}\right)\right]$$

Letting $x=1/2$ in (E.64) gives us

(E.64a) $$\log\Gamma(3/2) = -\frac{1}{2}(1+\log 2) + \sum_{k=0}^{\infty}\left[\left(\frac{2k+1}{2}\right)\log\left(\frac{2k+3}{2k+1}\right) - k\log\left(1+\frac{1}{k}\right)\right]$$

$$= \log\left(\frac{\sqrt{\pi}}{2}\right)$$

Differentiation of (E.64) results in (after a little algebra)

(E.65) $$\psi(1+x) = \frac{\Gamma'(1+x)}{\Gamma(1+x)} = \log x + \sum_{k=0}^{\infty}\left[\log\left(1+\frac{1}{x+k}\right) - \frac{1}{x+k+1}\right]$$

$$= \log x + \sum_{k=0}^{\infty}\left[\log\left(1+\frac{1}{x+k}\right) - \frac{1}{x+k} + \left(\frac{1}{x+k} - \frac{1}{x+k+1}\right)\right]$$

$$= \log x + \frac{1}{x} + \sum_{k=0}^{\infty}\left[\log\left(1+\frac{1}{x+k}\right) - \frac{1}{x+k}\right]$$

and we therefore have

(E.66) $$\psi(x) = \log x + \sum_{k=0}^{\infty}\left[\log\left(1+\frac{1}{x+k}\right) - \frac{1}{x+k}\right]$$

in accordance with [126, p.14]. We may note from (E.66) that

(E.66a) $$\lim_{x\to\infty}[\psi(x) - \log x] = 0$$

A further differentiation of (E.65) results in

$$\psi'(1+x) = \frac{\Gamma(1+x)\Gamma''(1+x) - \left[\Gamma'(1+x)\right]^2}{\Gamma^2(1+x)} = \frac{1}{x} + \sum_{k=0}^{\infty}\left[-\frac{1}{(x+k)(x+k+1)} + \frac{1}{(x+k+1)^2}\right]$$

We will see in (E.69) that

$$\frac{1}{1+x} = \sum_{k=1}^{\infty} \frac{1}{(k+x)(k+1+x)}$$

and therefore we see that

$$\sum_{k=0}^{\infty} \frac{1}{(x+k)(x+k+1)} = \sum_{k=1}^{\infty} \frac{1}{(x+k)(x+k+1)} + \frac{1}{x(1+x)}$$

$$= \sum_{k=1}^{\infty} \frac{1}{(x+k)(x+k+1)} + \frac{1}{x} - \frac{1}{1+x} = \frac{1}{x}$$

Hence we obtain

$$\psi'(1+x) = \sum_{k=0}^{\infty} \frac{1}{(x+k+1)^2} = \varsigma(2, x+1)$$

By taking the logarithm of (E.12) we have

(E.67) $$\log\Gamma(x) = -\log x + \sum_{k=1}^{\infty} \left[ x \log\left(1+\frac{1}{k}\right) - \log\left(1+\frac{x}{k}\right) \right]$$

We also have by letting $x \to x+1$

(E.68)

$$\log\left[x\Gamma(x)\right] = \log\Gamma(1+x) = -\log(1+x) + \sum_{k=1}^{\infty} \left[ (x+1)\log\left(1+\frac{1}{k}\right) - \log\left(1+\frac{x+1}{k}\right) \right]$$

Therefore we see that

$$\log\Gamma(1+x) - \log\Gamma(x) = \log x$$

$$= \log x - \log(1+x) + \sum_{k=1}^{\infty} \left[ \log\left(1+\frac{1}{k}\right) - \log\left(1+\frac{x+1}{k}\right) + \log\left(1+\frac{x}{k}\right) \right]$$

$$= \log x - \log(1+x) + \sum_{k=1}^{\infty} \left[ \log\left(1+\frac{1}{k}\right) - \log\left(1+\frac{1}{k+x}\right) \right]$$

$$= \log x - \log(1+x) + \sum_{k=1}^{\infty} \left[ -\frac{1}{k} + \log\left(1+\frac{1}{k}\right) + \frac{1}{k} - \log\left(1+\frac{1}{k+x}\right) \right]$$

and hence using (E.6h) we have

(E.69) $$\log(1+x) = -\gamma + \sum_{k=1}^{\infty}\left[\frac{1}{k} - \log\left(1+\frac{1}{k+x}\right)\right]$$

$$= -\gamma + \sum_{m=0}^{\infty}\left[\frac{1}{m+1} - \log\left(1+\frac{1}{m+1+x}\right)\right]$$

We therefore see that

(E.69a) $$\log x = -\gamma + \sum_{k=0}^{\infty}\left[\frac{1}{k+1} - \log\left(1+\frac{1}{k+x}\right)\right]$$

which concurs with the corresponding formula in [126, p.14].

Differentiating (E.69) results in the telescoping series

(E.70) $$\frac{1}{1+x} = \sum_{k=1}^{\infty}\frac{1}{(k+x)(k+1+x)}$$

From (E.69a) we have

$$H_N - \log N = \gamma - \sum_{k=0}^{\infty}\left[\frac{1}{k+1} - \log\left(1+\frac{1}{k+N}\right)\right] + H_N$$

and thus we see that

$$\lim_{N\to\infty}\left(\sum_{k=0}^{\infty}\left[\frac{1}{k+1} - \log\left(1+\frac{1}{k+N}\right)\right] - H_N\right) = 0$$

**Theorem:**

$$1 = 1$$

**Proof:**

The above identity is well known: I rediscovered it while trying to derive something else.

We have the gamma function

$$\Gamma(t) = \int_0^{\infty} x^{t-1}e^{-x}dx$$

and therefore we see that

$$\Gamma'(t)=\int_0^\infty x^{t-1}e^{-x}\log x\,dx \text{ and } \Gamma'(1)=-\gamma=\int_0^\infty e^{-x}\log x\,dx$$

Hence we have

$$\gamma^2=\int_0^\infty e^{-x}\log x\,dx\int_0^\infty e^{-y}\log y\,dy$$

$$=\int_0^\infty\int_0^\infty e^{-(x+y)}\log x\log y\,dx\,dy$$

With the substitution $x=u^2$ and $y=v^2$ this becomes

$$=16\int_0^\infty\int_0^\infty e^{-(u^2+v^2)}\log u\log v\,du\,dv$$

and converting to polar coordinates results in

$$\gamma^2=16\int_0^\infty\int_0^\infty e^{-r^2}\log(r\cos\theta)\log(r\sin\theta)r^3\cos\theta\sin\theta dr\,d\theta$$

$$=16\int_0^\infty\int_0^\infty e^{-r^2}r^3\log^2 r\cos\theta\sin\theta\,dr\,d\theta$$

$$+16\int_0^\infty\int_0^\infty e^{-r^2}r^3\log r\log\cos\theta\cos\theta\sin\theta\,dr\,d\theta$$

$$+16\int_0^\infty\int_0^\infty e^{-r^2}r^3\log r\log(r\sin\theta)\cos\theta\sin\theta\,dr\,d\theta$$

$$+16\int_0^\infty\int_0^\infty e^{-r^2}r^3\log\cos\theta\log\sin\theta\cos\theta\sin\theta\,dr\,d\theta$$

The derivative of the beta function gives us

$$\frac{\partial}{\partial x}B(x,y)=\int_0^1 t^{x-1}\log t(1-t)^{y-1}dt$$

and with the substitution $t=\sin^2\theta$ we get at $x=y=1$

$$\left.\frac{\partial}{\partial x}B(x,y)\right|_{(x,y)=(1,1)}=4\int_0^{\pi/2}\log\sin\theta\cos\theta\sin\theta\,d\theta$$

We also have

$$\frac{\partial}{\partial x} B(x,y) = B(x,y)\left[\psi(x) - \psi(x+y)\right]$$

and hence

$$\left.\frac{\partial}{\partial x} B(x,y)\right|_{(x,y)=(1,1)} = \psi(1) - \psi(2) = -1$$

Similarly, we may show that

$$\left.\frac{\partial}{\partial y} B(x,y)\right|_{(x,y)=(1,1)} = \psi(1) - \psi(2) = -1 = 4\int_0^{\pi/2} \log\cos\theta \cos\theta \sin\theta \, d\theta$$

We also have

$$\frac{\partial^2}{\partial x \partial y} B(x,y) = B(x,y)\left[\psi(x) - \psi(x+y)\right]\left[\psi(y) - \psi(x+y)\right] - B(x,y)\psi'(x+y)$$

and this gives us

$$\left.\frac{\partial^2}{\partial x \partial y} B(x,y)\right|_{(x,y)=(1,1)} = \left[\psi(1) - \psi(2)\right]^2 - \psi'(2) = 2 - \varsigma(2)$$

With the substitution $t = r^2$ we get

$$\int_0^\infty e^{-r^2} r^3 \log^2 r \, dr = \frac{1}{2}\int_0^\infty e^{-t} t \log^2 t \, dt = \frac{1}{2}\Gamma''(2)$$

We know from (E.22) that

$$\Gamma''(2) = \varsigma(2,2) + \psi^2(2)$$

$$= \varsigma(2) - 1 + (1-\gamma)^2$$

Similarly we have

$$\int_0^\infty e^{-r^2} r^2 \log^2 r \, dr = \frac{1}{4}\int_0^\infty e^{-t} t \log t \, dt = \frac{1}{4}\Gamma'(2) = \frac{1}{4}(\gamma - 1)$$

and, unfortunately, combining everything together simply results in $\gamma^2 = \gamma^2$.

Alternatively, the Wolfram Integrator gives us

$$\int \log(a\cos x)\log(a\sin x)\cos x\sin x\,dx = A(x)+B(x)+C(x)$$

where $A(x) = -\frac{1}{2}\log^2\tan\left[\left(\frac{x}{2}\right)-i\right]+\log\sec^2\left(\frac{x}{2}\right)\log\tan\left[\left(\frac{x}{2}\right)-i\right]$

$$-\log\tan\left[\frac{1}{2}-\frac{i}{2}\tan\left(\frac{x}{2}\right)\right]\log\tan\left[\left(\frac{x}{2}\right)-i\right]+\log\tan\left[-\frac{1}{2}(1+i)\left\{\tan\left(\frac{x}{2}\right)-1\right\}\right]\log\tan\left[\left(\frac{x}{2}\right)-i\right]$$

$$+\log\tan\left[\frac{1}{2}(1-i)\left\{\tan\left(\frac{x}{2}\right)+1\right\}\right]\log\tan\left[\left(\frac{x}{2}\right)-i\right]-\log\tan^2\left[\left(\frac{x}{2}\right)-1\right]\log\tan\left[\left(\frac{x}{2}\right)-i\right]$$

$$+\log\left[1-i\tan\left(\frac{x}{2}\right)\right]\log\tan\left(\frac{x}{2}\right)+\log\left[1+i\tan\left(\frac{x}{2}\right)\right]\log\tan\left(\frac{x}{2}\right)+\log\sec^2\left(\frac{x}{2}\right)\log\tan\left[i+\left(\frac{x}{2}\right)\right]$$

$$-\log\left[\frac{1}{2}+i\frac{1}{2}\tan\left(\frac{x}{2}\right)\right]\log\left[i+\tan\left(\frac{x}{2}\right)\right]+\log\left[\left\{-\frac{1}{2}+i\frac{1}{2}\right\}\left\{\tan\left(\frac{x}{2}\right)-1\right\}\right]\log\left[i+\tan\left(\frac{x}{2}\right)\right]$$

$$+\log\left[\left\{\frac{1}{2}+i\frac{1}{2}\right\}\left\{\tan\left(\frac{x}{2}\right)+1\right\}\right]\log\left[i+\tan\left(\frac{x}{2}\right)\right]-\log^2\tan\left[\left(\frac{x}{2}\right)+i\right]\log\left[i+\tan\left(\frac{x}{2}\right)\right]$$

$B(x) =$

$$-\frac{1}{4}\cos 2x+\frac{1}{4}(\cos 2x+1)\log\cos x+\frac{1}{4}\cos 2x\log a+\frac{1}{4}(\cos 2x-1)\log\sin x$$

$$+\log\sec^2(x/2)\log a+\log\sec^2(x/2)\left[\log\sin x-\log\tan(x/2)\right]$$

$$+\frac{1}{2}(1-\cos 2x)\log a\log\cos x-\log a\log\left[1-\tan^2(x/2)\right]$$

$$+\frac{1}{2}(1-\cos 2x)\log^2 a+\frac{1}{2}(1-\cos 2x)\log a\log\sin x+\frac{1}{4}(\cos 2x-1)\log a$$

$$+\frac{1}{2}(1-\cos 2x)\log\cos x\log\sin x-\log\tan(x/2)\log\left[\tan(x/2)+1\right]$$

$$-\log\sin x\log\left[1-\tan^2(x/2)\right]+\log\tan(x/2)\log\left[1-\tan^2(x/2)\right]$$

and $C(x) =$

$$Li_2\left[1-\tan\left(\frac{x}{2}\right)\right]+Li_2\left[\frac{1}{2}(1+i)-\frac{1}{2}(1-i)\tan\left(\frac{x}{2}\right)\right]-Li_2\left[\frac{1}{2}-\frac{1}{2}i\tan\left(\frac{x}{2}\right)\right]$$

$$-Li_2\left[\frac{1}{2}+\frac{1}{2}i\tan\left(\frac{x}{2}\right)\right]-Li_2\left[-\tan\left(\frac{x}{2}\right)\right]+Li_2\left[-i\tan\left(\frac{x}{2}\right)\right]$$

$$+Li_2\left[i\tan\left(\frac{x}{2}\right)\right]+Li_2\left[-\left\{\frac{1}{2}+\frac{1}{2}i\right\}\left\{i+\tan\left(\frac{x}{2}\right)\right\}\right]+Li_2\left[\frac{1}{2}(1-i)\tan\left(\frac{x}{2}\right)+\frac{1}{2}(1+i)\right]$$

$$+Li_2\left[\frac{1}{2}(1+i)\tan\left(\frac{x}{2}\right)+\frac{1}{2}(1-i)\right]$$

With $x=\pi/2$ we see that

$$C(\pi/2)=Li_2(0)+Li_2(i)-Li_2\left[\frac{1}{2}(1-i)\right]-Li_2\left[\frac{1}{2}(1+i)\right]-Li_2(-1)$$

$$+Li_2(-i)+Li_2(i)+Li_2\left[-\frac{1}{2}(1+i)^2\right]+2Li_2(1)$$

Using Euler's identities for the dilogarithm we see that

$$Li_2\left[\frac{1}{2}(1-i)\right]+Li_2\left[\frac{1}{2}(1+i)\right]=\varsigma(2)-\log\left[\frac{1}{2}(1-i)\right]\log\left[\frac{1}{2}(1+i)\right]$$

We note from (4.4.67b) that

$$Li_2(i)+Li_2(i)=\frac{1}{2}Li_2(-1)=-\frac{\pi^2}{24}$$

Since $\sqrt{i}=\pm\frac{1+i}{\sqrt{2}}$ we also have $Li_2\left[-\frac{1}{2}(1+i)^2\right]=Li_2(-i)$. Hence we obtain

$$C(\pi/2)=3\varsigma(2)-\log\left[\frac{1}{2}(1-i)\right]\log\left[\frac{1}{2}(1+i)\right]$$

Similarly we find that

$$C(0)=-\varsigma(2)+2Li_2(1/2)-2\varsigma(2)+2\log\left[\frac{1}{2}(1-i)\right]\log\left[\frac{1}{2}(1+i)\right]$$

and, with a lot more work (involving L'Hôpital's rule), we obtain the same answer.

□

Appendix E is now closed with a useful trick.

We may alternatively write (E.7) as

$$H_n = \int_0^1 \frac{1-(1-x)^n}{x}\,dx = -n\int_0^1 (1-x)^n \log x\,dx$$

and we then have

$$H_n - \log n = -n\int_0^1 (1-x)^n \log x\,dx - n\int_0^1 (1-x)^n \log n\,dx$$

$$= -n\int_0^1 (1-x)^n \log(nx)\,dx$$

With the substitution $u = nx$ this becomes

$$= -\int_0^n \left(1-\frac{u}{n}\right)^{n-1} \log u\,du$$

and we obtain

$$\gamma = -\lim_{n\to\infty}\int_0^n \left(1-\frac{u}{n}\right)^{n-1} \log u\,du$$

It may then be rigorously demonstrated that

$$\lim_{n\to\infty}\int_0^n \left(1-\frac{u}{n}\right)^{n-1} \log u\,du = \int_0^\infty e^{-u} \log u\,du = \Gamma'(1)$$

□

## APPENDIX F

## SOME ELEMENTARY ASPECTS OF RIEMANN'S FUNCTIONAL EQUATION FOR THE ZETA FUNCTION

In 1768 Euler asserted that

$$\frac{1-2^{n-1}+3^{n-1}-4^{n-1}+5^{n-1}-...}{1-2^{-n}+3^{-n}-4^{-n}+5^{-n}-...}=\frac{(n-1)!\left(2^n-1\right)}{\left(2^{n-1}-1\right)\pi^n}\cos\left(\frac{n\pi}{2}\right)$$

which he verified for $n=1$ and $n=2k$. A little later, in 1906 to be precise, Landau wrote Euler's identity as

$$\frac{\lim_{x\to1}\sum_{n=1}^{\infty}(-1)^{n+1}n^{s-1}x^{n-1}}{\lim_{x\to1}\sum_{n=1}^{\infty}(-1)^{n+1}n^{-s}x^{n-1}}=\frac{\Gamma(s)\left(2^s-1\right)}{\left(2^{s-1}-1\right)\pi^s}\cos\left(\frac{n\pi}{2}\right)$$

proved its validity, and showed its equivalence to the functional equation for the Riemann zeta function [57]. Great achievements for both mathematicians!

The functional equation for the Riemann zeta function is shown below

(F.1) $$\varsigma(1-s)=2(2\pi)^{-s}\Gamma(s)\cos(\pi s/2)\varsigma(s)$$

and we referred to this in (4.4.199).

Seven different proofs of this are given in Titchmarsh's treatise [129] and one of these is shown towards the end of this Appendix.

Employing Euler's reflection formula (C.1) for the gamma function

$$\Gamma(s)\Gamma(1-s)=\frac{\pi}{\sin\pi s}$$

we obtain an equivalent form

(F.1a) $$\varsigma(s)=2(2\pi)^{s-1}\Gamma(1-s)\sin(\pi s/2)\varsigma(1-s)$$

Using the Hasse/Sondow identity, it was previously shown in (3.11a), (3.11b) and (3.11c) that

(F.2) $$\varsigma(0)=-\frac{1}{2}\ ,\ \varsigma(-1)=-\frac{1}{12}\text{ and }\varsigma(-2)=0$$

With regard to $\varsigma(0)$, see also the comment following equation (6.48) in Volume V.

Formally, with $s=0$, (F.1) suggests that

(F.3) $$\varsigma(0)=\lim_{s\to0}\left[\frac{\varsigma(1-s)}{2(2\pi)^{-s}\Gamma(s)\cos(\pi s/2)}\right]=-\frac{1}{2}$$

and therefore

$$\lim_{s\to0}\left[\frac{\varsigma(1-s)}{\Gamma(s)}\right]=1$$

Using L'Hôpital's rule this implies that

$$\lim_{s\to0}\left[\frac{\varsigma'(1-s)}{\Gamma'(s)}\right]=-1$$

In fact, Whittaker and Watson [135, p.266] show that

$$\lim_{s\to0}\left[\frac{\varsigma(s,a)}{\Gamma(1-s)}\right]=-1$$

which we have also shown in (4.3.204a) in Volume II(b).

As is well-known, the zeta function has a pole at $s=1$. With $s=2$ in (F.1) we get as before

(F.4) $$\varsigma(-1)=2(2\pi)^{-2}\Gamma(2)\cos\pi\,\varsigma(2)=-\frac{1}{12}$$

Letting $s=2n$ we get

$$\varsigma(1-2n)=2(2\pi)^{-2n}\Gamma(2n)\cos(n\pi)\varsigma(2n)$$

Therefore using Euler's identity (1.7) for $\varsigma(2n)$ we have

(F.4a) $$\varsigma(1-2n)=-\frac{B_{2n}}{2n}\quad\forall\ n\geq1$$

Putting $s=2n+1$ in (F.1) we obtain the so-called trivial zeros of the Riemann zeta function because $\cos\left(\frac{(2n+1)\pi}{2}\right)=0$ for all $n$

(F.4b) $$\varsigma(-2n)=0\quad\forall\ n\geq1$$

We see from (F.1) that $\varsigma(0)=2(2\pi)^{-1}\lim_{s\to1}[\cos(\pi s/2)\varsigma(s)]$ and hence

(F.4c) $$\lim_{s\to 1}[\cos(\pi s/2)\varsigma(s)] = -\frac{\pi}{2}$$

Applying logarithmic differentiation to the functional equation it is easily shown that

(F.5) $$\frac{\varsigma'(1-s)}{\varsigma(1-s)} = \log(2\pi) - \frac{\Gamma'(s)}{\Gamma(s)} + \frac{\pi}{2}\tan\left(\frac{\pi s}{2}\right) - \frac{\varsigma'(s)}{\varsigma(s)}$$

Taking the limit as $s \to 1$ we have

$$\frac{\varsigma'(0)}{\varsigma(0)} = \log(2\pi) - \gamma + \lim_{s\to 1}\left[\frac{\pi}{2}\tan\left(\frac{\pi s}{2}\right) - \frac{\varsigma'(s)}{\varsigma(s)}\right]$$

We know from Titchmarsh [129, p.20] that in the neighbourhood of $s = 1$

(F.5a) $$\frac{\pi}{2}\tan\left(\frac{\pi s}{2}\right) = -\frac{1}{s-1} + O\left(|s-1|\right)$$

(F.5b) $$\frac{\varsigma'(s)}{\varsigma(s)} = \frac{-1\left\{1/(s-1)^2\right\} + k + \ldots}{\left\{1/(s-1)\right\} + \gamma + k(s-1)\ldots} = -\frac{1}{s-1} + \gamma + \ldots$$

(where $k$ is a constant). See also (4.4.99a). We then see that

(F.5bi) $$\lim_{s\to 1}\left[\frac{\pi}{2}\tan\left(\frac{\pi s}{2}\right) - \frac{\varsigma'(s)}{\varsigma(s)}\right] = \gamma$$

Therefore we have

(F.5c) $$\frac{\varsigma'(0)}{\varsigma(0)} = \log(2\pi)$$

and hence

(F.6) $$\varsigma'(0) = -\frac{1}{2}\log(2\pi) = -\log\sqrt{2\pi}$$

This was derived in a different way in (4.3.116a) in Volume II(a).

Employing (F.5) with $s = 2$ we have

$$\frac{\varsigma'(-1)}{\varsigma(-1)} = \log(2\pi) - \frac{\Gamma'(2)}{\Gamma(2)} - \frac{\varsigma'(2)}{\varsigma(2)}$$

We have from (4.3.16)

$$\psi(n)=\frac{\Gamma'(n)}{\Gamma(n)}=H_{n-1}^{(1)}-\gamma \text{ and hence } \psi(2)=\frac{\Gamma'(2)}{\Gamma(2)}=1-\gamma$$

We therefore get

(F.7) $$\varsigma'(-1)=\frac{1}{12}(1-\gamma-\log 2\pi)+\frac{1}{2\pi^2}\varsigma'(2)$$

Equivalently we have

(F.7a) $$\varsigma'(-1)=\frac{1}{12}(1-\gamma-\log 2\pi)-\frac{1}{2\pi^2}\sum_{n=1}^{\infty}\frac{\log n}{n^2}$$

Putting $s=2n$ in (F.5) we obtain

(F.8) $$2n\varsigma'(1-2n)-\left[H_{2n-1}^{(1)}-\gamma-\log 2\pi\right]B_{2n}=\frac{\varsigma'(2n)}{\varsigma(2n)}B_{2n}$$

From the definition of a derivative and using (F.4b) we have

$$\varsigma'(-2n)=\lim_{h\to 0}\left[\frac{\varsigma(-2n+h)-\varsigma(-2n)}{h}\right]=\lim_{h\to 0}\left[\frac{\varsigma(-2n+h)}{h}\right]$$

Using the form of the functional equation given in (F.1a) we get

$$\frac{\varsigma(-2n+h)}{h}=2(2\pi)^{-2n+h-1}\Gamma(1+2n-h)\frac{\sin\left[\frac{\pi}{2}(-2n+h)\right]}{h}\varsigma(2n+1-h)$$

and accordingly we get

$$\varsigma'(-2n)=2(2\pi)^{-2n-1}(2n)!\varsigma(2n+1)\lim_{h\to 0}\left\{\frac{\sin\left[\frac{\pi}{2}(-2n+h)\right]}{h}\right\}$$

Since $\lim_{h\to 0}\left\{\frac{\sin\left[\frac{\pi}{2}(-2n+h)\right]}{h}\right\}=(-1)^n\frac{\pi}{2}$ we easily obtain

(F.8a) $$\varsigma'(-2n)=(-1)^n\frac{(2n)!}{2(2\pi)^{2n}}\varsigma(2n+1)$$

which we also saw in (4.3.112j) in Volume II(a).

Therefore for example we get

(F.8b) $$\varsigma'(-2) = -\frac{\varsigma(3)}{4\pi^2}$$

Rather than using logarithmic differentiation we may differentiate (F.1) directly to obtain an expression which does not contain the term $\varsigma(1-s)$

$$-\varsigma'(1-s) = 2(2\pi)^{-s}\left[-\log(2\pi)\Gamma(s)\cos(\pi s/2)\varsigma(s) + \Gamma'(s)\cos(\pi s/2)\varsigma(s)\right.$$

$$\left. + \Gamma(s)\cos(\pi s/2)\varsigma'(s) - \Gamma(s)\frac{\pi}{2}\sin(\pi s/2)\varsigma(s)\right]$$

and with $s = 3$ we simply recover (F.8b).

Differentiating again we obtain

(F.8c)
$$-\frac{1}{2}\frac{d}{ds}\left\{\frac{\varsigma'(1-s)}{(2\pi)^s}\right\} =$$

$$\left[-\log(2\pi)\Gamma'(s)\cos(\pi s/2)\varsigma(s) - \log(2\pi)\Gamma(s)\cos(\pi s/2)\varsigma'(s) - \log(2\pi)\Gamma(s)\frac{\pi}{2}\sin(\pi s/2)\varsigma(s)\right.$$

$$+\Gamma''(s)\cos(\pi s/2)\varsigma(s) + \Gamma'(s)\cos(\pi s/2)\varsigma'(s) - \Gamma'(s)\frac{\pi}{2}\sin(\pi s/2)\varsigma(s)$$

$$+\Gamma'(s)\cos(\pi s/2)\varsigma'(s) + \Gamma(s)\cos(\pi s/2)\varsigma''(s) - \Gamma(s)\frac{\pi}{2}\sin(\pi s/2)\varsigma'(s)$$

$$\left. -\Gamma'(s)\frac{\pi}{2}\sin(\pi s/2)\varsigma(s) \; -\Gamma(s)\frac{\pi}{2}\sin(\pi s/2)\varsigma'(s) \; -\Gamma(s)\left(\frac{\pi}{2}\right)^2\cos(\pi s/2)\varsigma(s)\right]$$

Hardy [129, p.16] gave the following functional equation for the alternating zeta function

(F.8d) $$\varsigma_a(-s) = \left(1-\left[2^{-s}-1\right]^{-1}\right)\pi^{-s-1}s\Gamma(s)\sin(\pi s/2)\varsigma_a(1+s)$$

$$= 2\frac{\left[2^{-s-1}-1\right]}{\left[2^{-s}-1\right]}\pi^{-s-1}s\Gamma(s)\sin(\pi s/2)\varsigma_a(1+s)$$

We therefore have for example

(F.8d) $$\varsigma_a(-1) = \frac{3}{\pi^2}\varsigma_a(2) \text{ and } \varsigma_a(-2) = 0$$

Since $\varsigma_a(s) = (1-2^{1-s})\varsigma(s)$ we get $\varsigma_a(2) = \frac{1}{2}\varsigma(2)$ and hence, as in (3.11b), we obtain

(F.8e) $$\varsigma_a(-1) = \frac{1}{4}$$

Logarithmic differentiation results in

(F.8f)

$$-\frac{\varsigma_a'(-s)}{\varsigma_a(-s)} = -\frac{2^{-s-1}\log 2}{\left[2^{-s-1}-1\right]} + \frac{2^{-s}\log 2}{\left[2^{-s}-1\right]} - \log\pi + \frac{1}{s} + \frac{\Gamma'(s)}{\Gamma(s)} + \frac{\pi}{2}\cot(\pi s/2) + \frac{\varsigma_a'(1+s)}{\varsigma_a(1+s)}$$

and with $s = 1$ we obtain

(F.8g) $$-\frac{\varsigma_a'(-1)}{\varsigma_a(-1)} = -\frac{2}{3}\log 2 - \log\pi + 1 - \gamma + \frac{\varsigma_a'(2)}{\varsigma_a(2)}$$

Since $\varsigma_a(s) = (1-2^{1-s})\varsigma(s)$ (which is valid for s greater than zero) we have

$$\varsigma_a'(s) = (1-2^{1-s})\varsigma'(s) + 2^{1-s}\varsigma(s)\log 2$$

$$\varsigma_a'(2) = \frac{1}{2}\varsigma'(2) + \frac{1}{2}\varsigma(2)\log 2$$

(F.8h) $$\frac{\varsigma_a'(2)}{\varsigma_a(2)} = \frac{\varsigma'(2)}{\varsigma(2)} + \log 2$$

From (F.7) we have

(F.8i) $$\frac{\varsigma'(2)}{\varsigma(2)} = 12\left[\frac{1}{2\pi^2}\varsigma'(2)\right] = 12\left[\varsigma'(-1) - \frac{1}{12}(1-\gamma-\log 2\pi)\right]$$

and hence we obtain

(F.8j) $$\varsigma_a'(-1) = -3\varsigma'(-1) - \frac{1}{3}\log 2$$

Since $\varsigma_a(s) = (1-2^{1-s})\varsigma(s)$ we have (assuming that this identity remains valid for s less than 0)

$$\varsigma_a'(s) = (1-2^{1-s})\varsigma'(s) + 2^{1-s}\varsigma(s)\log 2$$

$$\varsigma_a'(-1) = -3\varsigma'(-1) + 4\varsigma(-1)\log 2$$

and since by (3.11b) $\varsigma(-1) = -\frac{1}{12}$ we have as before

(F.8k) $$\varsigma_a'(-1) = -3\varsigma'(-1) - \frac{1}{3}\log 2$$

Since

$$\frac{1}{s} + \frac{\Gamma'(s)}{\Gamma(s)} = \frac{\Gamma(s)+\Gamma'(s)}{s\Gamma(s)} = \frac{1}{\Gamma(s+1)}\frac{d}{ds}\left[s\Gamma(s)\right] = \frac{\Gamma'(s+1)}{\Gamma(s+1)}$$

we have

$$\lim_{s\to 0}\left[\frac{1}{s} + \frac{\Gamma'(s)}{\Gamma(s)}\right] = \lim_{s\to 0}\left[\frac{\Gamma'(s+1)}{\Gamma(s+1)}\right] = \Gamma'(1) = -\gamma$$

Therefore, we may take the limit as $s \to 0$ of the following equation

$$-\frac{\varsigma_a'(-s)}{\varsigma_a(-s)} = -\frac{2^{-s-1}\log 2}{\left[2^{-s-1}-1\right]} + \frac{2^{-s}\log 2}{\left[2^{-s}-1\right]} - \log\pi + \frac{1}{s} + \frac{\Gamma'(s)}{\Gamma(s)} + \frac{\pi}{2}\cot(\pi s/2) + \frac{\varsigma_a'(1+s)}{\varsigma_a(1+s)}$$

to obtain

$$-\frac{\varsigma_a'(0)}{\varsigma_a(0)} = \log 2 - \log\pi - \gamma + \frac{\varsigma_a'(1)}{\varsigma_a(1)} + \lim_{s\to 0}\left(\frac{\pi}{2}\cot(\pi s/2) + \frac{\log 2}{\left[1-2^{s}\right]}\right)$$

Since $t^{-s}$ is a monotonic decreasing function of $t$ we have

$$(n+1)^{-s} < \int_n^{n+1} t^{-s}dt < n^{-s}$$

Summing from $n=1$ to $\infty$ we get

$$\varsigma(s) - 1 < \int_1^{\infty} t^{-s}dt < \varsigma(s)$$

We have $\int_1^{\infty} t^{-s}dt = \frac{1}{s-1}$ and hence $1 < (s-1)\varsigma(s) < s$ .We therefore obtain

$$\lim_{s\to 1}(s-1)\varsigma(s) = 1$$

We may also see this from the Hasse formula (3.12)

$$\varsigma(s) = \frac{1}{s-1}\sum_{n=0}^{\infty}\frac{1}{n+1}\sum_{k=0}^{n}\binom{n}{k}\frac{(-1)^k}{(k+1)^{s-1}}$$

for then we have

$$\lim_{s\to 1}(s-1)\varsigma(s) = \sum_{n=0}^{\infty}\frac{1}{n+1}\sum_{k=0}^{n}\binom{n}{k}(-1)^k = \sum_{n=0}^{\infty}\frac{1}{n+1}\delta_{n,0} = 1$$

The following proof of the Riemann functional equation is based on the fifth method given by Titchmarsh in his book [129, p.24].

In (4.4.38) we proved that

$$\Gamma(s)\varsigma(s) = \int_0^{\infty}\frac{x^{s-1}}{e^x-1}dx \qquad , \mathrm{Re}(s) = \sigma > 1$$

For $\sigma > 1$ this may be written as

$$\Gamma(s)\varsigma(s) = \int_0^1\left[\frac{1}{e^x-1} - \frac{1}{x}\right]x^{s-1}dx + \frac{1}{s-1} + \int_1^{\infty}\frac{x^{s-1}}{e^x-1}dx$$

and this holds by analytic continuation for $\sigma > 0$. Also, for $1 > \sigma > 0$, we have

$$\frac{1}{s-1} = \int_0^1\frac{x^{s-1}}{x}dx$$

and therefore

(F.9) $$\Gamma(s)\varsigma(s) = \int_0^{\infty}\left[\frac{1}{e^x-1} - \frac{1}{x}\right]x^{s-1}dx \qquad 1 > \sigma > 0$$

Similarly we may also obtain

$$\Gamma(s)\varsigma(s) = \int_0^1\left[\frac{1}{e^x-1} - \frac{1}{x} + \frac{1}{2}\right]x^{s-1}dx - \frac{1}{2s} + \int_1^{\infty}\frac{x^{s-1}}{e^x-1}dx$$

and this holds by analytic continuation for $\sigma > -1$. We have

$$-\frac{1}{2s} = \int_1^{\infty}\frac{x^{s-1}}{x}dx \qquad , 0 > \sigma > -1$$

Therefore we obtain

(F.10) $$\Gamma(s)\varsigma(s)=\int_0^\infty\left[\frac{1}{e^x-1}-\frac{1}{x}+\frac{1}{2}\right]x^{s-1}dx \qquad ,0>\sigma>-1$$

Substituting $t=-\frac{1}{2}ix$ in the following identity

$$\sin t=t\prod_{n=1}^\infty\left(1-\frac{t^2}{n^2\pi^2}\right)$$

and then performing logarithmic differentiation gives us the formula for $x\neq 2m\pi i$

$$\frac{1}{e^x-1}=\frac{1}{x}-\frac{1}{2}+2x\sum_{n=1}^\infty\frac{1}{4n^2\pi^2+x^2}$$

Hence we get

$$\Gamma(s)\varsigma(s)=2\int_0^\infty x\sum_{n=1}^\infty\frac{1}{4n^2\pi^2+x^2}x^{s-1}dx=2\sum_{n=1}^\infty\int_0^\infty\frac{x^s}{4n^2\pi^2+x^2}dx$$

$$=2\sum_{n=1}^\infty(2n\pi)^{s-1}\frac{\pi}{2\cos(s\pi/2)}=\frac{2^{s-1}\pi^s}{\cos(s\pi/2)}\varsigma(1-s)$$

Sebah and Gourdon give a slightly easier proof on their website [119a] and this is set out below (this proof is due to Hans Rademacher (1892-1969); see [110aa]).

We first of all apply the Euler-Maclaurin summation formula to the function $f(x)=x^{-s}$ to obtain

(F.11) $$\varsigma(s)=\frac{1}{s-1}+\frac{1}{2}+\sum_{r=2}^q\frac{B_r}{r!}s(s+1)...(s+r-2)-\frac{1}{q!}s(s+1)...(s+q-1)\int_1^\infty B_q\left(\{x\}\right)x^{-s-q}dx$$

where $B_q\left(\{x\}\right)$ are the Bernoulli polynomials evaluated at the fractional part $\{x\}=x-[x]$. The integral converges for $\operatorname{Re}(s)>1-q$ and, since $q$ can be arbitrarily large, it can therefore be seen that $\varsigma(s)$ is analytic on the whole complex plane with a simple pole at $s=1$. It can instantly be seen from this equation that $\varsigma(0)=-\frac{1}{2}$ (and this corrects the misprint in the paper by Sebah and Gourdon). With $s=-n$ where $n$ is a positive integer, and choosing $q=n+1$ we find

$$\varsigma(-n)=-\frac{1}{n+1}+\frac{1}{2}-\sum_{r=2}^{n+1}\frac{B_r}{r!}n(n-1)...(n-r+2)$$

$$= -\frac{1}{n+1}\sum_{r=0}^{n+1}\binom{n+1}{r}B_r$$

where we have used the fact (A.7) that $B_{2r+1} = 0$ for $r \geq 2$. Using (A.6) we have

$$= -\frac{B_{n+1}}{n+1}$$

We therefore see that

(F.12a) $\qquad \varsigma(-2n) = 0$

(F.12b) $\qquad \varsigma(1-2n) = -\frac{B_{2n}}{2n}$

As mentioned by Rademacher [110aa, p.81], the zeros of $\varsigma(s)$ at $s = -2n$ are often called the "trivial zeros" because they are so easily found.

We now use (F.11) with $q = 3$

(F.13) $\varsigma(s) = \frac{1}{s-1} + \frac{1}{2} + \frac{B_2}{2}s - \frac{1}{6}s(s+1)(s+2)\int_1^\infty B_3\left(\{x\}\right)x^{-s-3}dx$

which is valid for $\text{Re}(s) > -2$.

From (A.15) we have

$$B_3(x) = x\left(x - \frac{1}{2}\right)(x-1)$$

We now evaluate the following integral where $\text{Re}(s) < -1$ using integration by parts and (A.13)

$$\int_0^1 B_3(x)x^{-s-3}dx = -B_3(x)\frac{x^{-s-2}}{s+2}\Bigg|_0^1 + \frac{3}{s+2}\int_0^1 B_2(x)x^{-s-2}dx$$

$$= \frac{3}{s+2}\int_0^1 B_2(x)x^{-s-2}dx$$

Similarly we have

$$\int_0^1 B_2(x)x^{-s-2}dx = -B_2(x)\frac{x^{-s-1}}{s+1}\Bigg|_0^1 + \frac{2}{s+1}\int_0^1 B_1(x)x^{-s-1}dx$$

$$= -\frac{1}{6}\frac{1}{s+1} + \frac{2}{s+1}\int_0^1 B_1(x)x^{-s-1}dx$$

Finally, we have

$$\int_0^1 B_1(x)x^{-s-1}dx = -B_1(x)\frac{x^{-s}}{s}\Bigg|_0^1 + \frac{1}{s}\int_0^1 B_0(x)x^{-s}dx$$

$$= -\frac{1}{2s} + \frac{1}{s(1-s)}$$

Therefore we have

$$\int_0^1 B_3(x)x^{-s-3}dx = \frac{3}{s+2}\left[-\frac{1}{6}\frac{1}{s+1} + \frac{2}{s+1}\left\{-\frac{1}{2s} + \frac{1}{s(1-s)}\right\}\right]$$

$$= \frac{6}{s(s+1)(s+2)}\left[-\frac{s}{12} - \frac{1}{2} - \frac{1}{(s-1)}\right]$$

Combining this with (F.13) we get for $-2 < \mathrm{Re}(s) < -1$

(F.14) $$\varsigma(s) = -\frac{1}{6}s(s+1)(s+2)\int_0^\infty B_3\left(\{x\}\right)x^{-s-3}dx$$

We now use the Fourier series

$$B_3\left(\{x\}\right) = 12\sum_{n=1}^\infty \frac{\sin 2\pi nx}{(2\pi n)^3}$$

and obtain

$$\varsigma(s) = -2s(s+1)(s+2)\int_0^\infty \sum_{n=1}^\infty \frac{\sin 2\pi nx}{(2\pi n)^3 x^{s+3}}dx$$

$$= -2s(s+1)(s+2)\sum_{n=1}^\infty \frac{1}{(2\pi n)^3}\int_0^\infty \frac{\sin 2\pi nx}{x^{s+3}}dx$$

$$= -2s(s+1)(s+2)\sum_{n=1}^\infty \frac{1}{(2\pi n)^{1-s}}\int_0^\infty \frac{\sin t}{t^{s+3}}dt$$

The integral $\int_0^\infty \frac{\sin t}{t^{s+3}}dt$ is well-known and an easy evaluation using contour integration may be found in [24a]: the result is

(F.15) $$\int_0^\infty \frac{\sin t}{t^{s+3}}dt = \Gamma(-s-2)\sin\frac{\pi(s+2)}{2}$$

Thus

$$\varsigma(\text{s}) = -2s(s+1)(s+2)\ \Gamma(-s-2)\sin\frac{\pi(s+2)}{2}\ \frac{1}{(2\pi)^{1-s}}\sum_{n=1}^\infty \frac{1}{n^{1-s}}$$

Therefore we obtain

$$\varsigma(\text{s}) = 2^s\pi^{s-1}\Gamma(1-s)\varsigma(1-s)\sin\frac{\pi s}{2}$$

This is the functional equation and has been proved for complex values in the strip $-2 < \text{Re}(s) < -1$ : by analytic continuation it is valid throughout the whole of the complex plane.

More than 30 years have elapsed since I studied contour integration, and hence an alternative evaluation of the integral (F.15) is set out below using mainly real analysis.

From (4.4.57c) we have

$$\frac{1}{x^p} = \frac{1}{\Gamma(p)}\int_0^\infty u^{p-1}e^{-xu}du$$

Therefore we obtain

$$\int_0^\infty \frac{\sin x}{x^p}dx = \frac{1}{\Gamma(p)}\int_0^\infty\int_0^\infty u^{p-1}e^{-xu}\sin x\, dudx$$

(F.16) $$= \frac{1}{\Gamma(p)}\int_0^\infty u^{p-1}du\int_0^\infty e^{-xu}\sin x\, dx$$

By elementary calculus we have

$$\int_0^M e^{-xu}\sin x\, dx = \text{Im}\int_0^M e^{-xu}e^{ix}\, dx = \text{Im}\int_0^M e^{x(i-u)}\, dx = \text{Im}\left[\frac{e^{-xu}(\cos x + i\sin x)}{i-u}\right]_0^M$$

$$= \left[ \frac{e^{-xu}(-\cos x + u \sin x)}{1+u^2} \right]_0^M$$

Hence we get

$$\int_0^\infty e^{-xu} \sin x \, dx = \frac{1}{1+u^2}$$

and substituting this in (F.16) we have

$$\text{(F.17)} \qquad \int_0^\infty \frac{\sin x}{x^p} dx = \frac{1}{\Gamma(p)} \int_0^\infty \frac{u^{p-1}}{1+u^2} du$$

Letting $v = u^2$ in the last integral we get

$$\int_0^\infty \frac{u^{p-1}}{1+u^2} du = \frac{1}{2} \int_0^\infty \frac{v^{(p-3)/2}}{1+v} dv$$

Upon using (C.2) we obtain

$$= \frac{1}{2} \frac{\pi}{\sin \frac{\pi(p-1)}{2}}$$

Therefore, letting $p = s+3$ we obtain

$$\int_0^\infty \frac{\sin x}{x^{s+3}} dx = -\frac{1}{2\Gamma(s+3)} \frac{\pi}{\sin(\pi s/2)}$$

The letter which Ramanujan wrote to Hardy in 1913 [76, p.351] contained the following divergent series

$$1+2+3+\ldots = -\frac{1}{12}$$

by which Ramanujan meant $\varsigma(-1) = -\frac{1}{12}$ (as seen above in (F.2)). From (A.2) we have

$$\frac{1}{1-e^t} = -\sum_{n=0}^\infty B_n \frac{t^{n-1}}{n!}$$

and expanding the left-hand side as a geometric series we obtain

$$1+e^t+e^{2t}+...=1+\sum_{n=0}^{\infty}\frac{1^n}{n!}t^n+\sum_{n=0}^{\infty}\frac{2^n}{n!}t^n+...$$

The coefficient of $t^{n-1}$ in this formal expansion is for $n>1$

$$\frac{1}{(n-1)!}(1^{n-1}+2^{n-1}+...)\approx\frac{1}{(n-1)!}\varsigma(-n+1)$$

and the coefficient of $t^0$ is formally

$$1+\frac{1}{0!}(1^0+2^0+...)$$

Thus we have an "explanation" of why

$$\varsigma(-n+1)=-\frac{B_n}{n}$$

For more details see the recent paper by Doyon et al. [55b].

We now recall the Hasse/Sondow formula (3.11)

$$\varsigma_a(s)=\sum_{n=0}^{\infty}\frac{1}{2^{n+1}}\sum_{k=0}^{n}\binom{n}{k}\frac{(-1)^k}{(k+1)^s}$$

and hence we have (with $m$ being a positive integer)

$$\varsigma_a(-m)=\sum_{n=0}^{\infty}\frac{1}{2^{n+1}}\sum_{k=0}^{n}\binom{n}{k}(-1)^k(k+1)^m=\sum_{n=0}^{\infty}\frac{G_n(m)}{2^{n+1}} \tag{F.18}$$

where

$$G_n(m)=\sum_{k=0}^{n}\binom{n}{k}(-1)^k(k+1)^m$$

Reference to (3.11aa) shows that

$$G_n(1)=x\frac{d}{dx}\left[x(1-x)^n\right]\Big|_{x=1}$$

$$G_n(2)=x\frac{d}{dx}\left\{x\frac{d}{dx}\left[x(1-x)^n\right]\right\}\Big|_{x=1}$$

and more generally

$$G_n(m) = \left( x\frac{d}{dx} \right)^m f_n(x) \Bigg|_{x=1}$$

where $f_n(x) = x(1-x)^n$. Therefore, using (3.97) of Volume I

$$\left( x\frac{d}{dx} \right)^m g(x) = \sum_{k=1}^{m} S(m,k) x^k g^{(k)}(x)$$

we obtain

$$G_n(m) = \sum_{k=1}^{m} S(m,k) x^k f_n^{(k)}(x) \Bigg|_{x=1}$$

We have by differentiation

$$f_n^{(1)}(x) = (1-x)^n - nx(1-x)^{n-1}$$

$$f_n^{(2)}(x) = -2n(1-x)^{n-1} + n(n-1)x(1-x)^{n-2}$$

We conjecture that

$$f_n^{(k)}(x) = (-1)^{k+1} kn(n-1)...(n-k+2)(1-x)^{n-k+1}$$

$$+(-1)^k n(n-1)...(n-k+1)(1-x)^{n-k}$$

and taking the next derivative we obtain

$$f_n^{(k+1)}(x) = (-1)^{k+2} (k+1)n(n-1)...(n-k+1)(1-x)^{n-k+1}$$

$$+(-1)^{k+1} n(n-1)...(n-k)(1-x)^{n-k-1}$$

Therefore, since the formula is true for $k = 2$, we have proved it by mathematical induction for all $k \geq 2$.

Hence we have for all $k \geq 2$

$$f_n^{(k)}(1) = (-1)^{k+1} kn(n-1)...(n-k+2)\delta_{n,k-1}$$

$$+(-1)^k n(n-1)...(n-k+1)\delta_{n,k}$$

and $\qquad f_n^{(1)}(1) = \delta_{n,0} - n\delta_{n,1}$

$$\varsigma_a(-m) = \sum_{n=0}^{\infty} \frac{G_n(m)}{2^{n+1}} = \sum_{n=0}^{\infty} \frac{1}{2^{n+1}} \sum_{k=1}^{m} S(m,k) f_n^{(k)}(1)$$

Accordingly we obtain

(F.19) $$\sum_{n=0}^{\infty} \frac{1}{2^{n+1}} \sum_{k=1}^{m} S(m,k) f_n^{(k)}(1) = -\frac{B_{m+1}}{m+1}$$

I suspect that the left-hand side can be simplified, but I shall leave that exercise to others.

We have the Hasse identity (3.12)

$$\varsigma(s) = \frac{1}{s-1} \sum_{n=0}^{\infty} \frac{1}{n+1} \sum_{k=0}^{n} \binom{n}{k} \frac{(-1)^k}{(k+1)^{s-1}}$$

Guillera and Sondow [75aa] have recently shown that

(F.20) $$B_m(x) = \sum_{n=0}^{m} \frac{1}{n+1} \sum_{k=0}^{n} \binom{n}{k} (-1)^k (k+x)^m$$

and using (A.14a) we have $B_{2m}(1) = B_{2m}(0) = B_{2m}$. Hence we have (see also (A.23))

(F.21) $$B_{2m}(1) = B_{2m} = \sum_{n=0}^{2m} \frac{1}{n+1} \sum_{k=0}^{n} \binom{n}{k} (-1)^k (k+1)^{2m}$$

We may also obtain the trivial zeros using (A.14b)

(F.22) $$\varsigma(-2m) = -\frac{1}{2m+1} \sum_{n=0}^{\infty} \frac{1}{n+1} \sum_{k=0}^{n} \binom{n}{k} (-1)^k (k+1)^{2m+1} = B_{2m+1}(1) = 0$$

We also have the Hasse identity for the Hurwitz zeta function

$$\varsigma(s,a) = \frac{1}{s-1} \sum_{n=0}^{\infty} \frac{1}{n+1} \sum_{k=0}^{n} \binom{n}{k} \frac{(-1)^k}{(k+a)^{s-1}}$$

and accordingly we have

(F.22) $$\varsigma(-m+1,a) = -\frac{1}{m} \sum_{n=0}^{\infty} \frac{1}{n+1} \sum_{k=0}^{n} \binom{n}{k} (-1)^k (k+a)^m = -\frac{B_m(a)}{m}$$

Using the Euler-Maclaurin summation formula [14a], Hardy [76aa, p.333] showed that the Riemann zeta function could be expressed as follows

(F.23a) $$\varsigma(s) = \lim_{n\to\infty}\left[\sum_{k=1}^{n}\frac{1}{k^s} - \frac{n^{1-s}}{1-s} - \frac{1}{2}n^{-s}\right] \qquad \text{Re}(s) > -1$$

(F.23b) $$\varsigma(s) = \lim_{n\to\infty}\left[\sum_{k=1}^{n}\frac{1}{k^s} - \frac{n^{1-s}}{1-s} - \frac{1}{2}n^{-s} + \frac{1}{12}sn^{-s-1}\right] \qquad \text{Re}(s) > -3$$

and for $\text{Re}(s) > -5$

(F.23c) $$\varsigma(s) = \lim_{n\to\infty}\left[\sum_{k=1}^{n}\frac{1}{k^s} - \frac{n^{1-s}}{1-s} - \frac{1}{2}n^{-s} + \frac{1}{12}sn^{-s-1} - \frac{1}{720}s(s+1)(s+2)n^{-s-3}\right]$$

Letting $s = 0$ in (F.23a) we immediately obtain (F.2) $\varsigma(0) = -\frac{1}{2}$. Similarly, letting $s = -1$ in (F.23b) we see that $\varsigma(-1) = -\frac{1}{12}$ which we proved earlier in (F.2).

Differentiating the above identities results in for $\text{Re}(s) > -1$

(F.24a) $$\varsigma'(s) = \lim_{n\to\infty}\left[-\sum_{k=1}^{n}\frac{\log k}{k^s} + \frac{n^{1-s}(1-s)\log n - n^{1-s}}{(1-s)^2} + \frac{1}{2}n^{-s}\log n\right]$$

and with $s = 0$ we obtain

$$\varsigma'(0) = \lim_{n\to\infty}\left[-\sum_{k=1}^{n}\log k + \left(n + \frac{1}{2}\right)\log n - n\right]$$

Hence, using the Stirling approximation (A.11) we see that $\varsigma'(0) = -\frac{1}{2}\log(2\pi)$.

For $\text{Re}(s) > -3$ we have

(F.24b)
$$\varsigma'(s) = \lim_{n\to\infty}\left[-\sum_{k=1}^{n}\frac{\log k}{k^s} + \frac{n^{1-s}(1-s)\log n - n^{1-s}}{(1-s)^2} + \frac{1}{2}n^{-s}\log n - \frac{1}{12}sn^{-s-1}\log n + \frac{1}{12}n^{-s-1}\right]$$

and with $s = -1$ we get

(F.24c) $$\varsigma'(-1) = \lim_{n\to\infty}\left[-\sum_{k=1}^{n}k\log k + \left(\frac{n^2}{2} + \frac{n}{2} + \frac{1}{12}\right)\log n - \frac{n^2}{4} + \frac{1}{12}\right]$$

We now note from (4.4.226) in Volume III that the Glaisher-Kinkelin constant is defined by

(F.24d) $$\log A=\lim_{n\to\infty}\left[\sum_{k=1}^{n}k\log k-\left(\frac{n^2}{2}+\frac{n}{2}+\frac{1}{12}\right)\log n+\frac{n^2}{4}\right]$$

and hence we see that

(F.24e) $$\log A=\frac{1}{12}-\varsigma'(-1)$$

For $\mathrm{Re}(s)>-5$ we have

(F.24f)

$$\varsigma'(s)=\lim_{n\to\infty}\left[\begin{array}{l}-\sum_{k=1}^{n}\frac{\log k}{k^s}+\frac{n^{1-s}(1-s)\log n-n^{1-s}}{(1-s)^2}+\frac{1}{2}n^{-s}\log n-\frac{1}{12}sn^{-s-1}\log n+\frac{1}{12}n^{-s-1}\\ +\frac{1}{720}s(s+1)(s+2)n^{-s-3}\log n-\frac{1}{720}(3s^2+6s+2)n^{-s-3}\end{array}\right]$$

and with $s=-2$ we can easily show that

(F.24g) $$\log B=\lim_{n\to\infty}\left[\sum_{k=1}^{n}k^2\log k-\left(\frac{n^3}{3}+\frac{n^2}{2}+\frac{n}{6}\right)\log n+\frac{n^3}{9}-\frac{n}{12}\right]$$

where we used $\log B$ in (6.84) and

(F.24h) $$\log B=-\varsigma'(-2)=\frac{\varsigma(3)}{4\pi^2}\,.$$

Similarly, we find that

(F.24i) $$\log C=\lim_{n\to\infty}\left[\sum_{k=1}^{n}k^3\log k-\left(\frac{n^4}{4}+\frac{n^3}{2}+\frac{n^2}{4}-\frac{1}{120}\right)\log n+\frac{n^4}{16}-\frac{n^2}{12}\right]$$

and

(F.24j) $$\log C=-\varsigma'(-3)-\frac{11}{720}$$

With regard to (F.24a) we could determine $\varsigma''(0)$

$$\varsigma''(s)=\lim_{n\to\infty}\left[\sum_{k=1}^{n}\frac{\log^2 k}{k^s}+\frac{(1-s)^2[-n^{1-s}(1-s)\log^2 n]+2[n^{1-s}(1-s)\log n-n^{1-s}]}{(1-s)^4}-\frac{1}{2}n^{-s}\log^2 n\right]$$

(F.24k) $$\varsigma''(0)=\lim_{n\to\infty}\left[\sum_{k=1}^{n}\log^2 k-n\log^2 n+2n\log n-2n-\frac{1}{2}\log^2 n\right]$$

and compare the result with (4.3.234) in Volume II(b).

$$\varsigma''(0) = \gamma_1 + \frac{1}{2}\gamma^2 - \frac{1}{24}\pi^2 - \frac{1}{2}\log^2(2\pi)$$

The equation (F.24k) is contained in Ramanujan's Notebooks [21, Part I, p.203].

Finally, for reference we note the following approximations contained in [16]

$$\varsigma(2) = 1.644934066848... \qquad \varsigma(3) = 1.202056903159...$$

$$\varsigma(4) = 1.082323233711... \qquad \log 2 = 0.693147180559...$$

$$Li_4(1/2) = 0.517479061673...$$

**ACKNOWLEDGEMENTS**

I would like to wholeheartedly thank those mathematicians who make their papers freely available on the internet: without them, this paper would never have been written (or at least, not by me!). It is a real shame that more papers are not widely disseminated via arXiv, http://arxiv.org/archive/math , the CoLab Document Server, CiteSeer, http://citeseer.ist.psu.edu/cs , or similar repositories. Another extremely useful electronic aid is the Wolfram Integrator which enabled me to "play" with so many integrals without undergoing the drudgery of their evaluation.

More specifically, I would like to thank Mark Coffey and Jonathan Sondow for unhesitatingly providing me with copies of, and references to, various relevant papers.

In addition, I would also like to posthumously thank Paul Erdös for showing us that it is perfectly acceptable for a mathematician to write a paper using the first person singular.

Finally, I express my deep gratitude to my long-suffering wife, Sarah, and our three children who had to tolerate our dining room table being almost hidden by a mountain of papers for more than three years while I conducted research for this series of papers (and, according to my wife, just because there are three $n$'s in my surname $= n^3$). Unfortunately, the only zeta with which my children are familiar is the one usually associated with the names Catherine and Jones!

**ERRORS AND OMMISIONS**

During the course of my research I have noted various misprints and errors in a number of other papers; my work is not immune from such frailties and shall therefore be most grateful if readers will report all of my mistakes to me by email so that a more accurate version of this series of papers may be subsequently produced.

After graduating with a rather mediocre degree in mathematics from Queen's University Belfast in 1974, the author qualified as a chartered accountant and is currently employed as the chief financial officer of a quoted investment company in London. He specialises in tax and corporate finance and, in his spare time, he now dabbles in mathematics (after a very prolonged absence of nearly 30 years from that field of previously little endeavour).

A tax case, Smith v Schofield (249 STC 1992), was heard in the UK Court of Appeal in 1992 concerning the interpretation of some fairly obscure legislation regarding tax relief for the effects of inflation. Learned counsel put forward some relatively simple mathematical formulae in evidence before the court, only to be rebuked with the following comments from one of the judges [p.256]

"A great deal of the material - some 35 pages of it – has been included in the appeal bundle. Apart from allowing counsel to refer to the graphs in support of their arguments, we declined to consider it. Nor were we prepared to consider a detailed criticism of the Crown's calculations and methods, prepared by an accountant, which Mr Park sought to introduce by way of a riposte. Simple algebraic formulae, of the kind set out in para 11 of Sch 5, may be helpful in illustrating the meaning of the words used, but to my mind it would be absurd to seek for the intention of Parliament in page after page of abstruse mathematical calculations, all of them founded on arbitrary and controversial assumptions. If material of this sort is ever to be introduced in a tax appeal it should in my judgement be classified as expert evidence, and should be produced before the Special or General Commissioners where it can be tested by cross-examination of the author, rather than on the hearing of an appeal by way of case stated in the High Court."

The formula actually contained in the legislation referred to above was very simple. I did not have any involvement with this case but it immediately struck me that this was hardly British judiciary acting at its best!

Whilst on this occasion the UK court was reluctant to consider the mathematical intricacies, I am reminded that the position was rather different in the US in 1897 when the Indiana State Legislature [23aa, p.247] was champing at the bit to pass a bill stating that $\pi = 4$ (exactly)!

The author certainly finds mathematics to be more lucid than the morass of tax law which he had to decipher on a daily basis!

Donal F. Connon
Elmhurst
Dundle Road
Matfield
Kent TN12 7HD
dconnon@btopenworld.com